\documentclass[psamsfonts, reqno, 10pt]{amsart}
\usepackage{amsfonts,amsmath, amsthm, amssymb, latexsym, epsfig, verbatim}
\usepackage[all]{xy}
\usepackage{xspace}
\usepackage[mathscr]{eucal}

	\advance \topmargin by -\headheight
	\advance \topmargin by -\headsep

	\evensidemargin \oddsidemargin
\renewcommand{\include}{\input} % this will change all \includes to \inputs

\RequirePackage{amsfonts}
\RequirePackage{amsmath}
\RequirePackage{amsthm}
\RequirePackage{amssymb}
\RequirePackage{latexsym}
\RequirePackage{epsfig}
\RequirePackage{verbatim}
\RequirePackage[mathscr]{eucal}
\RequirePackage[all]{xy}

	\newcommand{\setof}[2]{\ensuremath{\left\{ #1 \: : \: #2 \right\}}}
	
	\newcommand{\set}[1]{\mathbb #1}
	
	\newcommand{\ftn}[3]{ #1 : #2 \rightarrow #3 }

	\newcommand{\innerprod}[2]{\left\langle #1,#2 \right\rangle} 
	
	\newcommand{\norm}[1]{\left\|#1\right\|}
	\newcommand{\indmap}[1]{\delta_{1}^{#1}}
	\newcommand{\expmap}[1]{\delta_{0}^{#1}}
	
	\newcommand{\mat}[4]{\left( \begin{array}{cc} #1&#2\\#3&#4\\ \end{array} \right)}

	\newcommand{\voline}[1]{ \overline{V} ( #1 ) }
	
	\newcommand{\vstar}[1]{ V_{*} ( #1 ) }
	\newcommand{\vtilde}[1]{( \vstar{ #1}, \totalkv{#1})}
	\newcommand{\vtildeh}[1]{[#1]}
	\newcommand{\hypon}[1]{k ( #1  )_{+}}
	\newcommand{\hypo}[1]{k ( #1  )}
	\newcommand{\unit}[1]{1_{ #1 }}
	\newcommand{\unitize}[1]{\widetilde{#1}}
	\newcommand{\matcont}[2]{M_{#1} \left( C(#2) \right)}
	\newcommand{\matalg}[2]{M_{#1} \left( #2 \right)}
	\newcommand{\md}[1]{\overline{#1}}
	\newcommand{\quot}[1]{Q(#1)}
	\newcommand{\ideal}[1]{I(#1)} 
	\newcommand{\sixk}{\ensuremath{\mathbf{K}}\xspace}
	\newcommand{\kmod}[3]{K_{#1}( #2 ; \Z / #3 \Z ) }
	\newcommand{\hclass}[1]{\left[ #1 \right]_{k}}
	\newcommand{\kclass}[1]{\left[ #1 \right]_{0}}
	\newcommand{\totalk}[1]{\underline{K} ( #1 )} 
	\newcommand{\totalkv}[1]{\widetilde{K} ( #1 )} 
	\newcommand{\mc}[1]{\mathcal{#1}}
	
	\newcommand{\Z}{\ensuremath{\mathbb{Z}}\xspace}
	\newcommand{\C}{\ensuremath{\mathbb{C}}\xspace}
	
	\newcommand{\R}{\ensuremath{\mathbb{R}}\xspace}
	
	\newcommand{\PN}{\ensuremath{\mathcal{P}}\xspace}
	\newcommand{\K}{\ensuremath{\mathcal{K}}\xspace}
	\newcommand{\dmap}{\ensuremath{\mathnormal{d}}\xspace}
	
	\newcommand{\mcalE}{\ensuremath{\mathcal{E}}}
	
	\newcommand{\submcalE}{\ensuremath{\mathcal{E}_{0}}}

	\newcommand{\mbbT}{\ensuremath{\mathbb{T}}}
	\newcommand{\hzero}{\ensuremath{\mathcal{S}}\xspace}
	\newcommand{\kthy}{\ensuremath{\mathnormal{K}}-theory\xspace}

		\newcommand{\kl}{\ensuremath{\operatorname{KL}}}
	\newcommand{\kk}{\ensuremath{\operatorname{KK}}}

	\newcommand{\spec}{\ensuremath{\operatorname{sp}}}
	\newcommand{\Hom}{\ensuremath{\operatorname{Hom}}}
	
	\newcommand{\Aut}{\ensuremath{\operatorname{Aut}}}
	
	\newcommand{\innerauto}{\ensuremath{\operatorname{Ad}}}
	\newcommand{\diag}{\ensuremath{\operatorname{diag}}}
	
	\newcommand{\id}{\ensuremath{\operatorname{id}}}
	\newcommand{\dist}{\ensuremath{\operatorname{dist}}}
	\newcommand{\bootstrap}{\ensuremath{\mathcal{N}}\xspace}
	\newcommand{\Extab}{\ensuremath{\operatorname{Ext}}_{\Z}^{1}}
	\newcommand{\pext}{\ensuremath{\operatorname{Pext}}_{\Z}^{1}}
	\newcommand{\Ext}{\ensuremath{\operatorname{Ext}}}
	\newcommand{\ext}{\ensuremath{\operatorname{ext}}_{\Z}^{1}}

	\newcommand{\purelyinf}{separable amenable purely infinite simple C*-algebra\xspace}
	\newcommand{\multialg}[1]{M(#1)\xspace}
	\newcommand{\corona}[1]{M(#1)/#1\xspace}

	\newcommand{\ran}{\ensuremath{\operatorname{ran}}\xspace}
	\newcommand{\dirlim}{\displaystyle \lim_{\longrightarrow}}

\newcommand{\subaealgs}{\ensuremath{\mathcal{A}\submcalE}-algebras\xspace}
\newcommand{\aealgs}{\ensuremath{\mathcal{A}\mcalE}-algebras\xspace}
	\newcommand{\aealg}{\ensuremath{\mathcal{A}\mcalE}-algebra\xspace}
	\newcommand{\ealgs}{\ensuremath{\mcalE}-algebras\xspace}
	\newcommand{\ealg}{\ensuremath{\mcalE}-algebra\xspace}

\newcommand{\subaealg}{\ensuremath{\mathcal{A}\submcalE}-algebra\xspace}
	\newcommand{\subealgs}{\ensuremath{\submcalE}-algebras\xspace}
	\newcommand{\subealg}{\ensuremath{\submcalE}-algebra\xspace}

	\def\ATalg{A\mbbT-algebra\xspace}
	\def\ATalgs{A\mbbT-algebras\xspace}

	\theoremstyle{plain}
	\newtheorem{thm}{Theorem}[section]
	\newtheorem{lemma}[thm]{Lemma}
	\newtheorem{theorem}[thm]{Theorem}
	\newtheorem{proposition}[thm]{Proposition}
	\newtheorem{corollary}[thm]{Corollary}

	\theoremstyle{definition}
	\newtheorem{definition}[thm]{Definition}
	
	\newtheorem{remark}[thm]{Remark}
	\newtheorem{notation}[thm]{Notation}

	\numberwithin{equation}{section}
	\numberwithin{figure}{section}

\begin{document}
	\title[Extensions of circle algebras by purely infinite C*-algebras]%
	{A classification theorem for direct limits of extensions of circle algebras by purely infinite C*-algebras}
	\author{Efren Ruiz}
        \address{Department of Mathematics \\    
        University of Toronto \\
        40 St. George St. \\
        Toronto, Ontario M5S 2E4 CANADA}
        \email{eruiz@math.utoronto.edu}
        \date{\today}
	
%AMS info

	\keywords{Classification, Amenable C*-algebras}
	\subjclass[2000]{Primary: 46L35}

	\begin{abstract}
	We give a classification theorem for a class of C*-algebras which are direct limits of finite direct sums of \subealgs.  The invariant consists of the following: (1) the set of Murray-von Neumann equivalence classes of projections;  (2) the set of homotopy classes of hyponormal partial isometries; (3) a map $d$; and (4) total $K$-theory.
	\end{abstract}
        \maketitle

%\begin{center}
%\textbf{Draft - Not for Distribution}
%\end{center}      
         
\section*{INTRODUCTION}
The program of classifying all amenable C*-algebras was initiated by George Elliott.  He's work involving AF-algebras has provided a foundation for other classification theorems.  Elliott \cite{af} proved that the class of all AF-algebras are classified by their dimension groups.  In 1993, Elliott \cite{at} classified another class of C*-algebras.  He showed the class of all \ATalgs are classified by their $K_{1}$ group and their ordered group $K_{0}$.  This remarkable result lead to many classification results for separable amenable C*-algebras.  See \cite{classpr} for references.  Most of these results involved simple C*-algebras with stable rank one or purely infinite simple C*-algebras.  

In 1997, Lin and Su \cite{gentoe} classified a class of (not necessarily simple) separable amenable C*-algebras with real rank zero that can be expressed as a direct limit of generalized Toeplitz algebras.  The invariant used by Lin and Su consists of the following three objects:  (1) $V(A)$, the set of Murray-von Neumann equivalence classes of projections; (2) $\hypon{A}$, certain equivalence classes of hyponormal partial isometries; and (3) a map $\ftn{d}{\hypon{A}}{V(A)}$.  Denote this invariant by $\vstar{A}$.  Lin and Su showed $V_{*}$ is a complete invariant for the class of all unital separable amenable C*-algebras with real rank zero that can be expressed as a direct limit of generalized Toeplitz algebras.

$V_{*}$ was then used by Lin \cite{inftoeplitz} to classify another class of (not necessarily simple) amenable separable C*-algebras.  Its basic building blocks consists of C*-algebras that are finite direct sums of corners of unital essential extensions of $\matcont{k}{S^{1}}$ by a non-unital Cuntz algebra $\mc{O}_{m} \otimes \K$.  Lin proved that all unital C*-algebras with real rank zero that can be expressed as a direct limit of these building blocks are classified by $V_{*}$.  Note by results of Cuntz (\cite{cuntzpureinf} and \cite{kthypureinf}), $\mc{O}_{m} \otimes \K$ is a simple amenable separable C*-algebra with $K_{0}(\mc{O}_{m} \otimes \K) = \Z / m \Z$ and $K_{1}(\mc{O}_{m} \otimes \K) = 0$.  

In this paper, we consider C*-algebras which can be expressed as direct limits of finite direct sums of corners of unital essential extensions of $\matcont{k}{S^{1}}$ by a \purelyinf $I$ in the bootstrap class \bootstrap.  A C*-algebra of this form will be called an \aealg.  We will classify a subclass of the class of all unital \aealgs with real rank zero.   A C*-algebra in this subclass will be called an \subaealg.  We will show that the class of all unital \subaealgs with real rank zero are classified by $V_{*}$ and total $K$-theory.  In addition, we will show that all unital \ATalgs with real rank zero and all the C*-algebras classified in \cite{inftoeplitz} are \subaealgs with real rank zero.  Also, by a result of Kirchberg \cite{kirchpureinf} (and independently by Phillips \cite{phillipspureinf}), the class of all unital separable amenable purely infinite simple C*-algebras with torsion free $K_{1}$ in the bootstrap class \bootstrap is contain in the class of all unital \subaealgs with real rank zero.  

The paper is organized as follows.  In section 1, we give some definitions and basic properties.  In section 2, we give perturbation lemmas that are used throughout the paper.  In section 3, we introduce the invariant.  In section 4 and section 5, we prove a uniqueness theorem and an existence theorem.  In section 6, we prove our main result.

\vspace{.7cm}

\noindent\textbf{Acknowledgements.}  The author would like to thank his Ph.D advisor, Huaxin Lin, for his guidance and support during the entire process of this paper.  Most of this work appeared in the author's Ph.D thesis.  The author would also like to thank N. Christopher Phillips for many useful conversations.     
 
\section{DEFINITIONS AND BASIC PROPERTIES}\label{definitions}
\subsection{PRELIMINARIES}
(1) Let $X$ be a compact subset of $S^{1}$.  If $z(t) = t$ for all $t \in X$, then $z$ is called the standard unitary generator for $C(X)$.\\
(2) For a unital C*-algebra $E$, set $U(E) = \setof{u \in E}{u \ \mbox{unitary}}$ and set \\
$U_{0}(E) = \setof{u \in U(E)}{ \exists f \in U(C([0,1], E)) \ \mbox{such that} \ f(0)= 1_{E} \ \mathrm{and} \ f(1) = u}$.\\
(3) $C$ is a unital C*-subalgebra of a unital C*-algebra $A$ if $C \subset A$ and $\unit{C} = \unit{A}$.   \\
(4) The set of all $*$-homomorphisms from $A$ to $B$ will be denoted by $\Hom(A,B)$.\\
(5) We write $p \sim q$ if $p$ is Murray-von Neumann equivalent to $q$.  \\
(6) An extension $0 \to B \overset{i}{\to} E \overset{\pi}{\to} A \to 0$ of C*-algebras is said to be \textbf{essential} if for every nonzero ideal, $J$ of $E$, we have $J \cap i(B) \neq 0$.  If $E$ is a unital essential extension of $A$ by $B$, then $E$ can be identified as a unital C*-subalgebra of $\multialg{B}$.  Also, recall that every essential extension can be described by its Busby invariant $\ftn{\tau_{E}}{A}{\corona{B}}$.\\
(7) Suppose $E$ is an extension of $A$ by $B$.  Then the exponential map and the index map associated to $E$ will be denoted by $\expmap{E}$ and $\indmap{E}$ respectively.\\
(8) Let \bootstrap be the bootstrap category of \cite{uct} and let \PN be the class of all $I \in \bootstrap$ such that $I$ is a non-unital \purelyinf with $K_{*}(I)$ finitely generated.  Let $\PN_{1}$ be the class of all unital separable amenable purely infinite simple C*-algebras in \bootstrap.\\
(9)  All semigroups have identities and all homomorphisms preserve identities.

\subsection{$\mc{E}$-ALGEBRAS}
\begin{definition} \label{def D : ealgs}
Suppose $E = p M_{n} (E') p$ for some $n \in \Z_{>0}$ and for some projection $p \in M_{n} (E') $. Suppose $E'$ is a unital essential extension of $C(S^{1})$ by a C*-algebra $I \in \PN$.  Let $\pi$ be the quotient map induced by the extension $E'$.  If $\pi(p) = \unit{ \matcont{n}{S^{1}} }$ or $p \in I$, then $E$ is said to be an \textbf{$\mc{E}$-algebra}.  If $I \in \PN$ and either $K_{1}(I) = 0$ or $K_{1}(I)$ is torsion free and $\ker \indmap{E'} \neq \{0\}$, then $E$ is called an \textbf{$\mc{E}_{0}$-algebra}. 
\end{definition}

Proposition \ref{def P : her ealg} will show that $p \matalg{n}{E'} p$ is an \ealg for any projection $p$.  

\begin{notation}
Let $E$ be an \ealg.  If $E$ is not a unital purely infinite simple C*-algebra, then $E$ is a unital essential extension of $M_{n} ( C(S^{1}) )$ by $I$, for some $I \in \PN, \ n \in \Z_{>0}$.  Denote $I$ by $\ideal{E}$ and denote $\matcont{n}{S^{1}}$ by $\quot{E}$.  If $E$ is a unital purely infinite simple C*-algebra, then set $\ideal{E} = E$ and set $\quot{E} = 0$.  If $E = \bigoplus_{i=1}^{k} E_{i}$, where each $E_{i}$ is an \ealg, then set $\ideal{E} = \bigoplus_{i=1}^{k} \ideal{E_{i}}$ and set $\quot{E} = \bigoplus_{i = 1}^{k} \quot{E_{i}}$.  It is clear that $0 \to \ideal{E} \to E \to \quot{E} \to 0$ is exact.  Hence $E \in \bootstrap$ since $\ideal{E}$ and $\quot{E}$ are in \bootstrap.  
\end{notation}
 
\begin{proposition}\label{def P : induce hom} 
Suppose $E_{1}$ and $E_{2}$ are finite direct sums of \ealgs.  Let $\ftn{i_{k}}{\ideal{E_{k}}}{E_{k}}$ be the inclusion map and let $\ftn{\pi_{k}}{E_{k}}{\quot{E_{k}}}$  be the quotient map for $k = 1,2$.  If $\varphi \in \Hom(E_{1}, E_{2})$, then there exist unique $\ideal{\varphi} \in \Hom(\ideal{E_{1}} , \ideal{E_{2}})$ and $\quot {\varphi} \in \Hom(\quot{E_{1}}, \quot{E_{2}})$ such that $i_{2} \circ \ideal{\varphi} = \varphi \circ i_{1}$ and $\quot{\varphi} \circ \pi_{1} = \pi_{2} \circ \varphi$.
\end{proposition}

\begin{proof}
Since $\ideal{E_{1}} $ is a finite direct sum of purely infinite simple C*-algebras and $\quot{E_{2}}$ is a finite C*-algebra, we have $\pi_{2} \circ \varphi \circ i_{1} = 0$.  Therefore, $\ran (\varphi \circ i_{1}) \subset \ran(i_{2})$.  Hence, there exists $\ideal{\varphi} \in \Hom(\ideal{E_{1}} , \ideal{E_{2}})$ such that $i_{2} \circ \ideal{\varphi} = \varphi \circ i_{1}$.  The rest of the proposition now follows.
\end{proof}

\begin{remark}
Suppose $E_{1}$ is an \ealg and $\ker \varphi \neq \{ 0 \}$.  Then $\ran (\varphi \circ i_{1}) = \{ 0 \} $.  Hence, there exists $ \varphi_{Q_{1}} \in \Hom( \quot{E_{1}} ,  E_{2} )$ such that $\varphi_{Q_{1}} \circ \pi_{1} = \varphi$ and $\pi_{2} \circ \varphi_{Q_{1}} = \quot{\varphi}$.
\end{remark}
  
\begin{definition}\label{def D : aealgs}
If $E$ is a direct limit of finite direct sums of \ealgs, then $E$ is called an \textbf{$\mc{A}\mc{E}$-algebra}.  If $E$ is a direct limit of finite direct sums of \subealgs, then $E$ is called an \textbf{$\mc{A}\mc{E}_{0}$-algebra}.
\end{definition}

Let $E = \dirlim ( E_{n}, \varphi_{n, n+1})$ be an \aealg.   Set $\ideal{E} = \dirlim ( \ideal{E_{n}} , \ideal{\varphi_{n, n+1}} )$ and set $\quot{E} = \dirlim( \quot{E_{n}}, \quot{\varphi_{n, n+1}) }$.  Then $0 \rightarrow \ideal{E} \rightarrow E \rightarrow \quot{E} \rightarrow 0$ is exact.

\begin{proposition}\label{def P : induce hom aealg}
Proposition \ref{def P : induce hom} is still true when $E_{1}$ and $E_{2}$ are \aealgs. 
\end{proposition}

\begin{proof}
Note that $\quot{E_{1}}$ and $\quot{E_{2}}$ are finite C*-algebras and $\ideal{E_{1}}, \ideal{E_{2}}$ are direct limits of finite direct sums of purely infinite simple C*-algebras.  Hence, arguing as in Proposition \ref{def P : induce hom}, we get the desired conclusion.
\end{proof}

\begin{proposition}\label{def P : liftingproj}
Suppose $E$ is an \ealg such that $\quot{E} \cong \matcont{k}{S^{1}}$.  If $e$ is a projection in $E$ not in $\ideal{E}$, then every projection $p \in \matalg{n}{\quot{e E e}}$ lifts to a projection in $\matalg{n}{e E e}$ for all $n \in \Z_{>0}$.  Consequently, $\expmap{e E e} = 0$.  
\end{proposition}

\begin{proof}
By \cite[1.2]{corona1}, $RR(\ideal{E}) = 0$.  The proposition now follows from Lemma 2.5, Lemma 2.8, and Remark 2.9 in \cite{quasidiag}.
\end{proof}

\begin{proposition}\label{def P : lifting matrix unit}
Let $0 \to B \to E \to A \to 0$ be a unital extension.  Suppose

(1)  every hereditary C*-subalgebra of $B$ has an approximate identity consisting of projections and

(2)  every projection in $A$ lifts to a projection in $E$.

\noindent Supppose $C$ is a unital C*-subalgebra of $A$ such that $C \cong \matalg{n}{\C}$.   Let $\{\overline{ e}_{ij} \}_{i, j = 1}^{n}$ be a system of matrix units for $C$.  Then, there exists a system of matrix units $\{e_{ij} \}_{i, j = 1}^{n}$ in $E$ such that $\unit{E} - \sum_{i = 1}^{n} e_{ii} \in B$.
\end{proposition}

Effros \cite[9.8]{dim} proved the above proposition when $B$ is an AF-algebra.  The key components of the proof are (1) and (2).  

\begin{corollary}\label{def C : liftmatunit}
Suppose $E$ is an \ealg such that $\quot{E} \cong \matcont{l}{S^{1}}$.  Let $e \in E$ be a projection not in $\ideal{E}$.  Suppose $\{q_{ij}\}_{i,j = 1}^{n}$ is a system of matrix units for $\matalg{n}{\C} \subset \quot{e E e} \cong \matcont{n}{S^{1}}$.  Then there exists a system of matrix units $\{e_{ij}\}_{i,j= 1}^{n}$ in $e E e$ such that $\pi(e_{ij}) = q_{ij}$ and $e - \sum_{i=1}^{n} e_{ii} \in \ideal{e E e}$.
\end{corollary}

\begin{proof}
Since $\ideal{e E e} \in \PN$, by \cite[1.2]{corona1} and \cite[2.6]{heralgs}, $\ideal{e E e}$ satisfies (1) in Proposition \ref{def P : lifting matrix unit}.  By Proposition \ref{def P : liftingproj}, $e E e$ satisfies (2) in Proposition \ref{def P : lifting matrix unit}.  
\end{proof}

\begin{proposition}\label{def P : absorbing extensions}
Let $A$ be a unital amenable C*-algebra in \bootstrap and let $I \in \PN$.  Let $\ftn{ \tau }{ A }{ \corona{I} }$ be a unital monomorphism.  Then, for any strongly unital trivial extension $\tau_{0}$, there exists $u \in U(\multialg{I})$ such that $\innerauto(\pi(u))\circ \tau = \tau \oplus \tau_{0}$.
\end{proposition}

\begin{proof}
Since $\ftn{ \tau_{0} }{ A }{ \corona{ I } }$ is a strongly unital trivial extension, there exists a unital homomorphism $\ftn{ \sigma }{ A }{ \multialg{ I } }$ such that $\pi \circ \sigma = \tau_{0}$.  Choose isometries $s_{1}, s_{2} \in \multialg{ I }$ such that $1 = s_{1} s_{1}^{*} + s_{2} s_{2}^{*}$.  Then $\tau \oplus \tau_{0}$ can be represented by $\pi (s_{1} ) \tau (a) \pi (s_{1}^{*}) + \pi (s_{2}) \tau_{0}(a) \pi (s_{2}^{*})$.  Let $E_{1} = \pi^{-1} ( \tau ( A ) )$.  Define $\ftn{ \eta }{ E_{1} }{\multialg{I} }$ by $\eta (x) = s_{1} x s_{1}^{*} + s_{2} (\sigma \circ \pi (x)) s_{2}^{*}$ for all $x \in E_{1}$.  Since $s_{2} (\sigma \circ \pi(\cdot) )s_{2}^{*}|_{I} = 0$, by \cite[8.3.1, pp. 125]{classentropy} (\cite[7(iii)]{kirchpureinf}), there exists $u \in U(\multialg{ I })$ such that $\eta(x) - u^{*} x u \in I$ for all $x \in E_{1}$.  Hence $\pi (u^{*} ) \tau (a) \pi (u) = ( \tau \oplus \tau_{0} )(a)$.
\end{proof}

\begin{corollary}\label{def C : isoealg}
Let $E$ and $E'$ be two \ealgs such that $Q = \quot{E} \cong \quot{E'}$ and $I = \ideal{E} \cong \ideal{E'}$.  If $[ \tau_{E} ] = [\tau_{E'}]$ in $\Ext(Q, I)$, then $E \cong E'$.
\end{corollary}

The next two lemmas are well-known and we state it without proof.

\begin{lemma}\label{def L : liftunitiso}
Let $I \in \PN$.  Let $E$ be a unital C*-subalgebra of $\multialg{I}$ which contains $I$ as an essential ideal.  Then every unitary in $E/I$ lifts to a non-unitary isometry in $E$.
\end{lemma} 

\begin{lemma}\label{def L : trivial ext indmap zero}
Let $I \in \PN$.  Suppose $E$ is a unital C*-algebra containing $I$ as an essential ideal.  Then, $u \in U(E/I)$ can be lifted to $v \in U(E)$ if and only if $\indmap{E}([ u ] ) = 0$.  Moreover, if $E$ is an \ealg and $\quot{E} = C(S^{1})$, then $E$ has a splitting if and only if $ [\tau_{E} ] = 0$ in $\Ext(C(S^{1}), \ideal{E})$.  
\end{lemma}

\begin{definition}\label{def D : quasidiag}
Suppose $I \in \PN$.  A unitary $U \in \multialg{I}$ is \textbf{quasi-diagonal} if for every $\epsilon > 0$, there exist a sequence $\{ e_{k} \}_{ k = 1}^{\infty}$ of mutually orthogonal projections and a dense sequence $\{ \lambda_{k} \}_{k = 1}^{\infty}$ in $\spec(U)$ such that 

(1) $\left\{ f_{n} = \sum_{k = 1}^{n} e_{k} \right\}_{n = 1}^{\infty}$ is an approximate identity of $I$ consisting of projections,

(2) $U - \sum_{k = 1}^{\infty} \lambda_{k} e_{k} \in I$, and

(3) $\norm{ U - \sum_{k = 1}^{\infty} \lambda_{k} e_{k} } < \epsilon$,

\noindent where the sums converge in the strict topology.
\end{definition}

\begin{proposition}\label{def P : diagonalization}
Suppose $I \in \PN$.  Suppose $E$ is a unital C*-subalgebra of $\multialg{I}$ that contains $I$ as an essential ideal.  Let $\ftn{\pi}{E}{E/I}$ be the quotient map.  If $U \in U(E)$ such that $\spec(U) = \spec(\pi(U))$, then $U$ is quasi-diagonal.
Consequently, if $E$ is an \ealg such that $E$ is a trivial extension, then $E$ is a quasi-diagonal extension and every $v \in U(E)$ is quasi-diagonal.
\end{proposition}

\begin{proof}
Choose a sequence $\{ \lambda_{k} \}_{k = 1}^{\infty}$ in $\spec(U)$ such that for all $n \in \Z_{>0}$, the sequence $\{ \lambda_{k} \}_{k = n}^{\infty}$ is dense in $\spec(U)$.  Let $\{ p _{k} \}_{k = 1}^{\infty}$ be a sequence of mutually orthogonal projections in $I$ such that $\left\{ q_{n} = \sum_{k = 1}^{n} p_{k} \right\}_{n = 1}^{\infty}$ is an approximate identity for $I$.  Set $U' = \sum_{k = 1}^{\infty} \lambda_{k} p_{k}$, where the series converges in the strict topology.  Then $\spec( \pi(U)  ) = \spec( \pi(U') ) = X \subset S^{1}$.  Define $\ftn{ \sigma_{1}, \sigma_{2} } { C(X) } { \multialg{I} }$ by $\sigma_{1}(z) = U$ and $\sigma_{2} (z) = U'$.  Let $E_{1} = C$*$( U , I)$ and let $E_{2} = C$*$( U' , I)$.  Let $s_{1}, s_{2} \in \multialg{I}$ be isometries in $\multialg{I}$ such that $1 = s_{1} s_{1}^{*} + s_{2} s_{2}^{*}$.  Define $\ftn{ \eta_{1}}{ E_{1} }{ \multialg{ I } }$ by $\eta_{1} (x) = s_{1} x s_{1}^{*} + s_{2} (\sigma_{2} \circ \pi (x)) s_{2}^{*}$ for all $x \in E_{1}$ and define $\ftn{\eta_{1}}{E_{2}}{\multialg{I}}$ by $\eta_{2} (a) = s_{1} (\sigma_{1} \circ \pi(a)) s_{1}^{*} + s_{2} a s_{2}^{*}$ for all $a \in E_{2}$.   

Note that $\eta_{1} (U) = \eta_{2} (U')$.  By \cite[8.3.1, pp. 125]{classentropy} (\cite[7(iii)]{kirchpureinf}), there exist unitaries $V, W \in \multialg{I}$ such that $\norm{ \eta_{1} (U) - V U V^{*} } < \frac{\epsilon}{2}$, $\norm{ \eta_{2} (U') - W U' W^{* } } < \frac{\epsilon}{2}$, $\eta_{1} (U) - V U V^{* } \in I$, and $\eta_{2} (U') - W U' W^{* } \in I$.  Therefore, $U - V^{*} W U' W^{*} V \in I$ and $\norm{ U - V^{*} W U' W^{*} V } < \epsilon$.  For all $k \in \Z_{>0}$, set $e_{k} = V^{*} W p_{k} W^{*} V$.  

The last statement follows from Lemma \ref{def L : trivial ext indmap zero} and the above result.
\end{proof} 

\begin{lemma}\label{def L : unitgen}
Suppose $E$ is an \ealg such that $\quot{E} = C(S^{1})$.  Let $\ftn{\pi}{E}{\quot{E}}$ be the quotient map.  Then $K_{1}(E) \cong K_{1}(\ideal{E}) \oplus \ran (\pi_{*,1})$, where $\ran (\pi_{*,1})$ is generated by $[ \pi(u) ] $ for some $u \in U(E)$.  
\end{lemma}

\begin{proof}
The lemma follows from the exactness of the six-term exact sequence in \kthy, $\expmap{E} = 0$ (Proposition \ref{def P : liftingproj}), and $K_{1}(C(S^{1}))= \langle [z] \rangle$.
\end{proof}

\begin{proposition}\label{def P : her ealg}
Let $E$ be an \ealg and let $p$ be a nonzero projection in $E$.  Then $p E p$ is an \ealg.  Moreover, if $E$ is an \subealg, then $p E p$ is an \subealg.
\end{proposition}

\begin{proof}
If $p \in \ideal{E}$, then $p E p$ is a unital separable purely infinite simple C*-algebra in \bootstrap such that $K_{*}(p E p) \cong K_{*}(p \ideal{E} p)$.  Hence, $p E p$ is an \ealg.  

Suppose $ p \notin \ideal{E}$.  Let $\{ \overline{e}_{ij} \}_{i,j=1}^{n}$ be system of matrix units of $\matcont{n}{S^{1}} \cong \quot{E}$.  By Corollary \ref{def C : liftmatunit}, there exists a system of matrix units $\{ e_{i j} \}_{i, j = 1}^{n} \subset p E p$ such that $\pi ( e_{i j} ) = \overline{e}_{i j}$ and $p - \sum_{ i=1 }^{n} e_{i i} \in p \ideal{E} p$.  Note that $e_{11 } p \ideal{E} p e_{11} \in \PN$.  Hence, $e_{11} p E p e_{11}$ is an \ealg with $\quot{e_{11} p E p e_{11}} \cong C(S^{1})$.  Then $ p E p$ is isomorphic to $q M_{n} (E') q$ and $q M_{n} (E') q$ is an \ealg, where $I' = e_{11} p \ideal{E} p e_{11}$, $E' = e_{11} p E p e_{11}$, and $q = \sum_{i = 1}^{n} e_{i i}$.  Note that $\ideal{E'} \cong \ideal{p E p}$. 

The last statement is clear from the above arguments.
\end{proof}

\section{PERTURBATION LEMMAS}\label{perturbation lemmas}

We now give several lemmas that will play an important role in many of the proofs in the sequel.  Let $A$ be a unital C*-algebra.  Denote the unit of $\matalg{m}{A}$ by $1_{m}$.  If $m=0$, then $1_{0} = 0$.  If $A$ is a non-unital C*-algebra, then $\unitize{A}$ will denote the C*-algebra obtained by adjoining a unit to $A$. 

\begin{lemma}\label{pert L : approx unitarily eq purely inf}
Let $A \in \mc{P}_{1}$ and let $\epsilon > 0$.  Suppose $u_{1}, u_{2} \in U(A)$ such that $[ u_{1} ] = [ u_{2} ]$ in $K_{1} (A)$ and $\sup \{\dist (\lambda , \spec( u_{i} ) )\ :\ \lambda \in S^{1} \} < \epsilon$ for $i=1,2$.  Then, there exists $W \in U(A)$ such that $\norm{ W^{*} u_{1} W- u_{2} } < 3 \epsilon$.
\end{lemma}

\begin{proof}
If $[u_{1}] = [ u_{2}] = 0$ in $K_{1}(A)$, then the conclusion follows from \cite[2.2]{inftoeplitz}.  If $[u_{1}] \neq 0$ in $K_{1}(A)$, then the conclusion follows from \cite[3]{normalelements}.
\end{proof}

If $E$ is a Hilbert $A$-module, let $\mc{L}_{A}(E)$ be the set of all module homomorphism $\ftn{T}{E}{E}$ for which $T$ has an adjoint $T^{*}$.   It is a well-known fact that $\mc{L}_{A}(E)$ is a C*-algebra with respect to the operator norm.  Let $x, y \in E$.  Define $\theta_{x,y} \in \mc{L}_{A}(E)$ by $\theta_{x,y} (z) = x \innerprod{y}{z}$.  Let $\mc{K}_{A}(E)$ be the C*-subalgebra of $\mc{L}_{A}(E)$ generated by the collection $\{ \theta_{x,y} \}_{x, y \in E}$.  Then $\mc{K}_{A}(E)$ is an ideal of $\mc{L}_{A}(E)$.    

\begin{definition}\label{pert D : hilbert}
Let $A^{(n)} = \bigoplus_{k = 1}^{n} A$ be the Hilbert $A$-module of orthogonal direct sum of $n$ copies of $A$.  $H_{A}$ will denote the following Hilbert $A$-module:
\begin{equation*}
\setof{ \{ a_{n} \} } { a_{n} \in A \ \mathrm{and} \ \left\{ \sum_{k = 1}^{n} a_{k}^{*} a_{k} \right\}_{n = 1}^{\infty} \ \mathrm{converges} \ \mathrm{in} \ \mathrm{norm} \  \mathrm{as} \ n \to \infty}.
\end{equation*}
\end{definition}

\begin{proposition}\label{pert P : hilbert}
$\mc{L}_{A}(A^{(n)}) \cong M_{n}(A) \cong \mc{K}_{A}(A^{(n)})$, $\mc{K}_{A}(H_{A}) \cong A \otimes \K$, and $\mc{L}_{A}(H_{A}) \cong \multialg{A \otimes \K}$ for any C*-algebra $A$.  
\end{proposition}

For a proof of the above proposition see \cite[15.2.12]{olsen}.  We will use the above identifications throughout.  Let 
{\scriptsize \begin{equation*}
 s_{n} = \left(
\begin{matrix}
0      		& 0      	& \dots  	& 0      	&  1 \\
1  		& 0      	& \dots  	& 0      	& 0 \\
0      		& 1 		& \dots  	& 0      	& 0 \\
\vdots 	& \vdots 	& \ddots 	& \vdots 	& \vdots \\
0      		& 0      	& \dots  	& 1  	& 0     
\end{matrix}
\right) \in \mc{L}_{A}(A^{(n)}).
\end{equation*}} 

\noindent Suppose $A$ is a unital C*-algebra.  Let $e_{i} = (0, \dots , 0, \unit{A}, 0 \dots, 0)$ be the element in $A^{(n)}$ that has $\unit{A}$ in the $i^{\emph{th}}$ coordinate.  Then $s_{n} (e_{i}) = e_{i+1}$ for all $i = 1, 2, \dots n-1$ and $s_{n} (e_{n}) = e_{1}$.  Note that $s_{n}$ is a unitary in $\mc{L}_{A}(A^{(n)}) \cong M_{n}(A)$.  Define the \textbf{standard unilateral shift} to be the element $S \in \mc{L}_{A}(H_{A})$ which sends $\{ a_{n} \}_{n = 1}^{\infty} \in H_{A}$ to $\{ 0, a_{1} , a_{2}, \dots \}$.  

\begin{lemma}\label{pert L : spectrum of a unitary}
Suppose that $A$ is a unital C*-algebra.  For any $u \in U(M_{m} (A))$ and for any integer $n > m$, set 
$w_{n} = \mat{ u } { 0 } { 0 }{1_{n-m}} s_{n}$.  Then 
\begin{equation*}
\sup\setof{ \dist(\lambda, \spec( w_{n} ) )}{ \lambda \in S^{1} } < \frac{ \pi } { n } + 2 \left( \frac{m}{n} \right)^{1/2}.
\end{equation*} 
\end{lemma}
  
The above lemma was proved by Lin \cite[2.4]{inftoeplitz} for the case $K_{1}(A) = 0$.  By the proof, we see that the assumption $K_{1}(A) = 0$ may be omitted.
  
\begin{lemma}\label{pert L : approx unitarily eq unitaries1}
Let $A \in \mc{P}_{1}$.  Then, for any $\epsilon > 0$ and for any $k \in \Z_{>0}$, there exists $N > k$ such that for all $n \geq N$ and for any $u, v \in U(\matalg{k}{A})$ with $[ u ] = [ v ] $ in $K_{1} (A)$, there exists $W \in U(\matalg{n}{A})$ such that 
\begin{equation*}
\norm{ W \diag( u , 1_{n-k} )s_{n} W^{*} - \diag (v , 1_{n-k} )s_{n} } < \epsilon.
\end{equation*}
\end{lemma}

\begin{proof}
The lemma follows from Lemma \ref{pert L : spectrum of a unitary} and Lemma \ref{pert L : approx unitarily eq purely inf}.
\end{proof}

\begin{lemma}(R\o rdam)\label{pert L : approx unitarily eq unitaries2}  
Let $A$ be a unital C*-algebra.  Let $l, k \in \Z_{>0}$ and let $m > k + l$.  Suppose $u \in U(\matalg{k}{A})$, $v\in U(A)$, and $w_{1} \in U(\matalg{m}{A})$.  
Suppose that 
\begin{equation*}
\norm{ w_{1} \diag (1_{l},  v,  1_{m - l - 1}) s_{m} w_{1}^{*} - \diag (1_{l},  u,  1_{m - k - l } ) s_{m} } \leq \frac{7 - 2 \pi}{l}.
\end{equation*}
Then, there exists $w \in U(\matalg{2m}{A})$ such that $w p_{2m} = p_{2m} w = p_{2m}$ and 
\begin{equation*}
\norm{ w \diag (1_{l}, v, 1_{2m - l - 1} ) s_{2m} w^{*} - \diag (1_{l} , u, 1_{2m - k - l} ) s_{2m} } \leq \frac{7}{l},
\end{equation*}
where $p_{2m} = \diag (0, 0 , \dots , 0 , 1) \in \matalg{2m}{A}$.
\end{lemma}

\begin{proof}
Let $n = 2 m$.  Define $X \in \matalg{n}{\C} \subset \matalg{n}{A} \cong \mc{L}_{A}(A^{(n)})$ as follows:
$X e_{i} = e_{2 i - 1}$ and $X e_{i + m} = e_{2 i}$ for $1 \leq i \leq m$,
where $\{ e_{i} \}_{i = 1}^{n}$ is the standard orthonormal basis for $\C^{n}$.  A direct computation shows that 
{\tiny
\begin{equation*}
X s_{n} X^{*} = 
\left(
\begin{matrix}
0      		& 0      	& \dots  	& 0      	&  v_{1} \\
1_{2}  	& 0      	& \dots  	& 0      	& 0 \\
0      		& 1_{2}  	& \dots  	& 0      	& 0 \\
\vdots 	& \vdots 	& \ddots 	& \vdots 	& \vdots \\
0      		& 0      	& \dots  	& 1_{2}  	& 0     
\end{matrix}
\right) \quad \mathrm{and} \quad
X 	\left(	\begin{matrix}
			s_{m} & 0 \\
			0 & s_{m} 
		\end{matrix}
	\right)  
X^{*} = 
\left(	\begin{matrix}
		0      		& 0      	& \dots  	& 0      	&  1_{2} \\
		1_{2}  	& 0      	& \dots  	& 0      	& 0 \\
		0      		& 1_{2}  	& \dots  	& 0      	& 0 \\
		\vdots 	& \vdots 	& \ddots 	& \vdots 	& \vdots \\
		0      		& 0      	& \dots  	& 1_{2}  	& 0     
	\end{matrix}
\right),
\end{equation*}
}

\noindent where 
$v_{1} = 	\left( 	\begin{matrix}
				0 & 1 \\
				1 & 0
			\end{matrix}
		\right)$.  Note that the eigenvalues of $v_{1}$ are $1$ and $-1$.  Hence, there exists a unitary $v_{2}$ such that $v_{2}^{l} = v_{1}$ and $\norm{ v_{2} - 1_{2} } \leq \frac{\pi}{l}$.  
		
Set $z = \diag (v_{1}, v_{2}^{l-1}, v_{2}^{l - 2}, \dots , v_{2}, 1_{2} , \dots, 1_{2} ) \in \matalg{n}{\C}$.  Then 
{\tiny
\begin{equation*}
z^{*} X s_{n} X^{*} z = 
\left(		\begin{matrix}
			0     		& 0     	&        	&   		&	        & 0      	& 1_{2}  \\
			v_{2} 	& 0     	&        	&   		&        	& 0      	& 0      \\
      					& v_{2} 	&        	&   		&        	& 0      	& 0      \\
      					&       	& \ddots 	&   		&        	& \vdots 	& \vdots \\
      					&       	&        	& 1_{2} 	&        	& 0      	& 0      \\
      					&       	&        	&   		& \ddots 	& \vdots 	& \vdots \\
			0     		&       	&        	&   		&        	& 1_{2}     & 0
		\end{matrix} \right)
\quad
\mathrm{and}\quad \norm{ X^{*} z^{*} X s_{n} X^{*} z X - 
\left(		\begin{matrix}
			s_{m} & 0 \\
			0 & s_{m} 
		\end{matrix}
\right) } \leq \frac{\pi}{ l}.
\end{equation*}
}

\noindent Let $z = \{ a( i , j ) \}_{i,j = 1}^{n}$ and let $X^{*} z X = \{ b( i , j ) \}_{i,j = 1}^{n}$.  Then 
\begin{equation*}
b(l + i, l + j ) = a(2 l + 2 i - 1, 2 l + 2 j - 1) \ \mathrm{if} \ i \neq j,\\
\end{equation*}
\begin{equation*}
b(i, m + l + j ) = a( 2 i - 1 , 2( l + j) ) = 0 \ \mathrm{if} \ 0 < i \leq l \ \mathrm{and} \ j > 0, \ \mathrm{and}\\
\end{equation*}
\begin{equation*}
b( m + l + i , j ) = a( 2( l + i ) , 2 j - 1 ) = 0 \ \mathrm{if} \ i > 0\ \mathrm{and} \ 0 < j \leq l.
\end{equation*}
Therefore 

{\scriptsize \begin{equation*}
X^{*} z X = 
\left(		\begin{matrix}
			B 	& 0       			& C     	& 0       \\
			0 	& 1_{m - l} 		& 0     	& 0       \\
			D 	& 0       			& 1_{l} 	& 0       \\
			0 	& 0       			& 0     	& 1_{m - l}
		\end{matrix} 
\right)\ \mathrm{and} \ \ 
X^{*} z^{*} X = 
\left(		\begin{matrix}
			B^{*} 	& 0       		& D^{*}     & 0       \\
			0 		& 1_{m - l} 	& 0     	& 0       \\
			C^{*} 	& 0       		& 1_{l} 	& 0       \\
			0 		& 0       		& 0     	& 1_{m - l}
		\end{matrix} 
\right),
\end{equation*}}

\noindent where $B, C,$ and $D$ are in $\matalg{l}{\C}$.  Hence, $\diag( 1_{l}, u, 1_{n - k - l} )$ and $ \diag( 1_{l}, v, 1_{n - l - 1})$ commute with $X^{*} z X$ and $X^{*} z^{*} X$.  Let $w = X^{*} z X \diag (w_{1} , 1_{m}) X^{*} z^{*} X$.  Then, we have
{\tiny
\begin{align*}
&\norm{ w \diag(1_{l}, v, 1_{n-l-1} )s_{n} w^{*} - \diag(1_{l}, u, 1_{n-k-l}) s_{n}} \leq \frac{\pi}{l} \\ 	
&\quad\quad+\norm{ X^{*} z X \diag(w_{1}, 1_{m}) \diag(1_{l}, v, 1_{n-l-1})\diag(s_{m}, s_{m})\diag(w_{1}^{*},1_{m})
												 X^{*} z^{*} X - \diag(1_{l}, u, 1_{n-k-l})s_{n}  }
						\\
&\hspace{.5in}\leq \frac{\pi}{l} + \frac{7 - 2 \pi }{l} 
		+ \norm{ X^{*} z X \diag(1_{l}, u, 1_{n-k-l}) \diag(s_{m}, s_{m}) X^{*} z X - \diag(1_{l}, u, 1_{n-k-l})s_{n}  }
							\\
 &\hspace{1.0in}\leq \frac{\pi}{l} + \frac{7 - 2 \pi}{l} + \frac{\pi}{l} = \frac{7}{l}.
\end{align*}
}

\noindent Finally, $w p_{n} = p_{n}$ since $X^{*} z X e_{n} = X^{*} z e_{n} = X^{*} e_{n} = e_{n}$.
\end{proof}

\begin{lemma}\label{pert L :  approx unitarily eq isometries}
Let $B \in \mc{P}_{1}$ and set $A = B \otimes \K$.  Let $E$ be the unital C*-subalgebra of $\multialg{A}$ generated by $S$ and $A$, where $S$ is the standard unilateral shift in $\multialg{A} \cong \mc{L}_{B}( H_{B} )$.  Let $\ftn{\pi}{E}{E/A}$ be the quotient map.  Suppose $T$ is an isometry in $E$ such that $\pi (T) = \pi (S)$ and $[ 1 - T T^{*} ] = [ 1 - S S^{*} ]$ in $K_{0} (A)$.  Then, for any $\epsilon  > 0$, there exists $w \in U(\unitize{A})$ such that $\norm{ w^{*} T w - S} < \epsilon$. 
\end{lemma}

\begin{proof}
Let $p_{1} = 1-SS^{*} = \unit{B}$ and $q_{1} = 1 - TT^{*}$.  Since $[ p_{1} ] = [ q_{1} ]$ in $K_{0}(A)$, by \cite[1.1]{corona2}, there exists a unitary $w_{00} \in \unitize{A} \subset E$ such that $w_{00}^{*} p_{1} w_{00} = q_{1}$.
By replacing $T$ with $w_{00} T w_{00}^{*}$, we may assume $p_{1} = q_{1}$.  Let $u =  p_{1} + T S^{*}$.  Then $u \in \unitize{A}$ and $u S = T$.  

Let $p_{i} = p_{1} \otimes e_{i i}$ and let $E_{k} = \sum_{i = 1}^{k} p_{i}$, where $\{ e_{i j} \}_{i, j = 1}^{\infty}$ is the standard system of matrix units for $\K$.  Choose $m' \in \Z_{>0}$ such that $\frac{7}{m'} < \epsilon$.  Since $\{ E_{n} \}_{n = 1}^{\infty}$ forms an approximate identity for $A$ consisting of projections, there exist $v \in U(\unitize{A})$ and a positive integer $L > m'$ such that (1) $p_{i} v = v p_{i} \ \mathrm{for} \ i = 2, \dots , k $; (2) $p_{1} v = v p_{1} = p_{1}$; (3) $v = v E_{L} + 1 - E_{L}$; and (4) $\norm{ v - u } < \frac{\epsilon}{2}$.  Suppose there exists $w \in U(\unitize{A})$ such that $\norm{ w^{*} S w - v S} < \frac{\epsilon}{2}$.  Then, $\norm{w^{*} S w - T} \leq \norm{ w^{*} S w - v S} + \norm{ v S - T } < \epsilon$.  Hence, by replacing $T$ by $vS$, we may assume $u=v$.  

Let $\{ e_{i} \}_{i = 1}^{\infty} = \{ (1, 0, \dots ) , (0 , 1, 0, \dots ), \dots \}$ be the standard orthonormal basis for $H_{B}$.  Define a scalar matrix $w_{0} \in \mc{L}_{B}(H_{B})$ as follows:
\begin{equation*}
\begin{array}{ccl}
w_{0} ( e_{i} ) 			= e_{ L + i + 1 }  	&\mathrm{ for } \quad 	&i = 1, 2, \dots , L + 1 \\
w_{0} ( e_{L + i + 1} ) 	= e_{i} 			&\mathrm{for} \quad 		&i = 1, 2, \dots, L + 1\\
w_{0} ( e_{ 2L + 2 + j } ) 	= e_{ 2 L + 2 + j} 	&\mathrm{for} \quad 		&j = 1, 2, \dots 
\end{array}
\end{equation*}
Note that $w_{0} \in U(\unitize{A})$ and $w_{0}^{*} S w_{0} E_{L} = S E_{L}$.  Let $u' = w_{0}^{*} u w_{0}$.  Then 
$u' = E_{L} + ( E_{2L + 2} - E_{L} ) u' ( E_{2L + 2} -E_{L} ) + 1 - E_{2L + 2}$.

Set $u'' = p_{1} + u' w_{0}^{*} S w_{0} S^{*}$.  It is easy to check that $u'' S = w_{0}^{*} T w_{0}$.  Clearly, $u'' p_{1} = p_{1} = p_{1} u''$.  Note, for $1 < i \leq L$, $u' w_{0}^{*} S w_{0} S^{*} (e_{i}) = u' w_{0}^{*} S w_{0} (e_{i - 1}) = u' w_{0}^{*} S (e_{L + i}) = u' w_{0}^{*} (e_{L + i + 1}) = u'(e_{i}) = e_{i}$. 
Also, if $i > 2L + 4$, then $u' w_{0}^{*} S w_{0} S^{*} (e_{i}) = u' w_{0}^{*} S w_{0} (e_{i-1}) = u' w_{0}^{*} S(e_{i-1}) = u' w_{0}^{*} (e_{i}) = u' (e_{i}) = e_{i}$.  Therefore, $u'' = E_{L} + ( E_{2L + 4} - E_{L} ) u'' ( E_{2L + 4} - E_{L} ) + ( 1 - E_{2L + 4} )$.  So, we may assume that $u = u''$.

Let $u_{1} = ( E_{2L + 4} - E_{L} ) u'' ( E_{2L + 4} - E_{L} )$.  Note that we may identify $u_{1}$ as an element in $U(\matalg{L+4}{B})$.  By Lemma \ref{pert L : approx unitarily eq unitaries1}, there exist $m > 2L + 4$, a unitary $w_{1} \in \matalg{m}{B}$, and a unitary $v \in B = p_{L + 1} A p_{L + 1}$ such that
\begin{equation*}
\norm{ w_{1} (E_{L} + v + E_{m} - E_{L + 1} ) s_{m} w_{1}^{*} - ( E_{L} + u_{1} +E_{m} -E_{2L + 4} ) s_{m}  } \leq \frac{7 - 2 \pi}{L}.
\end{equation*}
By Lemma \ref{pert L : approx unitarily eq unitaries2}, there exists $w_{2} \in U(\matalg{2m}{B})$ such that 
\begin{equation*}
\norm{ w_{2} (E_{L} + v + E_{2m} - E_{L + 1} )s_{2m} w_{2}^{*} - (E_{L} + u_{1} + E_{2m} - E_{2L + 4} ) s_{2m} } \leq \frac{7}{L} 
\end{equation*}
and $w_{2} p_{2m} = p_{2m} w_{2} = p_{2m}$.  Set $w_{3} = w_{2} + 1 - E_{2m}$.  Then $w p_{2m} = p_{2m}$.  Thus
\begin{equation*}
w_{3} (E_{L} + v + 1 - E_{L + 1} ) S w_{3}^{*} (1 - E_{2m - 1} ) = (E_{L} + v + 1 - E_{L + 1}) S (1 - E_{2m-1}).
\end{equation*}
Note that $w_{3}^{*} E_{2m - 1} = E_{2m - 1} w_{3}^{*} E_{2m - 1}$ and $(E_{L} + v + E_{2m} - E_{L + 1}) s_{2m} E_{2m - 1} = (E_{L} + v + 1 - E_{L + 1}) S E_{2m - 1}$.  So, we have $(E_{L} + v + E_{2m} - E_{L + 1}) s_{2m} w_{3}^{*} E_{2m - 1} = (E_{L} + v + 1 - E_{L + 1}) S w_{3}^{*} E_{2m - 1}$.

Let $w = w_{3} \diag( v^{*} , \dots, v^{*}, 1, \dots )$,
where $v^{*}$ is repeated $L - 1$ times.   Then $w \in U(\unitize{A})$ such that 
{\normalsize
%\small
\begin{align*}
\norm{ w S w^{*} - u S}  &= \norm{ (w_{3} (E_{L} + v + 1 - E_{L+1})Sw_{3}^{*} - u S)E_{2m - 1} }\\
&\leq
  \norm{ w_{2} (E_{L} + v + E_{2m} - E_{L + 1}) s_{2m} w_{2}^{*} - (E_{L} + u_{1} + E_{2m} - E_{L + 4}) s_{2m}  } \\ 
    &\leq  \frac{7}{L} < \epsilon. 
\end{align*}
}
\end{proof}

\begin{lemma}\label{pert L : approx unitarily purely inf isometry}
Let $A$ be a separable purely infinite simple C*-algebra and let $s, t$ be two non-unitary isometries in $A$ (in $\unitize{A}$ if $A$ does not have a unit).  Then, for any $\epsilon > 0$, there exists $W \in U(A) \ (\mathrm{or} \ U(\unitize{A}))$ such that $\norm{ W^{*} s W - t } < \epsilon$. 
\end{lemma}

\begin{proof}
Suppose $A$ is unital.  Then the conclusion follows from Lemma 2.9 in \cite{inftoeplitz}.  Suppose $A$ is non-unital.  Since $A$ is a separable C*-algebra with real rank zero, there exist a projection $e \in A$ and non-unitary isometries $s', t' \in e A e$ such that $\norm{ s' + 1 - e - s} < \frac{\epsilon}{3}$ and $\norm{t' + 1 - e - t} < \frac{\epsilon}{3}$.  By the unital case, there exists $w \in U(e A e)$ such that $\norm{ w^{*} s' w - t'}  < \frac{\epsilon}{3}$.  Let $W = w + 1 -e \in U(\unitize{A})$.  Then $\norm{ W^{*} s W - t } < \epsilon$.
\end{proof}

\begin{lemma}\label{pert L : approx unitarily eq unitaries}
Let $E$ be an \ealg.  Then, for any $\epsilon > 0$ and for any $U_{1}, U_{2} \in U(E)$ satisfying $\spec( U_{1} ) = \spec( U_{1} ) = S^{1}$, $\pi ( U_{1} ) = \pi ( U_{2} )$, and $ [ U_{1} ] = [ U_{2} ]$ in $K_{1} (E)$, there exists $W \in U(\unitize{I})$ such that $\norm{ W^{*} U_{1} W - U_{2} } < \epsilon$.
\end{lemma}
 
\begin{proof}
By Proposition \ref{def P : diagonalization}, we may assume $U_{i} = \sum_{k = 1}^{\infty} \lambda_{k}^{(i)} e_{k}^{(i)}$ for $i = 1,2$, where the sum converges in the strict topology.  Using a similar argument as in \cite[2.10]{inftoeplitz}, we get $W \in U(E)$ such that $\norm{ W^{*} U_{1} W - U_{2} } < \epsilon$.
\end{proof}
 
\begin{lemma}\label{pert L : approx unitarily eq isometries trivial}
Let $I \in \PN$.  Let $E$ be a unital C*-subalgebra of $\multialg{I}$ which contains $I$ as an essential ideal.  Suppose $S_{1}, S_{2}$ are two non-unitary isometries such that  $\pi ( S_{1} ) = \pi  ( S_{2} )$ and $[ 1 - S_{1} S_{1}^{*} ] = [ 1 - S_{2} S_{2}^{*} ] = 0$ in $K_{0} (I)$.  Then, for any $\epsilon > 0$, there exists a unitary $W \in \unitize{I}$ such that $\norm{ W^{*} S_{1} W - S_{2} } < \epsilon$.
\end{lemma}

\begin{proof}
By Lemma \ref{def L : trivial ext indmap zero}, there exists $U \in U(E)$ and $\pi (U) = \pi ( S_2 )$.  By Proposition \ref{def P : diagonalization}, we may assume that there exists a sequence of mutually orthogonal projections $\{ e_{n} \}_{n = 1}^{\infty}$ in $I$ such that $U = \sum_{n = 1}^{\infty} U_{n}$, where $U_{n} \in U(e_{n} I e_{n})$ and the sum converges in the strict topology. 

Let $T = \sum_{n = 2} ^{\infty} U_{n} + v$, where the sum converge in the strict topology and $v \in e_{1} I e_{1}$ such that $v^{*} v = e_{1}$ and $v v^{*} \neq 0$.  Note that $T^{*} T = 1$, $T T^{*} \neq 1$, and $\pi (T) = \pi (S_{1})$.  Using the same argument as in \cite[2.11]{inftoeplitz}, there exists $W_{1} \in U(\unitize{I})$ such that $\norm{W_{1}^{*} S_{1} W_{1} - T} < \frac{\epsilon}{2}$.  Since $\pi(T) = \pi(S_{2})$, we again get $W_{2} \in U(\unitize{I})$ such that $\norm{ W_{2}^{*} S_{2} W_{2} - T } < \frac{\epsilon}{2}$.  Let $W = W_{1} W_{2}^{*}$.  Then $\norm{ W^{*} S_{1} W - S_{2} }  < \epsilon$.
\end{proof}

\section{THE INVARIANT}\label{invariant}

\begin{definition}\label{inv D : m-v proj}
Let $A$ be a C*-algebra.  Let $V(A)$ be the set of all Murray-von Neumann equivalence classes of projections in matrices (of all sizes) over $A$.  If addition on $V(A)$ is given by $[ p ]  + [ q ] = [ p \oplus q ] \ \mbox{for all}  \ [ p ], [ q ] \in V(A)$, 
then $V(A)$ is an abelian semigroup. 
\end{definition}

If $\varphi \in \Hom(A, B)$, then $V(\varphi)$ will denote the induced homomorphism.    

\begin{lemma}\label{inv L : kthypureinf}
Let $A$ be a C*-algebra such that $A$ has an approximate identity consisting of projections.  Then the map from the Grothendieck group of $V(A)$ to $K_{0}(A)$ is an isomorphism.
\end{lemma}

For a proof of the above lemma, see \cite[5.5.5]{blackadarB}.  If $A$ is a non-unital separable purely infinite simple C*-algebra, then $A$ has an approximate identity consisting of projections.  Throughout this section, $\kclass{p}$ will denote the image of $[p]$ in $K_{0}(A)$.

\begin{lemma}\label{inv L : projealg}
Let $E$ be an \ealg with $\quot{E} = C(S^{1})$.  Then for every projection $p \in \matalg{n}{E}$ not in $\matalg{n}{\ideal{E}}$, there exist a projection $e(p) \in \ideal{E}$ and a positive integer $n(p) \leq n$ such that $p \sim 1_{n(p)} \oplus e(p)$.  \end{lemma}

\begin{proof}
Since $\quot{E} = C(S^{1})$, there exists a projection $q \in \matcont{n}{S^{1}}$ such that $\pi(q) = 1_{m}$ and $p \sim q$.   So, we may assume $\pi(p) = 1_{m} = \pi(1_{m})$.  Since $\ideal{E}$ has an approximate identity consisting of projections, there exists a partial isometry $v \in \matalg{n}{E}$ such that $v^{*} v \leq 1_{m}$, $v v^{*} \leq p$, and $\pi(v) = 1_{m}$.  Let $e$ be a nonzero projection in $\ideal{E}$ such that $\kclass{e} = 0$ in $K_{0}(\ideal{E})$.  By \cite[1.5]{kthypureinf}, there exists $e' \in \ideal{E}$ such that $(1_{m} - v^{*} v) \oplus e' \sim e$.  Hence, $p \sim v^{*} v \oplus (p - v v^{*}) \sim 1_{m} \oplus e' \oplus ( p - v v^{*})$.
\end{proof}

\begin{notation}
\noindent (1) Denote the disjoint union of $X_{1}$ and $X_{2}$ by $X_{1} \sqcup X_{2}$. 

\noindent (2) Suppose $E$ is an \ealg such that $\quot{E} \cong \matcont{n}{S^{1}}$.  If $x \in K_{0}(\ideal{E})$, then denote the image of $x$ in $K_{0}(\ideal{E}) / \ran \indmap{E}$ by $\md{x}$.
\end{notation} 

\begin{theorem}\label{inv T : computation invariant proj}
(1) Let $E$ be a separable purely infinite simple C*-algebra.  Consider the abelian semigroup $\{ 0 \} \sqcup K_{0}(E)$, where $0 + x = x$ for all $x \in K_{0}(E)$.  Then $V(E) \cong \{ 0 \} \sqcup K_{0}(E)$, where the isomorphism sends $[0]$ to $0$ and sends $[p]$ to $\kclass{p}$ for every nonzero projection $p$.  

(2) Let $E$ be an \ealg with $\quot{E} = C(S^{1})$.  Consider the abelian semigroup $V = \{ 0 \} \sqcup K_{0} (\ideal{E}) \sqcup \left(\Z_{>0} \oplus \left(K_{0} (\ideal{E})/ \ran \delta_{1}^{E} \right)\right)$, where addition is defined as follows:  Addition in $K_{0} (\ideal{E})$ and in $\Z_{>0}\oplus \left(K_{0} (\ideal{E})/ \ran \delta_{1}^{E}\right)$ are the usual addition in those semigroups.  If $x \in K_{ 0 } (\ideal{E})$ and $ ( z_{1} , \md{z_{2}} ) \in \Z_{>0} \oplus \left( K_{0} (\ideal{E})/ \ran \delta_{1}^{E}\right)$, then $x + ( z_{1} , \md{z_{2} } ) = ( z_{1} , \md { z_{2} } + \md{ x } )$.  Suppose $\ftn{\alpha}{V(E)}{V}$ is defined by $\alpha([p]) = 0$ if $p = 0$; $\alpha([p]) = \kclass{p}$ if $p \in \matalg{n}{\ideal{E}}\setminus\{ 0 \}$; and $\alpha([p])= (n(p), \kclass{e(p)})$ otherwise, where $p \sim 1_{n(p)} \oplus e(p)$ is the decomposition given in Lemma \ref{inv L : projealg}.  Then $\alpha$ is a well-defined isomorphism.  Moreover, the natural map from $V(\ideal{E})$ to $V(E)$ is injective and $\alpha$ sends its image onto $\{ 0 \} \sqcup K_{0}(\ideal{E})$.

(3) Let $E$ be a trivial extension with $\quot{E} = C(S^{1})$.  Then $\alpha$ in (2) gives an isomorphism from $V(E)$ onto $\{ 0 \} \sqcup K_{0}(\ideal{E}) \oplus \Z_{\geq 0}$.  Moreover, the natural map from $V(\ideal{E})$ to $V(E)$ is injective and $\alpha$ sends its image onto $\{ 0 \} \sqcup K_{0}(\ideal{E})$.

(4) Let $E$ be an \ealg.  Then the map $x$ to $\diag(x, 0)$ induces an isomorphism from $V(E)$ onto $V(\matalg{n}{E})$.  If $p$ is a projection in $\matalg{n}{E}$ not in $\matalg{n}{\ideal{E}}$, then the inclusion from $p \matalg{n}{E} p$ to $\matalg{n}{E}$ induces an isomorphism from $V( p \matalg{n}{E} p)$ onto $V(\matalg{n}{E})$.  Hence, $V(E)$ is a finitely generated abelian semigroup.
\end{theorem}

\begin{proof}
(1) is a consequence of \cite[1.4]{kthypureinf} and Lemma \ref{inv L : kthypureinf}.

It is easy to check that if $p \in \ideal{E}$, $q \in E$, and $p \sim q$, then $q \in \ideal{E}$.  Therefore,
\begin{align*}
V(E) 		&= V(\ideal{E}) \sqcup \setof{[ p ]}{p \in \matalg{n}{E} \setminus \matalg{n}{\ideal{E}} \mbox{ for some } n } \\
		&\cong \{ 0 \} \sqcup K_{0}(\ideal{E}) \sqcup \setof{[ p ]}{p \in \matalg{n}{E} \setminus \matalg{n}{\ideal{E}} \mbox{ for some } n }.
\end{align*}

(2) and (3) now follows from the exactness of the six-term exact sequence in \kthy and Lemma \ref{inv L : projealg}.  The last statement of the theorem is clear.        
\end{proof}

Let $A$ be a C*-algebra. An element $s \in A$ is called \textbf{hyponormal} if $s^{*}s \geq ss^{*}$.  For a unital C*-algebra $A$, let $S_{n}(A)$ be the set of all nonzero hyponormal partial isometries in $\matalg{n}{A}$.  Let $S(A) = \bigcup _{n = 1 } ^{\infty} S_{n} (A) $, where we embed $S_{n}(A)$ into $S_{n + 1}(A)$ by sending $s$ to $\diag(s, 1)$.  

\begin{definition}\label{inv D : hyponormal}
For $v_{1}, v_{2} \in S_{n}(A)$, we write $v_{1} \simeq v_{2} $ if and only if there exists $m \in \Z_{\geq 0}$ such that $(1_{n} - v_{1}^{*} v_{1} + v_{1} ) \oplus 1_{m} $ is homotopic to $(1_{n} - v_{2}^{*} v_{2} + v_{2} ) \oplus 1_{m} $ in $S_{n + m} (A)$.  Set $\hypon{A}  = S(A) / \simeq$.  Let $\hclass{s}$ denote the equivalence class represented by $s$.  If addition is defined by $\hclass{u} + \hclass{v} = \hclass{ u \oplus v }$, then $\hypon{A}$ becomes an abelian semigroup.  If $A$ is a non-unital C*-algebra, then $\hypon{A} = \hypon{\unitize{A}}$.  
\end{definition}

If $\varphi \in \Hom(A, B)$, then $\hypon{\varphi}$ will denote the induced homomorphism.  

\begin{lemma}\label{inv L : facthyp}  
Let $A$ be a unital C*-algebra.  Then the following hold.

(1) For all projections $p, q \in \matalg{n}{A}$, we have $\hclass{p} = \hclass{q}$.  

(2) If $s$ is an isometry and $u$ is a unitary, then $\hclass{s u } = \hclass{s} + \hclass{u} = \hclass{u s}$.

(3) $\hypon{A} = \setof{\hclass{s}}{s \ \mbox{is an isometry in} \ \matalg{n}{A} \ \mbox{for some} \ n \in \Z_{>0}}$.
\end{lemma}     

The proof of the above lemma is easy and we leave it to the reader.

\begin{lemma}\label{inv L : equiv of iso approx the same}
Let $A$ be a unital C*-algebra.  Suppose that $s$ and $v$ are isometries in $A$ such that $\norm { s - v } < \frac{2}{4\sqrt{2} + 1}$.  Then $\hclass{s} = \hclass{v}$ in $\hypon{A}$.
\end{lemma}

\begin{proof}
Suppose $s$ or $v$ is in $U(A)$.  Then, $s, v \in U(A)$ and $s^{*} v \in U_{0}(A)$ since $\norm{ s - v } < \frac{2}{4\sqrt{2} + 1} < 2$.  Thus, $\hclass{v} = \hclass{s}$ in $\hypon{A}$.  Suppose $s$ and $v$ are non-unitary isometries.  Note that $\norm { (1- s s^{*}) - (1 - v v^{*} ) } < \frac{4}{4\sqrt{2} + 1}$.  Hence, there exists $w \in U(A)$ such that $\norm{ 1 - w } < \frac{ 4 \sqrt{ 2}  } { 4 \sqrt{2} + 1} $, $ w^{*} (1 - v v^{*} )w = 1 -s s^{*}$, and $\norm { w^{*} v w - s } <  2$.  Let $x = w^{*} v w$.  Then $x$ is an isometry with $ 1- x x^{*} = 1 - s s^{*}$.  Also, $u = x^{*}s \in U(A)$, $x u = s$, and $\norm{1 - u}  < 2$.  Hence,by Lemma \ref{inv L : facthyp}, $\hclass{s} = \hclass{w^{*} v w} = \hclass{v} \mbox{ in } \hypon{A}$.
\end{proof}

\begin{lemma}\label{inv L : embdhyp}
Let $A$ be a unital C*-algebra.  If $\hclass{u} = \hclass{s}$, where $u \in U(A)$ and $s$ is an isometry in $A$, then $s \in U(A)$.  

If $\ftn{\iota}{K_{1}(A)}{\hypon{A}}$ be the natural map, then 
\begin{equation*}
\hypon{A} = \iota(K_{1} (A)) \sqcup \setof{ \hclass{s} }{ s \mbox{ non-unitary isometry in } S(A)}.
\end{equation*}
Moreover, $\iota$ is injective and $\iota$ is an isomorphism whenever $A$ has cancellation.
\end{lemma}

\begin{proof}
The first part of the lemma follows from Lemma \ref{inv L : equiv of iso approx the same}.  If $A$ has cancellation, then every isometry is a unitary.  Hence, by Lemma \ref{inv L : facthyp}, $\iota$ is surjective.
\end{proof}

\begin{lemma}\label{inv L : stdshift}
Let $E$ be an \ealg such that $\quot{E} = C(S^{1})$.  Then there exist a projection $e \in \ideal{E}$, an isomorphism $\ftn{\ideal{\eta}}{\ideal{E}}{e \ideal{E} e \otimes \K}$, a monomorphism $\ftn{\eta}{E}{\multialg{e\ideal{E} e \otimes \K}}$, and $y \in E$ such that $\pi(y) = z$, $\eta(y)$ is the standard unilateral shift for $e \ideal{E} e \otimes \K = I$, and the following diagram commutes:
\begin{equation*}
\xymatrix{ 0 \ar[r] & \ideal{E} \ar[r] \ar[d]^{\ideal{\eta}} & E \ar[r] \ar[d]^{\eta} & \quot{E} \ar[r] \ar[d]^{\quot{\eta}} & 0 \\
		0 \ar[r] & I \ar[r] & \multialg{I} \ar[r] & \corona{I} \ar[r] & 0 }
\end{equation*}
\end{lemma}

\begin{proof}
Let $T$ be a non-unitary isometry such that $\pi(T) = z$.  Set $e' = 1 - T T^{*}$.  By \cite[2.6]{heralgs} and \cite[3.2(i)]{corona1}, there exists an isomorphism $\ftn{\gamma}{\ideal{E}}{e' \ideal{E} e' \otimes \K}$.  Hence,  $\gamma$ extends to an isomorphism $\gamma$ from $\multialg{\ideal{E}}$ onto $\multialg{e' \ideal{E} e' \otimes \K}$.  Set $\eta' = \gamma \vert_{E}$.  Let $S'$ be the standard unilateral shift for $e' \ideal{E} e' \otimes \K = I'$.  Let $\tau$ and $\tau'$ be the Busby invariants for the extensions $C^{*}(\eta'(T), I')$ and $C^{*}(S', I')$ respectively.  Since $[1 - S' (S')^{*} ] = [ e' ] = \indmap{E}([z])$, we have $[\tau] = [\tau']$ in $\Ext(C(S^{1}), \ideal{E})$.  By Proposition \ref{def P : absorbing extensions}, there exists a unitary $W \in \multialg{e' \ideal{E} e' \otimes \K}$ such that $\pi(W S' W^{*} ) = z$.  Let $S = W S^{'} W^{*}$, let $e = W e' W^{*}$, and let $\eta = \innerauto(W) \circ \eta'$.      
\end{proof}

\begin{lemma}\label{inv L : notuniso}
Let $E$ be an \ealg such that $\quot{E} \cong \matcont{n}{S^{1}}$.  Suppose $S$ and $T$ are two non-unitary isometries. Then $[ \pi(S) ] = [ \pi(T) ]$ in $K_{1}(\quot{E})$ if and only if $\hclass{S} = \hclass{T}$ in $\hypon{E}$.
\end{lemma}

\begin{proof}
If $\hclass{S} = \hclass{T}$, then by Lemma \ref{inv L : embdhyp}, $[ \pi(S) ] = [ \pi(T) ]$ in $K_{1}(\quot{E})$. 
Suppose $[ \pi(S) ] = [ \pi(T) ]$ in $K_{1}(\quot{E})$.  Then there exists $w \in U_{0}(E)$ such that $\pi(S) = \pi(w T)$.  By Lemma \ref{inv L : facthyp}, $\hclass{w T} = \hclass{T}$.  Hence, we may assume $\pi(S) = \pi(T)$.

Let $0 < \epsilon < \frac{2}{4\sqrt{2} + 1}$.  Suppose $\kclass{1 -S S^{*}} \neq 0$.  Then $C^{*}(S, \ideal{E})$ is a unital essential extension of $C(S^{1})$ by $\ideal{E}$.  By Lemma \ref{inv L : stdshift}, $\pi(T) = \pi(S) = \pi(S_{1})$, where $S_{1}$ is the standard unilateral shift of $e \ideal{E} e \otimes \K \cong \ideal{E}$.   Applying Lemma \ref{pert L : approx unitarily eq isometries} to $T$ and $S_{1}$, then to $S$ and $S_{1}$, there exists $v \in U(\unitize{\ideal{E}})$ such that $\norm{ v^{*} S v - T} < \epsilon$.  Suppose $\kclass{1 - S S^{*} } = 0$.  Then by Lemma \ref{pert L : approx unitarily eq isometries trivial}, there exists $v \in U(\unitize{\ideal{E}})$ such that $\norm{ v^{*} S v - T} < \epsilon$.  Hence, in either case, by Lemma \ref{inv L : equiv of iso approx the same} and Lemma \ref{inv L : facthyp}, $\hclass{S} = \hclass{T}$.
\end{proof}

Note that in the proof of the above lemma we proved the following.
\begin{proposition}\label{inv P : isoapprox}
Let $E$ be an \ealg such that $\quot{E} \cong \matalg{n}{C(S^{1})}$.  Let $\epsilon > 0$.  If $S$ and $T$ are two non-unitary isomerties in $E$ such that $\pi(S) = \pi(T)$, then there exists $w \in U(\unitize{\ideal{E'}})$ such that $\norm{w^{*} S w - T} < \epsilon$.
\end{proposition}

\begin{theorem}  \label{inv T : computation invariant k}
(1) Let $E$ be a purely infinite simple C*-algebra.  Consider the abelian semigroup $K_{1} (E) \sqcup \{ \hzero \}$, where addition is defined as follows: If $x \in K_{1}(E)$, then $x + n \hzero  =  \hzero$ for all $n \in \Z$.  Addition in $K_{1}(E)$ is the usual one.  Define $\ftn{\alpha}{\hypon{E}}{ K_{1}(E) \sqcup \{ \hzero \}}$ by $\alpha(\hclass{s}) = [s]$ if $s$ is a unitary and $\alpha(\hclass{s}) = \hzero$ otherwise.  Then $\alpha$ is a well-defined isomorphism.  

(2) Let $E$ be an \ealg with $\quot{E} = C(S^{1})$ and $E$ is a trivial extension.  Consider the abelian semigroup
$G = \setof{ (x , y) }{ x \in \hypon{\ideal{E}}  \ \mathrm{and} \ y \in \Z }$
with coordinatewise addition.  Then there exists an isomorphism from $\hypon{E}$ onto $G$.  Moreover, the map from $\hypon{\ideal{E}}$ to $\hypon{E}$ is injective and $\alpha$ sends its image onto $\setof{(x,0)}{ x \in \hypon{\ideal{E}}}$.

(3) Let $E$ be an \ealg with $\quot{E} = C(S^{1})$ and $E$ is a non-trivial extension.  Consider the abelian semigroup $K_{1} (E) \sqcup \Z$ with addition defined as follows: if $x \in K_{1} (E)$ and $n \in \Z$, then $x+n = m+n$ where $ m = \pi_{*, 1} (x) \in K_{1} (\quot{E}) \cong \set{Z} $.  Addition in $K_{1} (E)$ and in $\set{Z}$ are the usual addition in those semigroups.  Then there exists an isomorphism $\alpha$ from $\hypon{E}$ onto $K_{1}(E) \sqcup \Z$.  Moreover, $\alpha$ sends the image of the natural injective map from $\hypon{\ideal{E}}$ to $\hypon{E}$ onto $K_{1}(\ideal{E}) \sqcup \{ 0 \}$.

(4) Let $E$ be an \ealg.  Then the map $x$ to $\diag(x, 1)$ induces an isomorphism from $\hypon{E}$ onto $\hypon{\matalg{n}{E}}$.  If $p$ is a projection in $\matalg{n}{E}$ not in $\matalg{n}{\ideal{E}}$, then the inclusion from $p \matalg{n}{E} p$ to $\matalg{n}{E}$ induces an isomorphism from $\hypon{ p \matalg{n}{E} p}$ onto $\hypon{M_{n}(E)}$.  Hence, $\hypon{E}$ is a finitely generated abelian semigroup.
\end{theorem}

\begin{proof}
(1) follows from Lemma \ref{pert L : approx unitarily purely inf isometry}, Lemma \ref{inv L : equiv of iso approx the same}, and Lemma \ref{inv L : facthyp}.  

(2)  Let $s$ be a non-unitary isometry in $\matalg{n}{E}$.  Since $E$ is a trivial extension and $\pi(s)$ is a unitary in $\matalg{n}{C(S^{1})}$, there exists $u \in U(\matalg{n}{E})$ such that $\pi (s) = \pi (u)$.  Hence, $s u^{*}$ is a non-unitary isometry in $\matalg{n}{\unitize{\ideal{E}}}$.  So, $\hclass{s} = \hclass{s u^{*}} + \hclass{u}$ in $\hypon{E}$.  Note that this decomposition is unique. 

Let $\beta$ be the isomorphism from $K_{1}(E)$ onto $K_{1}(\ideal{E}) \oplus \Z$.  Define $\ftn{\alpha}{\hypon{E}}{\hypon{\ideal{E}} \oplus \Z}$ by $\alpha(\hclass{w}) = \beta([u])$ if $w$ is a unitary.  If $w$ is a non-unitary isometry in $\matalg{n}{E}$, then by the above observation $\hclass{w} = \hclass{s} + \hclass{u}$ for some $u \in U(\matalg{n}{E})$ and non-unitary isometry $s \in \matalg{n}{\unitize{\ideal{E}}}$.  Define $\alpha([w]) = (\hclass{s}, [u])$, where $\hclass{s}$ is now considered as an element in $\hypon{\ideal{E}}$.  Then $\alpha$ is a well-defined isomorphism with the desired property. 

(3)  By Lemma \ref{inv L : notuniso}, $\setof{ \hclass{s} }{ s \mbox{ is a non-unitary isometry in } S(E)} \cong \Z$, where the isomorphism is induce by $\pi$.  Define $\ftn{\alpha}{\hypon{E}}{K_{1}(E) \sqcup \Z}$ by $\alpha(\hclass{s}) = \hclass{s}$ if $s$ is a unitary and $\alpha(\hclass{s}) = [\pi(s)]$ otherwise.  Then $\alpha$ is a well-defined isomorphism with the desired property.

The last statement of the theorem is clear.
\end{proof}

\begin{corollary}\label{inv C : iso k}
Let $E$ be an \ealg such that $\quot{E} \cong \matcont{n}{S^{1}}$.  If $\hypo{E}$ denotes the Grothendieck group of $\hypon{E}$, then $\ftn{\pi}{E}{\quot{E}}$ induces an isomorphism from $\hypo{E}$ onto $K_{1}(\quot{E})$.
\end{corollary}

\begin{proof}
By Lemma \ref{inv L : embdhyp}, $\hypon{\quot{E}} \cong \Z$ and by Theorem \ref{inv T : computation invariant k}, $\hypo{E} \cong \Z$.  It is now easy to see that $\pi$ induces an isomorphism from $\hypo{E}$ onto $K_{1}(\quot{E})$.     
\end{proof}

\begin{definition}\label{inv D : subvoline}
Let $A$ be a C*-algebra and let 
\begin{equation*}
\voline{A} = \setof{( [ u^{*} u ], \hclass{u} )}{ u \in S(A)} \subset V(A) \oplus \hypon{A}.
\end{equation*}
Define $\ftn { \dmap } { \hypon{A} } { V(A) }$ by $d ( \hclass{u} ) = [ u^{*} u - u u^{*} ]$.  A homomorphism $\ftn{\alpha}{\voline{A}}{\voline{B}}$ consists of two homomorphisms $\ftn{\alpha_{v}}{V(A)}{V(B)}$ and $\ftn{\alpha_{k}}{\hypon{A}}{\hypon{B}}$ such that if $\alpha(\hclass{s}) = \hclass{v}$, then $\alpha_{v}([s^{*} s]) = [v^{*}v]$.  
\end{definition}

\begin{lemma}\label{inv L : vinvo}
Let $A$ be a C*-algebra.  Define $\ftn{\Theta_{v}^{A}}{V(A)}{\voline{A}}$ by $\Theta_{v}^{A}([p]) = ([p], 0)$ and define $\ftn{\Theta_{k}^{A}}{\hypon{A}}{\voline{A}}$ by $\Theta_{k}^{A}(\hclass{s}) = ( [s^{*} s], \hclass{s})$.  Then $\Theta_{v}^{A}$ and $\Theta_{k}^{A}$ are injective homomorphisms.  Using $\Theta_{v}^{A}$ and $\Theta_{k}^{A}$, we may identify $V(A)$ and $\hypon{A}$ as subsemigroups of $\voline{A}$.
\end{lemma}

The proof of the above lemma is easy and we leave it for the reader.

\begin{lemma}\label{inv L : disindmap}
Let $E$ be a finite direct sum of \ealgs.  Let  $\ftn{\pi}{E}{\quot{E}}$ be the quotient map.  Then $\dmap = \indmap{E} \circ \hypon{\pi}$.
\end{lemma}

\begin{proof}
This follows from the definitions of $d$ and $\indmap{E}$.
\end{proof}

\begin{definition}\label{inv D : hom inv V}
Define $\vstar{A}$ to be the set of triples
\begin{equation*}
\setof{ ( [ u^{*} u ], \hclass{ u } , d( \hclass{ u } ) )}{ u \in S(A) } \subset V(A) \oplus \hypon{A} \oplus V(A).
\end{equation*}

Let $A$ and $B$ be two C*-algebras.  A homomorphism $\ftn{ \eta } {\vstar{A}}{\vstar{B}}$ is a homomorphism $\ftn{\eta}{ \voline{A} } { \voline{B} }$ for which the following diagram commutes:
\begin{equation*}\label{inv E : hom inv}
\xymatrix{
\hypon{A} \ar[r]^{d} \ar[d]_{\eta_{k}} 			&	V(A) \ar[d]^{\eta_{v}} 		\\
\hypon{B} \ar[r]^{d}						&	V(B)
}
\end{equation*} 
\end{definition}
If $\varphi \in \Hom(A,B)$, then $\varphi$ induces a homomorphism $\ftn{ \vstar{\varphi} }{ \vstar{A} }{ \vstar{B} }$. 

\begin{lemma}\label{inv L : homvo}
Let $E$ and $E'$ be two finite direct sums of \ealgs.  

(1) Suppose $\ftn{\alpha}{\voline{E}}{\voline{E'}}$ is a homomorphism.  Then, $\alpha_{v}$ maps $V(\ideal{E})$ to $V(\ideal{E'})$ and $\alpha_{k}$ maps $\hypon{\ideal{E}}$ to $\hypon{\ideal{E'}}$.  

(2) If $\ftn{\eta}{\vstar{E}}{\vstar{E'}}$ is a homomorphism, then $\eta$ induces a homomorphism from $\vstar{\ideal{E}}$ to $\vstar{\ideal{E'}}$.  Also, if $\ftn{\iota}{K_{1}(E)}{\hypon{E}}$ and $\ftn{\iota'}{K_{1}(E')}{\hypon{E'}}$ are the injective homomorphisms given in Lemma \ref{inv L : embdhyp}, then $\ftn{\eta_{k}}{\hypon{E}}{\hypon{E'}}$ maps $\iota(K_{1}(E))$ to $\iota'(K_{1}(E'))$.
\end{lemma} 

\begin{proof}
(1) Using the identifications in Theorem \ref{inv T : computation invariant proj} and Theorem \ref{inv T : computation invariant k} and all semigroup homomorphisms are assumed to preserve identities, one easily checks that $\alpha_{v}$ maps $V(\ideal{E})$ to $V(\ideal{E'})$ and $\alpha_{k}$ maps $\hypon{\ideal{E}}$ to $\hypon{\ideal{E'}}$.

(2) The first part of (2) follows from (1).  Let $\ftn{\pi}{E}{\quot{E}}$ and $\ftn{\pi'}{E'}{\quot{E'}}$ be the quotient maps.  Supppose $u \in U(\matalg{n}{E})$ and $\eta_{k}(\hclass{u}) = \hclass{s}$ for some isometry $s \in \matalg{m}{E'}$.  Then $0 \sim 1_{m} - s s^{*}$.  Hence, $s s^{*} = 1_{m}$.  Thus, $\eta_{k}$ maps $\iota(K_{1}(E))$ to $\iota'(K_{1}(E'))$.
\end{proof}

Let $E_{1}$ and $E_{2}$ be finite direct sums of \ealgs.  Then $E_{i}$ is an extension of $\quot{E_{i}}$ by $\ideal{E_{i}}$.  By Proposition \ref{def P : liftingproj}, $\expmap{E_{i}} = 0$.  So, the six-term exact sequence in \kthy associated to $E_{i}$ has the form
\begin{equation*}
0 \to K_{1} (\ideal{E_{i}})  \to K_{1} (E_{i}) \to K_{1} (\quot{E_{i}}) \to K_{0} (\ideal{E_{i}})  \to K_{0} (E_{i}) \to K_{0} (\quot{E_{i}}) \to 0
%\scalebox{.77}{\xymatrix{0 \ar[r] 	& K_{1} (\ideal{E_{i}})  \ar[r] & K_{1} (E_{i}) \ar[r]  & K_{1} (\quot{E_{i}}) \ar [r] & K_{0} (\ideal{E_{i}})  \ar[r] & K_{0} (E_{i}) \ar[r]	& K_{0} (\quot{E_{i}}) \ar[r]	 & 0} }
\end{equation*}
Denote this exact sequence by $\sixk(E_{i})$.  A map from $\sixk(E_{1})$ to $\sixk(E_{2})$ consists of six group homomorphisms $\alpha = \{ \alpha_{i} \}_{i=1}^{6}$ such that the following diagram commutes:

\begin{equation*}
\scalebox{.77}{
\xymatrix{
0 \ar[r] 	& K_{1} (\ideal{E_{1}})  \ar[r] \ar[d]_{\alpha_{1}}	& K_{1} (E_{1}) \ar[r] \ar[d]_{\alpha_{2}}  & K_{1} (\quot{E_{1}}) \ar [r] \ar[d]_{\alpha_{3}} 	& K_{0} (\ideal{E_{1}})  \ar[r] \ar[d]_{\alpha_{4}}	& K_{0} (E_{1}) \ar[r] \ar[d]_{\alpha_{5}}	& K_{0} (\quot{E_{1}}) \ar[r] \ar[d]_{\alpha_{6}}	& 0 \\
0 \ar[r] 	& K_{1} (\ideal{E_{2}})  \ar[r] 	& K_{1} (E_{2}) \ar[r] 		& K_{1} (\quot{E_{2}}) \ar [r] 	& K_{0} (\ideal{E_{2}}) \ar [r] 	& K_{0} (E_{2}) \ar[r] 		& K_{0} (\quot{E_{i}}) \ar[r] 	& 0
}}
\end{equation*}  

Let $E$ be a unital extension of $\quot{E}$ by $\ideal{E}$.  Set 
\begin{equation*}
\Gamma = \setof{ ( x, y) } { x \in K_{0}( \quot{E} )_{+}, y \in K_{1} (\quot{E}); \mbox{ if } x = 0, \mathrm{then} \ y = 0 }.
\end{equation*}
Let $K_{*} ( \quot{E})$ denote the graded group $K_{0} (\quot{E}) \oplus K_{1} (\quot{E})$ with the partial order generated by $\Gamma$.  By Proposition \ref{def P : induce hom}, any homomorphism $\ftn{\varphi}{E_{1}}{E_{2}}$ between two direct sums of \ealgs induces a map $\{ \alpha_{i} \}_{i=1}^{6}$ from $\sixk(E_{1})$ to $\sixk(E_{2})$ such that $\alpha_{6} \oplus \alpha_{3}$ preserves the order.

Using a similar method as in Section 1.16 in \cite{inftoeplitz} we get the following.

\begin{proposition}\label{inv P : hom inv}
Let $E_{1}$ and $E_{2}$ be two finite direct sums of \ealgs.  Suppose $\ftn{\eta}{\vstar{E_{1}}}{\vstar{E_{2}}}$ is a homomorphism.  Then $\eta$ induces a map $\{ \alpha_{i} \}_{i = 1}^{6}$ from $\sixk(E_{1})$ to $\sixk(E_{2})$ 
such that $\alpha_{6} \oplus \alpha_{3}$ preserves the order. 
\end{proposition}

Let $E$ be an \aealg.  It is easy to check that $\expmap{E} = 0$.  Hence, if $E$ and $E'$ are two \aealgs, then a map from $\sixk(E)$ to $\sixk(E')$ is defined exactly  the same way as for \ealgs.  Note that by Proposition \ref{def P : induce hom aealg}, any homomorphism $\ftn{\varphi}{E}{E'}$ between two \aealgs induces a map from $\sixk(E)$ to $\sixk(E')$.

\begin{proposition}\label{inv P : interwine diag invariant V}
Let $E = \dirlim ( E_{i} , \varphi_{i, i + 1})$ and $E' = \dirlim ( E_{i}' , \varphi_{i, i+1}' )$ be unital \aealgs.  Let $\ftn{\alpha}{\vstar{E}}{\vstar{E'}}$ be a homomorphism.  

(1)  $\alpha$ induces a commutative diagram
\begin{equation*}
\xymatrix{
\vstar{E_{1}}  \ar[r] \ar[d]_{\alpha^{(1)}}& \vstar{E_{2}}  \ar[r] \ar[d]_{\alpha^{(2)}}& \cdots  \ar[r] & \vstar{E} \ar[d]_{\alpha} \\
\vstar{E_{m_{1}}'} \ar[r]				  & \vstar{E_{m_{2}}'} 	\ar[r]				&\cdots  \ar[r] & \vstar{E'}		      }
\end{equation*}
for some increasing sequence of  natural numbers $\{ m_{k} \}_{k = 1}^{\infty}$.

(2)  There exists a unique map $\{ \alpha_{i} \}_{i = 1}^{6}$ from $\sixk(E)$ to $\sixk(E')$ induced by $\alpha$.

Furthermore, the map $\alpha_{6} \oplus \alpha_{3}$ from $K_{*} ( \quot{E} )$ to $K_{*} (\quot{E'})$ preserves the order.  If $\alpha_{v}([\unit{E}]) = [\unit{E'}]$, then $\alpha^{(i)}$ may be chosen such that $\alpha_{v}^{(i)} ([ \unit{E_{i}} ]) = [ \unit{E_{m_{i}}'} ]$ for all $i$.
\end{proposition}

\begin{proof}
Let $[0]$ denote the identity of $V(E)$, $V(E')$, $V(E_{n})$, and $V(E_{n}')$.  Note that $\vstar{E} = \dirlim ( \vstar{E_{i}} , \vstar{\varphi_{i, i + 1}})$ and $\vstar{E'} = \dirlim ( \vstar{E_{i}'} , \vstar{\varphi_{i, i + 1}'})$.  Denote the maps from $\hypon{E_{n}}$ to $V(E_{n})$,  from $\hypon{E_{n}'}$ to $V(E_{n}')$, from $\hypon{E}$ to $V(E)$, and from $\hypon{E'}$ to $V(E')$ by $d_{n}$, $d_{n}'$, $d$, and $d'$ respectively.  We will show that there exists $m_{1} \in \Z_{>0}$ and a homomorphism $\ftn{ \alpha^{(1)} } {\vstar{E_{1}}}{\vstar{E_{m_{1}}}}$ such that $\alpha \circ V_{*}(\varphi_{1,\infty}) = V_{*}(\varphi_{m_{1},\infty}') \circ \alpha^{(1)}$.  Note that we may assume $E_{1}$ has only one summand and $E_{1}$ is note a unital purely infinite simple C*-algebra. 

Case 1: Suppose $E_{1}$ is a non-trivial extension.  Then  $\hypon{E_{1}} \cong K_{1} (E_{1}) \sqcup \Z$ and $V(E_{1}) \cong \{ [0] \} \sqcup K_{0} ( \ideal{E_{1}} ) \sqcup \Bigl(\bigl(K_{0} (\ideal{E_{1}})/ \ran \indmap{E_{1}}\bigr) \oplus \Z_{>0} \Bigr)$.  Let $s_{1}$, $s_{2}$, and $s_{3}$ be non-unitary isometries in $E_{1}$ such that $\hclass{s_{1}} = 1$, $\hclass{s_{2}} = -1$, and $\hclass{s_{3}} = 1 - 1 = 0$.  

Suppose $\alpha_{v} \circ V(\varphi_{1,\infty}) \vert_{V(\ideal{E_{1}})} = [0]$.  For each $n \in \Z_{>0}$, choose $x_{n} \in V(E_{n}')$ such that $V(\varphi_{n,\infty}')(x_{n}) =  \alpha_{v} \circ V(\varphi_{1,\infty} )((0,1))$.  Define $\alpha_{1,n}$ from $V(E_{1})$ to $V(E_{n})$ by $\alpha_{1,n} \vert_{V(\ideal{E_{1}})} = [0]$ and $\alpha_{1,n} ((a, k)) = k x_{n}$.  It is clear that $\alpha_{1, n}$ is a homomorphism such that $\alpha_{v} \circ V(\varphi_{1,\infty}) = V(\varphi_{n, \infty}') \circ \alpha_{1, n}$.

Since $K_{1} (E_{1})$ is a finitely generated abelian group, there exists $n_{2} \in \Z_{>0}$ such that for all $n \geq n_{2}$, there exists a homomorphism $\ftn{\gamma_{1, n}}{ K_{1} (E_{1}) }{ K_{1}(E_{n}') }$ with 
$\alpha_{k} \circ \hypon{\varphi_{1, \infty}} \vert_{K_{1} (E_{1}) } = (\varphi_{n, \infty}' )_{*, 1} \circ \gamma_{1, n}$.  Note that $\ran (\alpha_{k} \circ \hypon{\varphi_{1,\infty}} ) \subset K_{1}(E')$ since $\alpha_{v} \circ V(\varphi_{1, \infty}) \vert_{V(\ideal{E_{1}})} = [0]$.  Let $(\alpha_{k} \circ \hypon{\varphi_{1,\infty}})(\hclass{s_{1}}) = \hclass{t_{1}'}$ and $(\alpha_{k} \circ \hypon{\varphi_{1,\infty}})(\hclass{s_{2}}) = \hclass{t_{2}'}$, where $t_{1}', t_{2}' \in U(\matalg{n}{E'})$.  Note that $\hclass{t_{1}'} + \hclass{t_{2}'} = 0$ in $K_{1}(E')$.  Therefore, there exists $n_{2} \in \Z_{>0}$ such that for all $n \geq n_{2}$, there exist $y_{n,i} \in K_{1}(E_{n}')$ with $(\alpha_{k} \circ \hypon{\varphi_{1,\infty}})(\hclass{s_{i}}) = \hypon{\varphi_{n,\infty}'}(y_{n,i})$ for $i = 1,2$ and $y_{n,1} + y_{n,2} = 0$ in $K_{1}(E_{n}')$.  Define $\ftn{ \beta_{1, n} }{ \hypon{E_{1}}}{ \hypon{E_{n}'} }$ by $\beta_{1, n} \vert_{ K_{1} (E_{1}) } = \gamma_{1, n}$ and $\beta_{1, n} ( \ell_{1} \hclass{s_{1}} + \ell_{2} \hclass{s_{2}} ) =  \ell_{1} y_{n,1} + \ell_{2} y_{n,2}$ for all $\ell_{1} , \ell_{2} \in \Z_{>0}$. 

Choose $m_{1} = \max\{ n_{1} , n_{2} \}$.  Set $\alpha^{(1)} = (\alpha_{1,m_{1}} , \beta_{1, m_{1}})$.  Then $\ftn{\alpha^{(1)}}{\voline{E_{1}}}{\voline{E_{m_{1}}'}}$ is a homomorphism.  Note that $(d_{m_{1}}' \circ \beta_{1, m_{1}})(\hclass{s_{i}}) = (\alpha_{1,m_{1}} \circ d_{1})(\hclass{s_{i}}) = [0]$.  Therefore, $\alpha^{(1)}$ is the desired homomorphism.

Suppose $\alpha_{v} \circ V(\varphi_{1,\infty}) \vert_{V(\ideal{E_{1}})} \neq [0]$.  It is easy to see that there exists $n_{1} \in \Z_{>0}$ such that for all $n \geq n_{1}$, there exists a homomorphism $\ftn{ \alpha_{1}^{ I, n}}{  \vstar{\ideal{E_{1}}}}{ \vstar{\ideal{E_{n}'}}}$ with $\vstar{\ideal{\varphi_{n,\infty}'} } \circ \alpha_{1}^{I,n} = \alpha \vert_{\vstar{\ideal{E_{1}}}} \circ \vstar{\ideal{\varphi_{1,\infty}}}$.  By Proposition \ref{inv P : hom inv}, $\alpha_{1}^{I,n}$ induces a homomorphism $\ftn{\lambda_{i,n}}{K_{i}(\ideal{E_{1}})}{K_{i}(\ideal{E_{n}'})}$ for $i=0,1$.  Also, $\alpha_{1}^{I,n}$ induces two 
homomorphisms $\ftn{\alpha_{1}^{I,n,v}}{V(\ideal{E_{1}})}{V(\ideal{E_{n}'})}$ and $\ftn{\alpha_{1}^{I,n,k}}{\hypon{\ideal{E_{1}}}}{\hypon{\ideal{E_{n}'}}}$.  Since $\ran(\alpha_{1}^{I,n, v} \circ d_{1})$ is finitely generated, we may choose $n_{1}$ such that $\ran(\alpha_{1}^{I,n,v} \circ d_{1}) \subset \ran d_{n}'$ for all $n \geq n_{1}$.  Since $V(E_{1})$ and $\hypon{E_{1}}$ are finitely generated, there exists $n_{2} \geq n_{1}$ such that $\ran(\alpha \circ \vstar{\varphi_{1,\infty}}) \subset \ran(\vstar{\varphi_{n,\infty}'})$ for all $n \geq n_{2}$.  Note that $\ran(\alpha_{1}^{I,n,v} \circ d_{1}) \subset \ran d_{n}'$ for all $n \geq n_{2}$.  Therefore, for all $n \geq n_{2}$, the homomorphism $\lambda_{0,n}$ induces a homomorphism $\ftn{ \widetilde{\lambda}_{0,n} }{ K_{0} (\ideal{E_{1}}) / \ran \indmap{E_{1}} }{ K_{0} (\ideal{E_{n}'}) / \ran \indmap{ E_{n}' }  }$ such that the diagram commutes
\begin{equation*}
\scalebox{.9}{\xymatrix{
K_{0} ( \ideal{E_{1}} ) \ar[r]^{\lambda_{0,n}} \ar[d]	& K_{0} (\ideal{E_{n}'}) \ar[r] 	& K_{0} (\ideal{E_{n}'}) / \ran \indmap{ E_{m_{1}}' }  \\
 K_{0} (\ideal{E_{1}}) / \ran \indmap{E_{1}}  \ar[rru]_{\widetilde{\lambda}_{0,n}} &                             &    
}}
\end{equation*}

\noindent For each $n \geq n_{2}$, choose $x_{n} \in V(E_{n}')$ such that $\alpha_{v} \circ V(\varphi_{1,\infty})((0,1)) = V(\varphi_{n,\infty}')(x_{n})$.  Define $\ftn{\alpha_{1,v,n}}{V(E_{1})}{V(E_{n}')}$ by $\alpha_{1,v,n} \vert_{V(\ideal{E_{1}})} = \alpha_{1}^{I, n, v}$ and $\alpha_{1,v, n}((a, k)) = \widetilde{\lambda}_{0,n}(a) + k x_{n}$ for $k \in \Z_{> 0}$.  Then, $\alpha_{1,v,n}$ is a homomorphism such that $V(\varphi_{n,\infty}') \circ \alpha_{1,v,n} = \alpha_{v} \circ V(\varphi_{1,\infty})$.

Since $K_{1}(E_{1})$ is finitely generated, there exists $n_{3} \geq n_{2}$ such that for all $n \geq n_{3}$ we have a homomophism $\ftn{\beta_{1,n}}{K_{1}(E_{1})}{K_{1}(E_{n}')}$ with $\alpha_{k} \circ \hypon{\varphi_{1,\infty}} \vert_{K_{1}(E_{1})} = (\varphi_{n,\infty}')_{*,1} \circ \beta_{1,n}$ and $\beta_{1, n} \vert_{K_{1}(\ideal{E_{1}})} = \lambda_{1,n}$.  Also, there exists $n_{4} \geq n_{3}$ such that for all $n \geq n_{4}$, there exist $y_{1,n}, y_{2,n} \in \hypon{E_{n}'}$ with $(\alpha_{k} \circ \hypon{\varphi_{1,\infty}})(\hclass{s_{i}}) = \hypon{\varphi_{n,\infty}'(y_{i,n})}$ and $y_{1,n} + y_{2,n}$ is the identity of the subsemigroup $\hypon{E_{n}'}\setminus K_{1}(E_{n}')$.  

Define $\ftn{\alpha_{1,k,n}}{\hypon{E_{1}}}{\hypon{E_{n}'}}$ by $\alpha_{1,k} \vert_{K_{1}(E_{1})} = \beta_{1,n}$ and $\alpha_{1,k,n}(\ell_{1} \hclass{s_{1}} + \ell_{2} \hclass{s_{2}}) = \ell_{1} y _{1,n} + \ell_{2} y_{2, n}$ for all $\ell_{1}, \ell_{2} \in \Z_{>0}$.  Then $\alpha_{1,k,n}$ is a homomorphism such that $\hypon{\varphi_{n,\infty}'} \circ \alpha_{1,k,n} = \hypon{\varphi_{1,\infty} } \circ \alpha_{k}$.  Note that there exists $m_{1} \geq n_{4}$ such that $(V(\varphi_{n_{4},m_{1}}') \circ \alpha_{1,v,n_{4}} \circ d_{1})(\hclass{s_{i}}) = (V(\varphi_{n_{4},m_{1}}') \circ d_{n_{4}}')(y_{i,n_{4}})$.  Hence, $\alpha^{(1)} = \left(V(\varphi_{n_{4},m_{1}}') \circ \alpha_{1,v,n_{4}}, \hypon{\varphi_{n_{4}, m_{1}}'} \circ \beta_{1,n_{4}} \right)$ is the desired homomorphism.

Case 2:  Suppose $E_{1}$ is a trivial extension.  Then $\hypon{E_{1}} \cong \hypon{\ideal{E_{1}}} \oplus \Z$ and $V(E_{1}) \cong \{ [0] \} \sqcup \left( K_{0}(\ideal{E_{1}}) \oplus \Z_{\geq 0}\right)$.  This case is proved in a similar fashion as in Case 1 but it is easier. 

Next, starting with $E_{2}$, there exist $m_{2}' \geq m_{1}$ and a homomorphism $\ftn{\beta^{(2)}}{\vstar{E_{2}}}{\vstar{E_{m_{2}'}'}}$ such that $\alpha \circ \vstar{\varphi_{2,\infty}} = \vstar{\varphi_{m_{2}', \infty}'} \circ \beta^{(2)}$.  Hence, there exists $m_{2} >  m_{2}'$ such that $\vstar{\varphi_{m_{2}', m_{2}}'}\circ \vstar{\varphi_{m_{1}, m_{2}'}'} \circ \alpha^{(1)} = \vstar{\varphi_{m_{2}', m_{2}}'} \circ \beta^{(2)} \circ \vstar{\varphi_{1,2}}$.  Let $\alpha^{(2)} = \vstar{\varphi_{m_{2}', m_{2}}'} \circ \beta^{(2)}$.  Then, the following diagram commutes:
\begin{equation*}
\scalebox{.9}{
\xymatrix{ \vstar{E_{1}} \ar[d]^{\alpha^{(1)}} \ar[rr]^{\vstar{\varphi_{1,2}}} & & \vstar{E_{2}} \ar[d]^{\alpha^{(2)}} \ar[rr]^{\vstar{\varphi_{2,\infty}} }& & \vstar{E} \ar[d]^{\alpha} \\
\vstar{E_{m_{1}}'} \ar[rr]_{\vstar{\varphi_{m_{1}, m_{2}}'}} & &\vstar{E_{m_{2}}'} \ar[rr]_{\vstar{\varphi_{m_{2},\infty}'}} & &\vstar{E'} }}
\end{equation*} 

Continuing this process, we get the desired result.
\end{proof} 

Let $C_{n}$ be the mapping cone of the degree $n$ map $\ftn{\theta_{n}}{C_{0}((0,1))}{C_{0}((0,1))}$.  Then $C_{n} \in \bootstrap$, $K_{0}(C_{n}) = \Z / n\Z$, and $K_{1}(C_{n}) = 0$.  The \textbf{total \kthy}of $A$ is defined to be $\totalk{A} = \bigoplus_{n = 0}^{\infty} \kmod{*}{A}{n}$, where $\kmod{*}{A}{n} = K_{*}(A \otimes C_{n})$ for $n \geq 2$, $\kmod{*}{A}{1} = K_{*}(A)$, and $\kmod{*}{A}{0} = 0$.  It is a $\Z / 2 \Z \times \Z_{\geq 0}$ graded group.  Let $\Lambda$ denote the category of Bockstein maps.  Denote the group of all $\Z / 2\Z \times \Z_{\geq 0}$ graded group homomorphisms which are $\Lambda$-linear by $\Hom_{\Lambda} ( \totalk{A} , \totalk{B})$.  See \cite[Section 4]{approxhom} for more details.   
  
Consider the following extension of C*-algebras $0 \to B \overset{\varphi}{\to}  E \overset{\psi}{\to}  A \to 0$.  Let $\varphi_{n} = \varphi \otimes \id_{C_{n}}$ and let $\psi_{n} = \varphi \otimes \id_{C_{n}}$.  Since $C_{n}$ is amenable, we have the following six-term exact sequence:
\begin{equation*}
\scalebox{.85}{
\xymatrix{ 
\kmod{0}{B}{n}	\ar[rr]^{ (\varphi_{n})_{*, 0} } &  & \kmod{0}{E}{n} \ar[rr]^{ (\psi_{n})_{*,0}} &  &\kmod{0}{A}{n} \ar[d]^{\expmap{E\otimes C_{n} }}\\
 \kmod{1}{A}{n} \ar[u]^{\indmap{ E \otimes C_{n}}} &  &\kmod{1}{E}{n} \ar[ll]^{( \psi_{n})_{*,1} } &  &\kmod{1}{B}{n} \ar[ll]^{ (\varphi_{n})_{*,1} }
}}
\end{equation*}

It is easy to check that if $E$ is a finite direct sum of \ealgs or an \aealg, then $\expmap{E \otimes C_{n} } = 0$ for all $n$.  
Hence, if $E$ is an \ealg or an \aealg, then set
{\footnotesize
\begin{equation*}
\totalkv{E} = \left( \totalk{B} , \ \totalk{E},\  \totalk{A},\  \bigoplus_{ n = 0 }^{\infty} ( (\varphi_{n})_{*,0} \oplus  (\varphi_{n})_{*,1}),\  \bigoplus_{ n = 0 }^{\infty} ( (\psi_{n})_{*} \oplus (\psi_{n})_{*,1}),\  \bigoplus_{ n = 0 }^{\infty}  \indmap{ E \otimes C_{n} } \right).
\end{equation*}
}

\begin{definition}\label{inv D : hom invariant totalk}
Let $E$ and $E'$ be two finite direct sums of \ealgs or $E$ and $E'$ are \aealgs.  Then a homomorphism $\ftn{ \eta } {\totalkv{E}}{\totalkv{E'}}$ is a system of $\Lambda$-linear maps, $\ftn{ \eta_{1} } { \totalk{\ideal{E}} } { \totalk{\ideal{E'}} }$, $\ftn{ \eta_{2} } { \totalk{ E } }{\totalk{E'}}$, and $\ftn{ \eta_{3} } { \totalk{\quot{E}} }{\totalk{\quot{E'}}}$, making the obvious diagrams commute.
\end{definition}

The invariant used to classify all unital \subaealgs with real rank zero is
$\vtilde{E}$.   A homomorphism $\ftn{\eta}{\vtilde{E_{1}}}{\vtilde{E_{2}}}$, where each $E_{i}$ is a finite direct sum of \ealgs, is a system of two homomorphisms $\ftn{\eta_{1}}{ \vstar{E_{1}} } { \vstar{E_{2}} }$ and $\ftn{\eta_{2}}{\totalkv{E_{1}}}{\totalkv{E_{2}}}$ such that the following diagrams commute:
{\scriptsize
\begin{equation*}
\xymatrix{
\totalk{\ideal{E_{1}}}	\ar[r]	\ar[d]_{\eta_{2}}	& K_{*} (\ideal{E_{1}}) \ar[d]_{(\alpha_{1}, \alpha_{4})} & \totalk{E_{1}}	\ar[r]	\ar[d]_{\eta_{2}}	& K_{*} (E_{1}) \ar[d]_{(\alpha_{2}, \alpha_{5})} &\totalk{\quot{E_{1}}}	\ar[r]	\ar[d]_{\eta_{2}}	& K_{*} (\quot{E_{1}}) \ar[d]_{(\alpha_{3}, \alpha_{6})}\\
	\totalk{\ideal{E_{2}}}	\ar[r]	& K_{*} (\ideal{E_{2}}) &\totalk{E_{2}} \ar[r]	& K_{*} (E_{2})  		&\totalk{\quot{E_{2}}}	\ar[r]				& K_{*} (\quot{E_{2}})
}
\end{equation*}
}

\noindent where the horizontal maps are the projection maps and $\{ \alpha_{i} \}_{ i = 1}^{6}$ are the unique maps (in Proposition \ref{inv P : hom inv}) induced by $\eta_{1}$.  If $E$ and $E'$ are two unital \aealgs, then a homomorphism from $\vtilde{E}$ to $\vtilde{E'}$ is defined similarly but now we use Proposition \ref{inv P : interwine diag invariant V} to get the unique maps $\{ \alpha_{i} \}_{i = 1}^{6}$.  By Proposition \ref{def P : induce hom} and by Proposition \ref{def P : induce hom aealg}, if $\varphi \in \Hom(E,E')$, $E$ and $E'$ are \ealgs or \aealgs, then $\varphi$ induces a homomorphism $\ftn{\vtildeh{\varphi}}{\vtilde{E}}{\vtilde{E}}$.

\section{THE UNIQUENESS THEOREM}\label{uniqueness}

\subsection{AUTOMORPHISMS}

\begin{definition}\label{uniq D : automorph group}
Let $E$ be a C*-algebra.  Denote the group of all automorphism of $E$ by $\Aut (E)$.  The topology on $\Aut(E)$ will be the norm topology, i.e. $\norm{\alpha} = \sup_{\norm{a} \leq 1} \norm{\alpha(a)}$.  Let $\Aut_{0}(E)$ denote the set of all $\alpha \in \Aut(E)$ that are in the same path component as $\id_{E}$ with the norm topology.
\end{definition}

\begin{theorem}\label{uniq T : connected id approx inner}(\cite[3.2]{stableapprox})
Let $E$ be a C*-algebra.  Then every automorphism $\alpha \in \Aut_{0}(E) \subset \Aut(E)$ is approximately inner. 
\end{theorem}    

Let $E$ be a separable C*-algebra.  Then the above theorem follows from the following facts:  (1) every $\alpha \in \Aut_{0}(E)$ is a product of derivable automorphisms (\cite[8.7.8]{autogrps}) and (2) every derivation is approximately inner (\cite[8.6.12]{autogrps}).  This observation was made by Lin in \cite[2.1]{idealcuntz}.  Lin then showed in \cite[3.2]{stableapprox} that the general case can be reduced to the separable case.  

\begin{corollary}\label{uniq C : aut same comp approx eq}
Suppose $\alpha, \beta \in \Aut (E)$.  Suppose that $\beta^{-1} \circ \alpha \in \Aut_{0} (E)$ and $\beta$ is approximately inner.  Then $\alpha$ is approximately inner.
\end{corollary}

\begin{figure}[h]
\scalebox{.5}{
\xymatrix{                               &  K_{0} (\ideal{E})  \ar[ddd] \ar[rr] &                             &  K_{0} (E) \ar[ddd] \ar[rr] &                               &  K_{0} (\quot{E}) \ar[ld]^{\expmap{E}} \ar[ddd] \\
           K_{1} (\quot{E})  \ar[ru]^{\indmap{E}} \ar[ddd]  &                               & K_{1} (E)  \ar[ll] \ar[ddd] &                              & K_{1} (\ideal{E}) \ar[ll] \ar[ddd] &                      \\
                                         &                               &                             &                              &                            &                       \\
                                        & K_{0} (\ideal{E} \rtimes_{\ideal{\alpha}} \Z) \ar[rr] \ar[ddd] &       & K_{0} ( E \rtimes_{\alpha} \Z ) \ar[rr] \ar[ddd] &    & K_{0} ( \quot{E} \rtimes_{\quot{\alpha}} \Z ) \ar[ld]^(.45){\expmap{E\rtimes_{\alpha} \ Z}} \ar[ddd]           \\
           K_{1} ( \quot{E} \rtimes_{\quot{\alpha}} \Z ) \ar[ur]^(.48){\indmap{E\rtimes_{\alpha} \Z}}  \ar[ddd] &     &  K_{1} ( E \rtimes_{\alpha} \Z ) \ar[ll] \ar[ddd] &       &  K_{1} ( \ideal{E} \rtimes_{\ideal{\alpha}} \Z ) \ar[ll]  \ar[ddd] &               \\
\\
                                       & K_{1} (\ideal{E}) \ar[rr] &   & K_{1} (E) \ar[rr] &  & K_{1}(\quot{E}) \ar[ld]^{\indmap{E}} \\
              K_{0} (\quot{E}) \ar[ur]^{\expmap{E}}        &                   &  K_{0} (E) \ar[ll] &  &  K_{0} (\ideal{E}) \ar[ll] & }}
\caption{K-Theory of the Crossed Product}\label{uniq F : cross prod}
\end{figure}
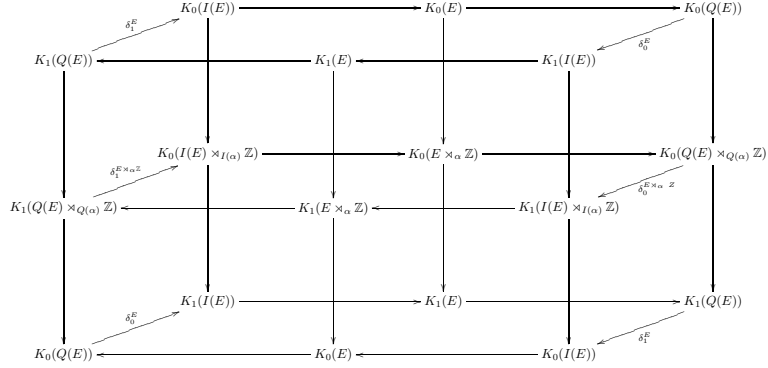

\begin{lemma}\label{uniq L : cdcross}
Let $E$ be an \ealg  and let $\alpha \in \Aut(E)$.  Then Figure \ref{uniq F : cross prod} is a commutative diagram.
\end{lemma}

\begin{proof}
This is clear from the naturality of the Pimsner-Voiculescu exact sequence.
\end{proof}

\begin{definition}\label{uniq D : decoy extension}
Let $E$ be a unital essential extension of $A$ by $I \in \PN$.  Note that we may assume $E$ is a unital C*-subalgebra of $\multialg{I}$.  Suppose $U \in \multialg{I}$ is a unitary such that $U x - x U \in I$ for all $x \in E$.  Define $E_{\alpha}$ to be the C*-subalgebra of $\multialg{I}$ generated by $E$ and $U$.  Define $A_{\alpha}$ to be the C*-subalgebra of $\corona{I}$ generated by $A$ and $\pi (U)$.  Note that $\alpha  = \innerauto (U) \in \Aut(E)$ such that $\pi \circ \alpha = \pi$.   
\end{definition}

\begin{lemma}\label{uniq L : cross prod to decoy}
Let $E$ be an \ealg such that $\quot{E} = C(S^{1})$ and let $U$ be a unitary in $\multialg{ I(E) }$ such that $ U x - x U \in \ideal{E}$ for all $x \in E$.  Set $\alpha = \innerauto (U)$.   Let $\ftn{ h } { E \rtimes_{\alpha} \Z}{E_{\alpha}}$ be the canonical surjective map.  Then $h$ gives the following commutative diagram:
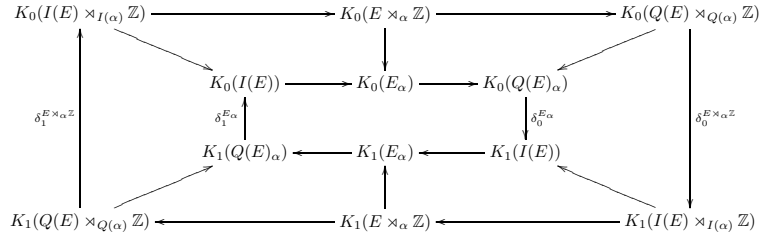
\begin{equation*}
\scalebox{.65}{
\xymatrix{ K_{0} ( \ideal{E} \rtimes _{ \ideal{\alpha} } \set{Z} ) \ar[rr] \ar[dr]    &     &  K_{0} (E \rtimes _{ \alpha } \set{ Z } ) \ar[rr] \ar[d] &  & K_{0} (\quot{E}  \rtimes _{ \quot{\alpha} } \set{Z} ) \ar[ddd]^{\expmap{E\rtimes_{\alpha} \Z}}  \ar[ld]\\
        &  K_{0} (\ideal{E})  \ar[r] & K_{0} (E_{ \alpha } ) \ar[r]  &K_{0} ( \quot{E}_{\alpha} ) \ar[d]^{\expmap{E_{\alpha}}} & \\
         & K_{1} ( \quot{E}_{\alpha} ) \ar[u]^{\indmap{E_{\alpha}}} & K_{1} ( E_{\alpha} ) \ar[l] & K_{1} (\ideal{E}) \ar[l] & \\
 K_{1} (  \quot{E} \rtimes _{ \quot{\alpha} } \set{Z} ) \ar[uuu]^{\indmap{E\rtimes_{\alpha} \Z}} \ar[ur] &   & K_{1} ( E \rtimes _{ \alpha } \set{ Z } ) \ar[ll] \ar[u] &   & K_{1} ( \ideal{E} \rtimes _{ \ideal{\alpha} } \set{Z} ) \ar[ll] \ar[lu]
} }
\end{equation*}
\end{lemma}

\begin{proof}
Note that $h$ sends $\ideal{E} \rtimes_{ \ideal{\alpha}} \Z$ to $\ideal{E}$.  Hence, we have the following commutative diagram:
\begin{equation*}
\scalebox{.9}{\xymatrix{
0 \ar[r]	& \ideal{E} \rtimes_{ \ideal{\alpha}} \Z \ar[r] \ar[d]^{\ideal{h}} 	& E \rtimes_{\alpha } \Z \ar[r] \ar[d]^{h} 	&\quot{E} \rtimes_{\quot{\alpha}} \Z \ar[r] \ar[d]^{\quot{h}} 	&0\\
0 \ar[r]	& \ideal{E} \ar[r]		& E_{\alpha} \ar[r]		& \quot{E}_{\alpha} \ar[r] 	& 0
}}
\end{equation*} 
The lemma now follows from the naturality of the six-term exact sequence in \kthy.
\end{proof}

\begin{lemma}\label{uniq L : decoyimtor}
Let $E$ be an \subealg and let $U \in \multialg{\ideal{E}}$ be a unitary such that $U x  - x U \in \ideal{E}$ for all $x \in E$.  Set $\alpha = \innerauto(U)$.  Suppose $\quot{E}_{\alpha} \cong \quot{E} \otimes C(S^{1})$ and $\alpha_{*,i} = (\id_{E})_{*,i}$ on $K_{i}(E)$ for $i= 0 ,1$.  Then $\expmap{E_{\alpha}} = 0$.  
\end{lemma}

\begin{proof}
Note that the homomorphism from $K_{0}(\quot{E} \rtimes_{\alpha} \Z)$ to $K_{0}(\quot{E}_{\alpha})$ in Lemma \ref{uniq L : cross prod to decoy} is an isomorphism.  If $K_{1}(\ideal{E}) = 0$, then it is clear that $\expmap{E_{\alpha}} = 0$.  Suppose $K_{1}(\ideal{E})$ is torsion free and $\ker \indmap{E} \neq \{ 0 \}$.  Hence, $\ran \indmap{E}$ is a torsion group.  Note that $\expmap{E} = 0$.  By performing a diagram chase in Figure \ref{uniq F : cross prod}, we see that $\ran \expmap{E \rtimes_{\alpha} \Z}$ is a torsion group.  Hence, by Lemma \ref{uniq L : cross prod to decoy}, $\ran \expmap{E_{\alpha}}$ is a torsion group.  Since $K_{1}(\ideal{E})$ is torsion free, $\expmap{E_{\alpha}} = 0$.
\end{proof}

\begin{lemma}\label{uniq L : continuous path unit}
Let $E$ be an \ealg with $\quot{E} = C(S^{1})$.  Let $U \in \multialg{\ideal{E}}$ be a unitary such that $U x - x U \in \ideal{E}$ for all $x \in E$.  Set $\alpha = \innerauto(U)$.  Then there exists a norm continuous path $V_{t} \in U(\multialg{\ideal{E}})$ such that 

(1) $V_{t} x - x V_{t} \in \ideal{E}$ for all $x \in E$ and for all $t \in [0, 1]$,

(2) $V_{0} = U$, and 

(3) $\quot{E}_{ \beta_{1} } \cong \quot{E} \otimes C(S^{1})$, where $\beta_{1} = \innerauto( V_{1})$.
\end{lemma}

\begin{proof}
Let $[x,y]$ denote the element $xy - yx$.  Note that $\quot{E}_{\alpha} \cong C(X)$, where $X$ is a compact subset of $S^{1} \times S^{1}$.  Let $\{ \xi_{n} \}_{n = 1}^{\infty}$ be a subset of $X$ such that for all $k \in \Z_{ > 0}$, $\{ \xi_{n} \}_{n = k}^{\infty}$ is dense in $X$ and let $\{e_{n}\}_{n = 1}^{\infty}$ be an approximate identity of $\ideal{E}$ consisting of projections.  For all $f \in C(X)$, let $\sigma_{0} (f) = \sum_{n = 1}^{\infty} f(\xi_{n})(e_{n} - e_{n-1})$, where the sum converges in the strict topology.  By \cite[3.2]{extrr0}, there exists an abelian AF-algebra $B \subset \corona{\ideal{E}}$ such that $\ran (\tau_{0}) \subset B$, where $\tau_{0} = \pi \circ \sigma_{0}$.  Hence, there exists a self-adjoint element $h_{1} \in \corona{\ideal{E}}$ such that $\tau_{0}(\pi(U)) = \exp(i h_{1})$ and $[h_{1}, \tau_{0}(b) ] = 0$.  Let $\tau = \tau_{\alpha} \oplus \tau_{0}$.  By Proposition \ref{def P : absorbing extensions}, there exists a unitary $Z \in \multialg{\ideal{E}}$ such that $\tau_{\alpha} = \innerauto(\pi(Z)) (\tau_{\alpha} \oplus \tau_{0})$.  Let $e_{1} = \innerauto(\pi(Z)) (\tau_{\alpha}(1) \oplus 0)$ and let $e_{2} = \innerauto(\pi(Z))(0 \oplus \tau_{0}(1))$.  Note that $e_{1}$ lifts to a projection $P_{1} \in \multialg{\ideal{E}}$.  If $P_{2} = 1 - P_{1}$, then $P_{1} + P_{2} = 1$ and $\pi(P_{i}) = e_{i}$.

Let $\ftn{\tau_{2}}{\quot{E} \otimes C([0,1])}{ e_{2} \corona{\ideal{E}} e_{2} }$ be a strongly unital essential trivial extension.
Let $g$ be a self-adjoint element in $C([0,1])$ such that $\spec(g) = [0, 2\pi]$.  Then the C*-algebra $C$ which is generated by $\quot{E} \otimes 1$ and $\exp(i (1\otimes g))$ is isomorphic to $\quot{E} \otimes C(S^{1})$.  Let $h_{2} = \tau_{2}(1 \otimes g)$.  If $\tau_{2}' = \tau_{2}\vert_{C}$, then $\tau_{2}'$ is a strongly unital trivial essential extension such that $[h_{2}, \tau_{2}'(x \otimes 1)] = 0$ for all $x \in \quot{E}$.  Note that by Proposition \ref{def P : absorbing extensions},  $\innerauto(\pi(Z)) \circ \tau_{0} \vert_{\quot{E}}$ is unitarily equivalent to $\tau_{2}' \vert_{\quot{E} \otimes 1}$.  Hence, by conjugating $\tau_{2}'$ by the image of a unitary in $P_{2} \multialg{\ideal{E}} P_{2}$, we may assume $\innerauto(\pi(Z)) \circ \tau_{0} \vert_{\quot{E}} = \tau_{2}' \vert_{\quot{E} \otimes 1}$.  Let $v_{t} = \innerauto(\pi(Z))(\pi(U) \oplus \exp(i(1-t)h_{1})) (e_{1} + \exp(i t h_{2}))$.  Then $v_{0} = \pi(U)$ and $C^{*}(v_{1}, \tau_{2}' (\quot{E} \otimes 1) ) \cong \quot{E} \otimes C(S^{1})$.  Hence, there exists a norm continuous path of unitaries $V_{t} \in \multialg{\ideal{E}}$ such that $V_{t} x  -x V_{t} \in \ideal{E}$ for all $x \in E$, $V_{0} = U$, and $\pi(V_{1}) = v_{1}$. 
\end{proof}

\begin{theorem}\label{uniq T : auto approx inner}
Let $E$ be an \subealg such that $\quot{E} = C(S^{1})$.  Let $U \in \multialg{\ideal{E}}$ be a unitary such that $ U x - x U \in \ideal{E}$ for all $x \in E$.  Let $\alpha = \innerauto( U )$.  If $\alpha_{*,i} = (\id_{E})_{*,i}$ on $K_{i} (E)$ for $i = 0,1$, then $\alpha$ is approximately inner.
\end{theorem}

\begin{proof}
By Lemma \ref{uniq L : continuous path unit} and Corollary \ref{uniq C : aut same comp approx eq}, we may assume $\quot{E}_{\alpha} =  \quot{E} \otimes C(S^{1}) = C(S^{1} \times S^{1})$.  Let $u = \pi (U)$.  Note that $\spec(u) = S^{1}$.  Let $H = \frac{1}{2} ( U + U^{*} )$ and let $W  = \exp (i H)$.  Since $\spec(U) = S^{1}$, we have $\spec(H) = [ -\pi, \pi]$.  Define $\beta = \innerauto (W)$.  Since $H$ is in the C*-algebra generated by $U$, we have $ W x - x W \in \ideal{E}$ for all $x \in E$.  Let $\beta_{t} = \innerauto ( \exp ( i H (1 - t) ) )$.  Then $\beta_{t}$ is an automorphism of $E$ such that $\beta_{0} = \beta$, $\beta_{1} = \id_{E}$, and $\quot{\beta_{t} } = \id_{\quot{E}}$.  

By Lemma \ref{uniq L : continuous path unit}, there exists a norm continuous path of unitaries $V_{t} \in \multialg{\ideal{E}}$ such that $V_{0} = W$, $V_{t} x - x V_{t} \in \ideal{E} \ \mbox{for all} \ x \in E$, and $\quot{E}_{\sigma_{1} } = C(S^{1} \times S^{1})$ where $\sigma_{1} = \innerauto(V_{1})$. Let $\tau_{\alpha}$ and $\tau_{\sigma_{1}}$ be the Busby invariants associated to the extensions $E_{\alpha}$ and $E_{\sigma_{1}}$ respectively.  

Since $\ideal{E} \in \PN$, there exists an isomorphism $\lambda_{i}$ from $K_{i}(\corona{\ideal{E}})$ to $K_{1-i}(\ideal{E})$ such that $\delta_{i}^{E_{\alpha}} = \lambda_{i} \circ (\tau_{\alpha})_{*,i}$.  By Lemma \ref{uniq L : decoyimtor}, $\expmap{E_{\alpha}} = 0$.  Hence, $(\tau_{\alpha})_{*,0} = 0$.  It is clear that $(\tau_{\sigma_{1}})_{*,0} = 0$.  

It is easy to check that the homomorphism $\ftn{ ( \tau_{\alpha} )_{*, 1} }{ K_{1} ( \quot{E}_{\alpha} ) } { K_{0} (\ideal{E}) }$ is completely determined by $(\tau_{\alpha} )_{*, 1} ( [z] ) = (\lambda_{1}^{-1} \circ \indmap{E}) ([z])$ and $(\tau_{\alpha} )_{*, 1} ( [ \pi (U) ] ) = 0$ and the homomorphism $(\tau_{\sigma_{1}})_{*,1}$ is completely determined by $(\tau_{\sigma_{1}} )_{*, 1} ( [z] ) = (\lambda_{1}^{-1} \circ \indmap{E}) ([z])$ and $(\tau_{\sigma_{1}} )_{*, 1} ( [ \pi (V_{1} ) ] ) = 0$.  Hence, $( \tau_{\alpha})_{*, 1} = ( \tau_{\sigma_{1} } )_{*, 1}$ on $K_{1} ( C(S^{1} \times S^{1}) )$.  

By Proposition \ref{def P : absorbing extensions}, there exists a unitary $Z \in \multialg{\ideal{E}}$ such that $\innerauto( \pi (Z)) \circ \tau_{\sigma_{1}}  = \tau_{\alpha}$.  So, $Z^{*} V_{1} Z - U \in \ideal{E}$.  Let $V = [ Z^{*} V_{1}  Z  ]^{*} U$.  Then $\pi (V) = 1$.  Hence, $\alpha$ is approximately inner since $\alpha = \innerauto( Z^{*} V_{1}  Z ) \circ \innerauto( V )$ and $\innerauto( Z^{*} V_{1} Z ) \in \Aut_{0} (E)$. 
\end{proof}

\begin{lemma}\label{uniq L : innerid}
Suppose $E_{1}$ and $E_{2}$ are \subealgs.  Suppose $E_{1}$ is a unital C*-subalgebra of $E_{2}$, $\quot{E_{1}} = C(S^{1})$, and $\ideal{E_{1}}$ is a nonzero hereditary C*-subalgebra of $\ideal{E_{2}}$.  Let $v \in \multialg{\ideal{E_{1}}}$ be a unitary such that $v x - x v \in \ideal{E_{1}}$ for all $x \in E_{1}$.  Set $\alpha = \innerauto(v)$.  

Suppose $\alpha_{*,1} ( n \xi ) = n \xi$ in $K_{1}(E_{2})$ for some $n \geq 1$, where $\xi$ is the generator of the copy of $\ran (\pi_{*,1})$ in $K_{1}(E_{1}) \cong K_{1}(\ideal{E_{1}}) \oplus \ran( \pi_{*,1})$.  Then $\alpha_{*,i} = (\id_{E_{1}})_{*,i}$ on $K_{i}(E_{1})$ for $i = 0 ,1$.
\end{lemma}

\begin{proof}
Let $\ftn{\iota_{i}}{\ideal{E_{i}}}{E_{i}}$ be the inclusion map and let $\ftn{\pi_{i}}{E_{i}}{\quot{E_{i}}}$ be the quotient map for $i = 0,1$.  Let $\ftn{j}{E_{1}}{E_{2}}$ be the inclusion map and let $\xi = [ w ]$ for some $w \in U(E_{1})$.  It is easy to check that $\alpha_{*,0} = (\id_{E_{1}})_{*,0}$.

Suppose $K_{1}(\ideal{E_{1}}) = 0$.  Then, $(\pi_{1})_{*,1}$ is injective.  Hence, $\alpha_{*,1} = (\id_{E_{1}})_{*,1}$ since $\quot{\alpha}_{*,1} = \id_{K_{1}(\quot{E_{1}})}$.  Suppose $K_{1}(\ideal{E_{1}})$ is torsion free and $\ker \indmap{E_{1}} \neq \{ 0 \}$.  Note that $j$ induces the following commutative diagram such that the rows are exact sequences and $\ideal{j}_{*,i}$ is an isomorphism for $i = 0,1$.  
\begin{equation*}
\scalebox{.9}{\xymatrix{ 0 \ar[r] & K_{1} ( \ideal{E_{1}} ) \ar[r] \ar[d]^{\ideal{j}_{*,1}}     & K_{1} (E_{1}) \ar[r] \ar[d]^{j_{*,1}}  & K_{1} ( \quot{E_{1}} ) \ar[d]^{\quot{j}_{*,1}} \ar[r] 	& K_{0} (\ideal{E_{1}}) \ar[d]^{\ideal{j}_{*,0}} \\
                 0 \ar[r] & K_{1} ( \ideal{E_{2}} ) \ar[r]           & K_{1} (E_{2}) \ar[r]     			& K_{1} ( \quot{E_{2}} ) \ar[r]    & K_{0}(\ideal{E_{2}}) }}                 
\end{equation*}

\noindent Since $K_{1}(\quot{E_{2}}) \cong \Z \cong K_{1}(\quot{E'})$, the map $\quot{j}_{*,1}$ is either injective or the zero map.  

Case (i). Suppose $ [ \quot{j}(z) ] = 0$ in $K_{1} ( \quot{E_{2}} )$.  Then $\quot{j}_{*,1}$ is the zero map.  Hence, $j_{*, 1} ( K_{1} ( E_{2} ) ) \subset K_{1} (\ideal{E_{2}})$.  By performing a diagram chase in the above diagram, we see that $ [ w ] = a_{1} + a_{2}$, for some $a_{1} \in K_{1} (\ideal{E_{1}})$ and $a_{2} \in \ker( j_{*,1})$.  Since $v$ is a unitary in $\multialg{\ideal{E_{1}}}$, we have $\alpha_{*, 1} (a_{1}) = a_{1}$.  Since $ \pi_{1} \circ \alpha = \pi_{1}$ and since $(\pi_{1})_{*,1} \vert_{\ker (j_{*,1})}$ is injective, by a diagram chase, $\alpha_{*, 1} (a_{2}) = a_{2}$.  Therefore, $\alpha_{*, 1} = (\id_{E_{1}})_{*,1}$ on $K_{1} (E_{1})$.  

Case (ii).  Suppose $ [ \quot{j}(z) ] \neq 0$.  By the Five Lemma, $j_{*, 1}$ is injective.  Therefore, $[ \alpha (w^{n}) ] = [w^{n} ] \neq 0$ in $K_{1} (E_{1})$.  Hence, by the exactness of the Pimsner-Voiculescu exact sequence, $[w^{n}]$ lifts to an element $x$ in $K_{0} (E_{1} \rtimes_{\alpha} \Z)$.  

Consider Figure \ref{uniq F : cross prod}, where $E$ is replaced by $E_{1}$.  Since $\pi_{1} \circ \alpha = \pi_{1}$ and since the unitary group of $\multialg{\ideal{E_{1}}}$ is connected,
\begin{equation*}
0 \to K_{0}(\quot{E_{1}}) \to K_{0}(\quot{E_{1}} \rtimes_{\quot{\alpha}} \Z) \to K_{1}(\quot{E_{1}}) \to  0 \quad \mathrm{and}
\end{equation*}
\begin{equation*}
0 \to K_{1}(\ideal{E_{1}}) \to K_{1}(\ideal{E_{1}} \rtimes_{\ideal{\alpha}} \Z) \to K_{0}(\ideal{E_{1}}) \to 0
\end{equation*}
are exact sequences. Note that $\expmap{E_{1}} = 0$, $K_{1}(\ideal{E_{1}})$ is torsion free, and $\ran \indmap{E_{1}}$ is a torsion group.  Hence, by an easy diagram chase in Figure \ref{uniq F : cross prod}, we see that $[w]$ lifts to an element in $K_{0} (E_{1} \rtimes_{\alpha} \Z)$.  Hence, by the exactness of the Pimsner-Voiculescu exact sequence, $\alpha_{*,1} = (\id_{E_{1}})_{*,1}$ on $K_{1}(E_{1})$. 
\end{proof}

\subsection{UNIQUENESS THEOREMS}

\begin{definition}\label{uniq D : weakly approx constant, approx unit eq}
Let $A = C(S^{1}) \otimes B$, where $B$ is a C*-algebra and let $\epsilon>0$.  We identify $C(S^{1}) \otimes B$ with $C( S^{1} , B)$.  A finite subset $\mc{F} \subset A$ is \textbf{weakly approximately constant to within $\epsilon$} if for any $t \in S^{1}$, there exists $U(t) \in U(A)$ such that $\norm{ U(t)^{*} f(t) U(t) - f(1) } < \epsilon$ for all $f \in \mc{F}$.
  
Suppose $C$ and $D$ are unital C*-algebras and $\varphi, \psi \in \Hom(C,D)$.  Let $\mc{G} \subset C$.  We say that $\varphi$ and $\psi$ are \textbf{approximately the same on $\mc{G}$ to within $\epsilon > 0$} if $\norm{ \varphi(f) - \psi(f) } < \epsilon$ for all $f \in \mc{G}$. 
\end{definition}

\begin{definition}\label{uniq D : stdgens}
Let $\{ e_{ij} \}_{i,j=1}^{n}$ be the standard system of matrix units in $\matalg{n}{\C} \subset \matcont{n}{S^{1}}$.  Let $z$ be the standard unitary generator of $C(S^{1}) = e_{11} \matalg{n}{C(S^{1})} e_{11}$.  Then $\{ z\} \sqcup \{ e_{ij}\}_{i,j=1}^{n}$ will be called the set of standard generators for $\matcont{n}{S^{1}}$. 
\end{definition}

\begin{theorem}\label{uniq T : quot app almost id codom not pure inf}
Let $E$ be an \subealg.  Let $\ftn{\pi}{E}{\quot{E}}$ denote the quotient map.  Let $\mc{F} = \{z_{i}\}_{i=1}^{n} \subset E$  such that $\pi(\mc{F})$ contains the set of standard generators for $\quot{E}$.  Let $\epsilon > 0$.  Then there exists $\delta > 0$ such that if $E'$ is an \subealg and $\ftn{\varphi_{i}}{E}{E'}$ is a unital monomorphism for $i = 1, 2$ with

(1) $\vtildeh{\varphi_{1}} = \vtildeh{\varphi_{2}}$ on $\vtilde{E}$ and 

(2) $\quot{ \varphi_{1}}$ and $\quot{\varphi_{2}}$ are approximately the same on $\pi(\mc{F})$ to within $\delta$, 

\noindent then $\varphi_{2}$ and $\varphi_{1}$ are approximately unitarily equivalent on $\mc{F}$ to within $\epsilon$.
\end{theorem}

\begin{proof}
Let $E = p \matalg{l}{E_{1}} p$ for some $p \in \matalg{l}{E}$ with $\pi(p) = 1_{l}$.  If $E$ is not isomorphic to $\matalg{l}{E_{1}}$, then define $\ftn{ \widetilde{\varphi}_{i} }{E \otimes M_{2} }{ E' \otimes M_{2}}$ by $ \widetilde{\varphi}_{i}  = \varphi_{i} \otimes \id$.  Since $\pi(p) = 1_{l}$, by Lemma \ref{inv L : projealg}, there exists a projection $e \in \ideal{E}$ such that $( e \oplus 1_{E} ) \matalg{2}{E} ( e \oplus 1_{E} ) \cong M_{l}(E_{1})$.  Since $[ \varphi_{1} (e) ] = [ \varphi_{2}(e) ]$ in $K_{0} (\ideal{E'})$, by \cite[1.1]{corona2}, there exists $W_{1} \in U(\unitize{\ideal{E'}})$ such that $W_{1}^{*} \varphi_{1} (e) W_{1} = \varphi_{2} (e)$.  Hence, we may assume $\varphi_{1} (e) = \varphi_{2}(e)$.  Let $E_{2} = ( e \oplus 1_{E} ) \matalg{2}{E} ( e \oplus 1_{E} ) $ and let $C = ( \varphi_{1}(e) \oplus 1_{E'} ) \matalg{2}{E'} ( \varphi_{1}(e) \oplus 1_{E'} )$.  Set $\psi_{i} = \widetilde{\varphi_{i}}\vert_{E_{2}}$.  Then $\psi_{i}$ is a unital monomorphism from $E_{2}$ to $C$.  Note that $\psi_{i} (\ideal{E_{2}}) \subset \ideal{C}$ since $\varphi_{i} (\ideal{E}) \subset \ideal{E'}$.  Hence, $\psi_{i}$ induces a homomorphism $\ftn{ \quot{\psi_{i}} } { \quot{E_{2}} } {\quot{C}}$.  Clearly, $\quot{\varphi_{i}} = \quot{ \psi_{i} }$.  Therefore, we can reduce the general case to the case $E = \matalg{l}{E_{1}}$ and hence to the case $l = 1$.

Suppose $E$ is a non-trivial extension.  By Lemma \ref{inv L : stdshift}, there exists a non-unitary isometry $S_{1} \in E$ such that $E$ is generated by $ (1 - S_{1} S_{1}^{*}) \ideal{E} (1 - S_{1} S_{1}^{*})$ and $S_{1}$.  Let $q = 1 - S_{1} S_{1}^{*}$.  Then, we may assume $z_{1} = S_{1}$ and $ z_{i} = g_{i} \in q \ideal{E} q$ for $i = 2, 3,\dots n$. 

Recall that $K_{1}(E) \cong K_{1}(\ideal{E}) \oplus \ran(\pi_{*,1})$.  By Lemma \ref{def L : unitgen}, there exists $\xi \in U(E)$ such that $[\xi]$ generates the copy of $\ran(\pi_{*,1})$ in $K_{1}(E)$.  Since $E$ is generated by $\mc{F}$, there exists $0 < \delta' < \min\left\{ \frac{1}{2}, \frac{\epsilon}{2}\right\}$ such that if $C$ is any unital C*-algebra and $\gamma$ and $\lambda$ are unital homomorphisms from $E$ to $C$ such that $\norm{ \gamma(x) - \lambda(x) } < \delta'$ for all $x \in \mc{F}$, then $\norm{\gamma(\xi) - \lambda(\xi)} < 1$.

Choose $0 < \rho < \delta'$ such that if $C$ is a unital C*-algebra, $S_{2}$ is an isometry in $C$, and $x \in C$ with $\norm{S_{2} - x} < \rho$, then $x^{*} x$ is invertible in $C$.  Also, $\rho$ can be chosen such that $\norm{ S_{2} - x|x|^{-1} } < \frac{\delta'}{100}$ and $\norm{ S_{2} S_{2}^{*} - x| x |^{-2} x^{*}} < \frac{\delta'}{100}$.

Let $0 < \delta < \rho$.  Then, there exists $a \in \ideal{E'}$ such that $\norm{ \varphi_{1}(z_{1}) - \varphi_{2}(z_{1}) + a } < \delta$ since $\norm{ \quot{\varphi_{1}}(z) - \quot{\varphi_{2}}(z) } < \delta$.  Let $z_{1}' = ( \varphi_{2} (z_{1}) - a)\vert \varphi_{2}(z_{1}) - a \vert^{-1}$.  Then $(z_{1}')^{*} z_{1}' = 1$ and $\norm{ ( 1 - \varphi_{1}(z_{1}) \varphi_{1}(z_{1})^{*} ) - d } < \frac{\delta'}{100}$, where $d = 1 - z_{1}'  (z_{1}')^{*}$.  Therefore, there exists $W' \in U(\unitize{\ideal{E'}})$ such that $(W')^{*} (1 - \varphi_{1}(z_{1}) \varphi_{1}(z_{1})^{*} ) W' = d$ and $\norm{ W' - 1 } < \frac{\delta'}{50}$.  Hence, $\norm{ (W')^{*} \varphi_{1} (z_{1}) W' - z_{1}'} < \frac{\delta'}{20}$ and $\pi' (z_{1}') = \pi' \circ \varphi_{2}(z_{1})$.

Set $S = z_{1}'$ and $T = \varphi_{2} (z_{1})$.  Since $[ 1 - S S^{*} ] = [ 1 -T T^{*} ]$ in $K_{0} (\ideal{E'})$, by \cite[1.1]{corona2}, there exists $W \in U(\unitize{\ideal{E'}})$ such that $\innerauto(W)(1 - S S^{*}) = 1 -T T^{*}$.  Note that $1 - S S^{*} \neq 0$ and $1 - T T^{*} \neq 0$ since $\varphi_{1}$ and $\varphi_{2}$ are unital monomorphisms and since $d([\pi'(S)]) = d([\pi'(T)])$ in $V(E')$.   By replacing $\innerauto(W') \circ \varphi_{1}$ with $\innerauto(W W') \circ \varphi_{1}$ and $S$ with $\innerauto(W)(S)$, we may assume that $1 - S S^{*}  =  1 - T T^{*} = d$.  

By \cite[2.6]{heralgs} and \cite[3.2(i)]{corona1}, we may write $\ideal{E'} = d' \ideal{E' } d' \otimes \K$ with $[ d ] = [ d' ]$ in $K_{0} (\ideal{E'} )$.  By \cite[1.1]{corona2}, there exists $V \in U(\unitize{\ideal{E}})$ such that $V^{*} d V = d'$.  Suppose we have found $u \in U(E')$ such that $\norm{ u^{*} V^{*} \varphi_{1}(x) V u -V^{*} \varphi_{2}(x) V} < \epsilon$ for all $x \in \mc{F}$.  Then $\norm{ w^{*} \varphi_{1}(x) w -  \varphi_{2}(x) } < \epsilon \ \mbox{for all} \ x \in \mc{F}$, where $w = \innerauto(V^{*} u V)$.  Hence, we may assume $V = 1$ and $\ideal{E'} = d \ideal{E'} d \otimes \K$.

By Proposition \ref{inv P : isoapprox}, there exists $w_{1} \in U(\unitize{\ideal{E'}})$ such that $\norm{ w_{1}^{*} T w_{1} - S } < \frac{ \delta}{16}$.  We also have $\norm{ w_{1}^{*} ( 1 - T T^{*} ) w_{1} - ( 1 - S S^{*} ) } < \frac{\delta}{8}$.  Therefore, there exists $w_{2}' \in U(\unitize{\ideal{E'}})$ such that $\norm{ 1 - w_{2}' } < \frac{\delta}{4}$ and $(w_{2}')^{*} w_{1}^{*} ( 1 - T T^{*} ) w_{1} w_{2}' = 1 - S S^{*}$.  Hence, we may assume $w_{1}^{*} ( 1- T T^{*} ) w_{1} = 1 - S S^{*}$ and $\norm{ w_{1}^{*} T w_{1} - S} < \frac{\delta}{2}$.

Note that $\innerauto(w_{1})\circ \varphi_{2}$ and $\innerauto(W') \circ \varphi_{1}$ map $q \ideal {E} q$ to $d \ideal{E'} d$ 
and $\vtildeh{\innerauto(w_{1}) \circ \varphi_{2}} = \vtildeh { \innerauto(W') \circ \varphi_{1} \vert_{q \ideal{E} q} }$ on $\totalk { q \ideal{E} q }$.  Hence, by \cite[4.10]{stableapprox}, there exists $u \in U(d \ideal{E'} d)$ such that 
$\norm{ w_{1}^{*}\varphi_{2} ( z_{i} ) w_{1} - u^{*} (W')^{*}\varphi_{1}(z_{i})  W' u } < \frac{\delta}{2}$ for all $i = 2, \dots, n$.

Let $X  = \setof{ S^{m} b (S^{*})^{n} }{ b \in d \ideal{E'} d \ \mathrm{and}  \ m, n \in \Z_{\geq 0}}$ and let $I$ be the closed linear span of the set $X$.  Let $E_{2}$ be the C*-algebra generated by $I$ and $S$.  Note that $S$ is the standard unilateral shift of $I$ and $\unit{E_{2}} = \unit{E'}$.  Hence, $E_{2}$ is a unital essential extension of $C(S^{1})$ by $I = \ideal{E_{2}}$.  Let $w_{2}= \sum_{n=0}^{\infty} S^{n} u (S^{*})^{n}$, where the sum converges in the strict topology.  Then $w_{2}$ is a unitary in $\multialg{\ideal{E_{2}}}$ and $w_{2}^{*} S w_{2} = S$.  Hence, $\alpha = \innerauto(w_{2}) \in \Aut(E_{2})$ and
\begin{equation*}
\norm{ w_{1}^{*} \varphi_{2} ( z_{i} ) w_{1} - w_{2}^{*} (W')^{*} \varphi_{1}(z_{i}) W' w_{2} }  < \frac{\delta}{2} < \delta'
\end{equation*}
for all $i = 2, \dots, n$ and 
\begin{equation*}
\norm{ w_{1}^{*} \varphi_{2} ( z_{1} ) w_{1} -  w_{2}^{*} (W')^{*} \varphi_{1}(z_{1}) W' w_{2}} < \frac{\delta}{16} + \frac{\delta'}{20} < \delta'.   
\end{equation*}
Hence, by the choice of $\delta'$ and by Lemma \ref{uniq L : innerid}, $(\innerauto(v))_{*,i} = \id_{K_{i}(E_{2})}$.  Thus, by Theorem \ref{uniq T : auto approx inner}, there exists $w_{3} \in U(E_{2})$ such that $\norm{w_{3}^{*} x w_{3} - w_{2}^{*} x w_{2} } < \frac{\epsilon}{2}$ for all $x \in \{ S \} \cup \setof{ (\innerauto(W') \circ \varphi_{1})(z_{i}) }{ i = 2, \dots, n}$.  Therefore, 
\begin{equation*}
\norm{ w_{1}^{*} \varphi_{2}(z_{i}) w_{1} - w_{3}^{*} (W')^{*}\varphi_{1}(z_{i}) W' w_{3}}  < \delta' + \frac{\epsilon}{2} < \epsilon
\end{equation*}
for all $i = 2, \dots, n$.  Also, we have
\begin{equation*}
\norm{w_{1}^{*} \varphi_{2})(z_{1}) w_{1} - w_{3}^{*} (W')^{*} \varphi_{1}(z_{1}) W' w_{3}}< \frac{\delta}{16} + \frac{\epsilon}{2} + \frac{\delta'}{20} < \epsilon.
\end{equation*}

Now, suppose $E$ is a trivial extension.  Let $0 < \delta < \min\left\{ \frac{1}{8} , \frac{\epsilon}{5} \right\}$.  By Proposition \ref{def P : diagonalization}, $E$ is generated by $\ideal{E}$ and $U = \sum_{n=1}^{\infty} \lambda_{n} p_{n}$.  Therefore, we may assume $z_{1} = (1 - p_{N} ) U (1 - p_{N})$ and $z_{2}, \dots , z_{n} \in p_{N} \ideal{E} p_{N}$, where $p_{N}$ is a projection in $\ideal{E}$.  Since $\vstar{\ideal{\varphi_{1}}} = \vstar{\ideal{\varphi_{2}}}$ on $\vstar{\ideal{E}}$, by \cite[1.1]{corona2}, there exists $W' \in U(\unitize{\ideal{E'}})$ such that $(W')^{*} \varphi_{1}(p_{N}) W' = \varphi_{2}(p_{N})$.  Hence, we may assume $\varphi_{1}(p_{N}) = \varphi_{2}(p_{N})$.  Let $q = \varphi_{1}(1 - p_{N})$.  

Note that there exists a $b \in q \ideal{E'} q$ such that $\norm{\varphi_{1}(z_{1}) - \varphi_{2}(z_{1}) + b } < \delta$.  Since $0 < \delta < 1/2$, $\vert \varphi_{2}(z_{1}) - b \vert$ is invertible.  Let $w  = (\varphi_{2}(z_{1}) - b)  \vert \varphi_{2}(z_{1}) - b \vert^{-1}$.  Then $w$ is a unitary in $q E' q$ and $\norm{ \varphi_{1}(z_{1}) - w } < 4 \delta$.  Note that $\pi'(w) = (\pi' \circ \varphi_{2})(z_{1})$ and $[w] = [ \varphi_{1}(z_{1}) ] = [ \varphi_{2}(z_{1}) ]$ in $K_{1}( q E' q)$.  Therefore, by Lemma \ref{pert L : approx unitarily eq unitaries}, there exists a unitary $U_{1} \in \C q + q\ideal{E'}q$ such that 
$\norm{ U_{1}^{*} w U_{1} - \varphi_{2}(z_{1}) } < \delta$.
Hence, $\norm{U_{1}^{*} \varphi_{1}(z_{1}) U_{1} - \varphi_{2}(z_{1}) } < 5 \delta < \epsilon$.  Since $\vtildeh{ \varphi_{1}\vert_{p_{N} \ideal p_{N}}} = \vtildeh{\varphi_{2}\vert_{p_{N} \ideal p_{N}}}$ on $\totalk{p_{N} \ideal{E} p_{N} }$, by \cite[4.10]{stableapprox}, there exists $u_{1} \in U((1 - q) \ideal{E'} ( 1 - q))$ such that $\norm{ u_{1}^{*} \varphi_{1}(z_{i}) u_{1} - \varphi_{2}(z_{i}) } < \epsilon$ for all $i = 2, \dots, n$.  Let $W = U_{1} + u_{1}$.  Then $W \in U(E')$ such that $\norm{ W^{*} \varphi_{1}(z_{i}) W - \varphi_{2}(z_{i}) } < \epsilon$ for all $i = 1, 2, \dots , n.$

By the proof, we see that $\delta$ is independent of the C*-algebra $E'$ and of the homomorphisms $\varphi_{1}$ and $\varphi_{2}$.
\end{proof}

\begin{theorem} \label{uniq T : uniq 1}
Let $E$ be an \subealg such that $\quot{E} \cong \matcont{l}{S^{1}}$ and let $\mc{F} \subset E$ be a finite subset of $E$ such that $ \pi (\mc{F}) $ contains the standard generators for $\quot{E}$.  Let $\epsilon > 0$.  Then there exists a $\delta > 0$ such that the following holds:

Let $E'$ be an \subealg.  Suppose $\varphi_{1}, \varphi_{2} \in \Hom(E, E')$ induce the same map on $\vtilde{E}$.  

(1)  If $\varphi_{1}$ and $\varphi_{2}$ are injective and the induced maps $\quot{\varphi_{1}}$ and $\quot{\varphi_{2}}$ are approximately the same to within $\delta$ on $\pi(\mc{F})$, then $\varphi_{1}$ and $\varphi_{2}$ are approximately unitarily equivalent on $\mc{F}$ to within $\epsilon$.

(2)  If $\varphi_{1}$ and $\varphi_{2}$ are not injective but $(\varphi_{1})_{Q}$ and $(\varphi_{2})_{Q}$ are injective and if the induced maps $\quot{\varphi_{1}}$ and $\quot{\varphi_{2}}$ are approximately the same on $\pi(\mc{F})$ to within $\delta$, then $\varphi_{1}$ and $\varphi_{2}$ are approximately unitarily equivalent on $\mc{F}$ to within $\epsilon$. 
\end{theorem}

\begin{proof}
It is easy to see that we may assume that $\varphi_{1}$ and $\varphi_{2}$ are unital homomorphisms.

(1) If $E'$ is not a purely infinite simple C*-algebra, then the conclusion follows from Theorem \ref{uniq T : quot app almost id codom not pure inf}.  If $E'$ is a purely infinite simple C*-algebra, then the conclusion follows from \cite[6.7]{sepBDF}.

(2) We will use the same notation as in the proof of Theorem \ref{uniq T : quot app almost id codom not pure inf}.  Since $\varphi_{i}$ is not injective, there exists a homomorphism $\ftn{(\varphi_{i})_{Q}}{\quot{E}}{E'}$ such that 
$\varphi_{i} = (\varphi_{i})_{Q} \circ \pi$ and $\pi' \circ (\varphi_{i})_{Q} = \quot{\varphi_{i}}$.  Note that $\spec(\varphi_{1} (z_{1}) ) = \spec(\varphi_{2} (z_{1}) ) = S^{1}$ since $(\varphi_{1})_{Q}$ and $(\varphi_{2})_{Q}$ are injective.  By Lemma \ref{pert L : approx unitarily eq unitaries}, there exists $U \in U(E')$ such that $\norm{ U^{*} \varphi_{2} (z_{1}) U - \varphi_{1}(z_{1}) } < \epsilon$.
\end{proof}

\begin{theorem}\label{uniq T : uniq 2}
Let $E$ and $E'$ be \subealgs such that $\quot{E} \cong \matcont{n}{S^{1}}$, let $\epsilon > 0$, and let $\mc{F} \subset E$ be a finite subset such that $\pi ( \mc{F} )$ is weakly approximately constant to within $\frac{\epsilon}{3}$.  Suppose $\varphi_{1}, \varphi_{2} \in \Hom(E,E')$ induce the same map on $\vtilde{E}$.  Then we have the following.

(1)  Suppose $\quot{E'} \cong \matalg{m}{C(S^{1})}$.  Let $\{ f_{i j} \}_{i,j = 1}^{m}$ be the standard system of matrix units for $M_{m}(\C) \subset \matalg{m}{C(S^{1})}$.  If $\varphi_{1}$ and $\varphi_{2}$ are not injective, $(\varphi_{1})_{Q}(z)$ and $(\varphi_{2})_{Q}(z)$ have finite spectra in $f_{11} \quot{E'} f_{11} \cong C(S^{1})$, then $\varphi_{1}$ and $\varphi_{2}$ are approximately unitarily equivalent on $\mc{F}$ to within $\epsilon > 0$.

(2)  If $\quot{\varphi_{1}}$ and $\quot{\varphi_{2}}$ are zero, then $\varphi_{1}$ and $\varphi_{2}$ are approximately unitarily equivalent on $\mc{F}$ to within $\epsilon$.
\end{theorem}

\begin{proof}
(1)  follows from the proof of Proposition 2.20(1) in \cite{inftoeplitz}.

(2)  Note that we may assume $\varphi_{1}$ and $\varphi_{2}$ are homomorphisms from $E$ to $p \ideal{E'} p$ for some projection $p \in \ideal{E'}$ and $\varphi_{1}(1) = p = \varphi_{2}(1)$.  If both $\varphi_{1}$ and $\varphi_{2}$ are injective, then the conclusion follows from \cite[6.7]{sepBDF}.  If $\varphi_{1}$ and $\varphi_{2}$ are not injective, then $\varphi_{i} = (\varphi_{i})_{Q} \circ \pi$.  Since $V(\varphi_{1}) = V(\varphi_{2})$ on $V(E)$, we may assume $(\varphi_{1})_{Q}$ and $(\varphi_{2})_{Q}$ both agree on $\matalg{n}{\C} \subset \matcont{n}{S^{1}} \cong \quot{E}$.  Therefore, we may assume $n = 1$.  Suppose $(\varphi_{1})_{Q}$ and $(\varphi_{2})_{Q}$ are injective.  Then the conclusion follows from Lemma \ref{pert L : approx unitarily eq purely inf}.  Hence, we may assume $(\varphi_{1})_{Q}$ is not injective.  Therefore, $(\varphi_{1})_{Q} (z)$ and $(\varphi_{2})_{Q}(z)$ are in the same path component as $p$ in $U(p \ideal{E'} p)$.  By \cite[2]{approxfinitespec}, there are $\xi_{1}, \dots, \xi_{l}, \zeta_{1}, \dots , \zeta_{L} \in S^{1}$ and mutually orthogonal projections $p_{1}, \dots , p_{l}, q_{1}, \dots , q_{L} \in p \ideal{E'} p$ such that 
{\footnotesize
\begin{equation*}
\norm{(\varphi_{1})_{Q}(f) - \sum_{k = 1}^{l} (\varphi_{1})_{Q}( f(\xi_{k}) 1)p_{k} } < \frac{\epsilon}{3}
\quad \mathrm{and} \quad
\norm{(\varphi_{2})_{Q}(f) - \sum_{j = 1}^{L} (\varphi_{2})_{Q}( f(\zeta_{j}) 1)q_{j} } < \frac{\epsilon}{3}
\end{equation*}
}

\noindent for all $f \in \pi( \mc{F} )$.  

Set $(\varphi_{1}')_{Q}(g) = \sum_{k = 1}^{l} (\varphi_{1})_{Q}(g(\xi_{k}) 1)p_{k}$ and $(\varphi_{2}')_{Q}(g) = \sum_{j = 1}^{L} (\varphi_{2})_{Q}(g(\zeta_{j}) 1)q_{j}$ for all $g \in \quot{E}$.  Since $\pi (\mc{F})$ is weakly approximately constant to within $\frac{\epsilon}{3}$, by part (1), $(\varphi_{1}')_{Q}$ and $(\varphi_{2}')_{Q}$ are approximately unitarily equivalent on $\pi (\mc{F})$ to within $\frac{\epsilon}{3}$.  Therefore, $\varphi_{1}$ and $\varphi_{2}$ are approximately unitarily equivalent on $\mc{F}$ to within $\epsilon$.
\end{proof}

\section{THE EXISTENCE THEOREM}\label{existence}
\begin{theorem}\label{ex T : existence thm V}
Let $E$ and $E'$ be two finite direct sums of \ealgs.  Let $\ftn{\alpha}{\vstar{E}}{\vstar{E'}}$ be a homomorphism such that $\alpha_{v} ( [ \unit{E} ] ) = [ P ] $ for some projection $P \in E'$.  Let $\{ \alpha_{i} \}_{i = 1}^{6}$ be the unique map from $\sixk (E)$ to $\sixk(E')$ induced by $\alpha$. (See Proposition \ref{inv P : hom inv}.)  Suppose $\psi \in \Hom(\quot{E},  \quot{E'})$ induces $\alpha_{3}$ and $\alpha_{6}$.  Then there exists $\varphi \in \Hom(E,E')$ such that $\varphi$ induces $\alpha$, $\pi ' \circ \varphi = \psi \circ \pi$, and $\varphi ( \unit{E} ) = P$.
\end{theorem}

\begin{proof}
By \cite[2.6]{heralgs} and \cite[1.2]{corona1}, $\ideal{E} = q I q \otimes \K$, for some $I \in \PN$ and for some projection $q \in I$.  Let $\{ e_{ij} \}_{i, j = 1}^{\infty}$ be the standard system of matrix units for $\K \subset q I q \otimes \K$.    
 
It is clear that we may assume $E'$ has only one summand.  Write $E = \bigoplus_{j = 1}^{k} E_{j}$, where each $E_{j}$ is an \ealg.  Denote the unit of $E_{j}$ by $1^{(j)}$.  Write $\psi = \bigoplus_{j = 1}^{k} \psi_{j}$, where each $\psi_{j}$ is in $\Hom(\quot{E_{j}}, \quot{E'})$.  Let $e_{j} = \psi (\unit{\quot{E_{j}}})$.  By \cite[2.5 and 2.9]{quasidiag} and by Proposition \ref{def P : liftingproj}, there exists a collection of mutually orthogonal projections $\{ d_{j} \}_{j = 1}^{k} \subset E'$ such that $\pi' (d_{j}) = e_{j}$ for all $j  = 1, 2, \dots, k$.  We may assume that $ P = d_{1} + \cdots + d_{k}$ and $\alpha_{v} ( [ 1^{(j) } ] ) = [d_{j}]$.  Thus, it is enough to show that for each $j$, there exists $\varphi_{j} \in \Hom(E_{j}, d_{j} E' d_{j})$ that induces $\alpha \vert_{\vstar{E_{j}}}$ such that $\pi ' \circ \varphi = \psi_{j} \circ \pi$ and $\varphi_{j} (\unit{E_{j}} ) = d_{j}$.  So, we may assume that $E$ has only one summand.  Moreover, since $d_{j} E' d_{j}$ is again an \ealg (Proposition \ref{def P : her ealg}), we may assume $P = \unit{E'}$.

Case (1). $E'$ is a unital amenable purely infinite simple C*-algebra in \bootstrap.  This case is an easy consequence of \cite[6.7]{sepBDF}. 

Case (2).  Suppose $E'$ is not a unital purely infinite simple C*-algebra.  Note that $E$ can not be a unital purely infinite simple C*-algebra.  Let $E = q \matalg{l}{E_{1}} q$ and let $E ' = P \matalg{k}{E_{2}} P$, where $E_{i}$ is an \ealg with $\quot{E_{i}} = C(S^{1})$.  We will show that we can reduce this general case to the case $l = 1$ and $q = 1$.  

First assume that $P \neq 1_{k}$.  Note that $e = 1_{l} - q$ is a projection in $\ideal{E}$ such that $e q  = 0$ and $e + q = 1_{l}$.  Since $\alpha_{v}(V(\ideal{E})) \subset V(\ideal{E'})$, there exists a projection $e' \in (1_{k} - P) \matalg{k}{\ideal{E_{2}}}(1_{k} - P)$ such that $[ e' ] = \alpha( [e ])$.  Let $E '' = ( P + e') \matalg{k}{E_{2}} (P + e')$.  Let $\pi ''$ denote the quotient map from  $E''$ onto $\quot{ E''}$ and let $\pi'$ denote the quotient map from $E'$ onto $\quot{E'}$.  

Suppose we have found $\varphi' \in \Hom(\matalg{l}{E_{1}}, E'')$ that induces $\alpha$ such that $\varphi' ( 1_{l}) = P + e'$ and $\pi' \circ \varphi' = \psi \circ \pi$.  Let $P' = \varphi' (q)$.  Then, by \cite[1.1]{corona2}, there exists $U \in U(E'')$ such that $U^{*} P' U = P$.  Denote the inclusion map from $E'$ to $E''$ by $j$.  It is easy to see that $\quot{j}$ is an isomorphism.  Choose $u$ such that $\quot{j}(u) = \pi''(U)$.  Since $\indmap{E''} ( [ \pi'' (U) ] ) = 0$ and since $\quot{j}_{*,1}$ is an isomorphism, we have $\indmap{E'}([ u ]) = 0$.  By Lemma \ref{def L : trivial ext indmap zero}, there exists $W \in U(E')$ such that $\pi ' ( W ) = u$.  Let $W_{0} = e' + W^{*}$ and let $V = U W_{0}$.  Take $\varphi = \innerauto (V) \circ \varphi' \vert_{E}$.  Then $\varphi$ induces $\alpha$, $\pi' \circ \varphi= \psi \circ \pi$, and $\varphi (\unit{E} ) = P$.  Hence, we may assume $E = \matalg{l}{E_{1}}$.

Assume that $P = 1_{k}$.  Suppose we have found $\varphi' \in \Hom(\matalg{l}{E_{1}}, \matalg{k}{E_{2}})$ such that $\varphi'$ induces $\alpha$, $\pi' \circ \varphi' = \psi \circ \pi$, and $\varphi'(1_{l}) = 1_{k}$.  Note that $\varphi'(q) \sim 1_{k} $.  Let $q' = \varphi' (q)$.  Then $1_{k} - q' \in \ideal{E'}$ and $[ 1_{k} - q' ] = 0$ in $K_{0}(\ideal{E'})$.  Let $\varphi = \iota \circ \varphi' \vert_{E}$, where $\ftn{\iota}{ q' E' q'}{ E' }$ is the homomorphism that sends $q' x q'$ to $q' x q' + 1 - q'$.  Since 
$q' \sim 1_{k}$ and $1 - q' \in \ideal{E'}$, it is easy to check that $\varphi$ is the desired homomorphism.

We will now show that we may assume $E = E_{1}$.  Let $\{ p_{i j} \}_{i, j = 1}^{l}$ be a system of matrix units for $\matalg{l}{\C} \subset \matalg{l}{E_{1}} = E$.  Let $\md{q}_{i j} = \psi ( \pi (p_{i j}) )$.  By Corollary \ref{def C : liftmatunit}, there exists a system of matrix units $\{ q'_{ij} \}_{i,j = 1}^{l} \subset E'$ such that $\pi ' ( q'_{ij}) = \md{q}_{ij}$ for all $i , j = 1, 2, \dots l$.  Note that $\alpha_{v} ( [ p_{ii} ] ) = \alpha_{v} ( [ p_{11} ] ),\ [ q'_{ii} ] = [ q_{11}' ]$, and $\alpha_{6} ( [ \pi ( p_{ii}) ] ) = [ \md{q}_{11}' ]$.  Since $\alpha_{v} ([ 1_{l} ]) = [ 1_{E'} ]$, if $d = 1_{E'} - \sum_{i = 1}^{l} q_{ii}' \neq 0$, then there are mutually orthogonal and mutually equivalent projections $a_{1}, a_{2}, \dots, a_{l} \in d E' d$ such that $\sum_{i = 1}^{l} a_{i} = d$.  So, there exists a system of matrix units $\{a_{ij} \}_{i , j = 1}^{l} \subset d E' d$ such that $a_{ii} = a_{i}$ for $i = 1, 2, \dots , l$.  Let $q_{ij} = q_{ij}' + a_{ij}$ for $i , j = 1, 2, \dots l$.  Then $\pi '( q_{ij} ) = \md{q}_{ij}$ and $\alpha ( [p_{ii}] ) = [ q_{ii} ]$.  It is now clear that it is enough to show that there exists a unital homomorphism $\ftn{\varphi}{p_{11} E p_{11}}{ q_{11} E' q_{11}}$ which induces $\alpha$ and $\pi '\circ \varphi = \psi \circ \pi$.  So, for the rest of the proof we will assume $E$ is an \ealg with $\quot{E} = C(S^{1})$ and $\alpha ( [ 1_{E} ]) = [ 1_{E'} ]$.      

Let $W$ be an isometry in $E'$ such that $\pi ' (W) = \psi (z)$.  Let $\underline{0}$ be the zero element in $K_{0} (\ideal{E}) \subset V(\ideal{E})$ which is represented by a nonzero projection.  We break case (2) into three sub-cases.  

(i)  Suppose $\alpha_{v} ( V(\ideal{E}) ) = \{ 0 \}$.  Then, there exists $\ftn{\alpha'}{\vstar{\quot{E}}}{\vstar{E'}}$ such that $\alpha' \circ \vstar{\pi} = \alpha$.  Note that $\hypon{\quot{E}} \cong K_{1} (\quot{E}) \cong \Z$ is generated by $[ z ]$.  So, there exists $W \in U(E')$ such that $ [ W ] = \alpha'( [ z ] )$ and $\pi'(W) = \psi(z)$.  Let $B$ be the C*-subalgebra generated by $W$.  Then there is a unital homomorphism $\ftn{ \varphi' } { \quot{E} } { B }$ such that $\varphi' (z) = W$.  Hence, $\varphi = \varphi' \circ \pi$ is the desired homomorphism.

(ii)  Assume that $E$ is a non-trivial extension and $\alpha_{v} ( V(\ideal{E}) ) \neq \{ 0 \}$.  Hence, $\alpha_{v} \circ d( [ S_{1}] ) \neq 0$.  By Lemma \ref{def L : liftunitiso}, we may assume $W$ is a non-unitary isometry.  We claim that there exists a unital homomorphism $\ftn{ \varphi' }{ E } { E '}$ such that $\ideal{\varphi'}$ induces $\alpha \vert_{\vstar{\ideal{E}}}$ and $\pi' \circ \varphi' = \psi \circ \pi $.  Let $f$ be a nonzero projection in $\ideal{E'}$ such that $\alpha_{v} ( [ e_{11} ] ) = [ f ]$.  By \cite[2.6]{heralgs} and \cite[1.2]{corona1}, we may assume $\ideal{E'} = f \ideal{E'} f \otimes \K$.  Hence $\ideal{E'}$ has an approximate identity consisting of projections, $\{ p_{i} \}_{i = 1}^{\infty}$, such that $ f_{i} = p_{i} - p_{i - 1}$ and $[ f_{i} ] = [ f ]$ for all $i \in \Z_{>1}$.  Therefore, by \cite[6.7]{sepBDF}, there exists a monomorphism $\ftn{\varphi_{1} } { \ideal{E} } { \ideal{E'}}$ such that $\varphi_{1} (e_{ii}) = f_{i}$ for all $i \in \Z_{>0}$ and $\varphi_{1}$ induces $\alpha \vert_{\vstar{\ideal{E}}}$.  By \cite[3.2]{inftoeplitz}, there exists a unital monomorphism $\ftn{\widetilde{\varphi_{1}}}{ \multialg{\ideal{E}}}{\multialg{\ideal{E'}}}$ extending $\varphi_{1}$.  Set 
$\varphi_{2} = \widetilde{\varphi_{1}} \vert_{E}$.  Then $C^{*}(W, \ideal{E'})$ and $C^{*}(\varphi_{2}(S_{1}), \ideal{E'})$ are equivalent unital essential extensions of $C(S^{1})$ by $\ideal{E'}$.  Then, by Proposition \ref{def P : absorbing extensions}, there exists a unitary $U \in \multialg{\ideal{E'}}$ such that $ \pi'(U) ( \pi' \circ \varphi_{2} (S_{1}) ) \pi'(U)^{*} = \pi' (W)$.  Then $\varphi' = \innerauto(U) \circ \varphi_{2}$ is the desired homomorphism.  This proves the claim.

We will now show that there exists a unital $\ftn{\varphi}{E}{E'}$ which induces $\alpha$ and $\pi' \circ \varphi = \psi \circ \pi$.
Suppose $\ker \indmap{E} = \{ 0 \}$.  Then $K_{1}(E) \cong K_{1}(\ideal{E})$, where the isomorphism is induced by the inclusion map from $\ideal{E}$ to $E$.  Hence, the homomorphism constructed above is the desired homomorphism.  Suppose $\ker \indmap{E} \neq \{ 0 \}$.  Let $e$ be a nonzero projection in $\ideal{E'}$ such that $ [e]$ is the zero element in $K_{0} ( \ideal{E'})$.  Note that the inclusion $\ftn{j}{ (\unit{E'} - e ) E' (\unit{E'} - e) } { E'}$ induces an isomorphism from $\vstar{(\unit{E'} - e ) E' (\unit{E'} - e)}$ onto $\vstar{E'}$.  Let $\pi'$ denote the quotient map from $E'$ onto $\quot{E'}$ and from $(\unit{E'} - e ) E' (\unit{E'} - e)$ onto $\quot{(\unit{E'} - e ) E' (\unit{E'} - e)}$.  Since $\quot{(\unit{E'} - e ) E' (\unit{E'} - e)} = \quot{E'}$, we have $\quot{j} = \id_{\quot{E'}}$.  Therefore, we may choose a non-unitary isometry $W' \in (\unit{E'} - e) E' (\unit{E'} - e)$ such that $\pi' ( W' ) = \psi( z)$.  By the above claim, there exists a unital monomorphism $\ftn{\varphi'}{E}{( \unit{E'} - e ) E' ( \unit{E'} - e) }$ such that $\ideal{\varphi'}$ induces $\alpha \vert_{\vstar{\ideal{E}}}$ and $\pi' \circ \varphi' = \psi \circ \pi$.  

By Lemma \ref{inv L : embdhyp}, $K_{1}(E)$ can be identified as a subsemigroup of $\hypon{E}$.  Since $\ker \indmap{E} \neq \{ 0 \}$, we have $K_{1}(E) \cong K_{1}(\ideal{E}) \oplus \Z$.  By Lemma \ref{def L : unitgen}, there exists $w \in U(E)$ such that $[w] = (0 , 1)$.  By Theorem \ref{inv T : computation invariant k}, $\hypon{e \ideal{E'} e} \cong K_{1}(e \ideal{E'} e) \sqcup \{ \hzero \}$.  It is easy to see there exists a homomorphism $\ftn{ \beta} {\vstar{E}}{\vstar{ e \ideal{E'} e} }$ such that $ \beta (x) = [ e ]$ for all $x \neq 0$ in $V(E)$, $\beta( y )  = \hzero$ for all $y \in \hypon{E} \setminus ( K_{1} (E) \cup K_{1} (\ideal{E} ))$, and $\beta ( [ w ]) = \alpha_{k}( [ w ] ) - (\varphi')_{*, 1} ( [ w ])$.  Therefore, by case (1), there exists a unital homomorphism $\ftn{ \varphi'' } { E } {e E' e}$ such that $\varphi''$ induces $\beta$.  Then $\varphi = \varphi' + \varphi''$ is the desired homomorphism.

(iii)  Suppose $E$ is a trivial extension and $\alpha (V( \ideal{E})) \neq \{ 0 \}$.  Since $E$ is a trivial extension, $\hypon{E} \cong \hypon{\ideal{E}} \oplus \Z$, where the copy of $\Z$ is generated by a unitary in $E$.  By Lemma \ref{inv L : homvo}, $\alpha(0,1) \in K_{1}(E')$.  So, by Lemma \ref{def L : unitgen}, there exists $W \in U(E')$ such that $\alpha ( 0 , 1 ) = [ W ]$ in $K_{1}(\ideal{E'}) \oplus \Z \cong K_{1} (E') \subset \hypon{E'}$.   

We first assume that $\spec ( \pi' (W) )  = S^{1}$.  As in the proof of case (2ii), there exist an approximate identity $ \{ q_{i} \}_{i = 1}^{\infty}$ of $\ideal{E'}$ consisting of projections such that $ [ q_{i} ] - [ q_{i-1} ] = \alpha ( [ e_{ii} ] )$ for all $i \in \Z_{>0}$ and a unital monomorphism $\ftn{\varphi' } {\multialg{\ideal{E}}}{\multialg{\ideal{E'}}}$ such that $\ideal{\varphi'}$ induces $\alpha \vert_{\vstar{\ideal{E}}}$ and $\varphi' (e_{ii} ) = d_{i} = q_{i} - q_{i - 1}$ for all $i \in \Z_{>0}$.  By Proposition \ref{def P : diagonalization}, we may assume $E$ is generated by $w$ and $\ideal{E}$, where $w = \sum_{k = 1 }^{\infty} \lambda_{k} e_{kk}$ (convergence is in the strict topology).  Let $V = \varphi' (w)$.  Note that $\spec ( \pi' (V) ) = S^{1}$.  Set $E_{1} = C^{*} (V, \ideal{E'} )$ and set $E_{2} = C^{*} (W , \ideal{E'})$.  Then $E_{1}$ and $E_{2}$ are trivial essential extensions of $C(S^{1})$ by $\ideal{E'}$.  Hence, by \cite[8.3.1, pp.125]{classentropy} (\cite[7(iii)]{kirchpureinf}), there exists a unitary $U \in \multialg{\ideal{E'}}$ such that $\norm{ U^{*} V U - W } < 1$ and $U^{*} V U - W \in \ideal{E'}$.  Let $\varphi = \left(\innerauto{U} \circ \varphi' \right)\vert_{E}$.  Then $\varphi$ is the desired homomorphism.

Suppose $\spec( \pi' (W) ) = X \neq S^{1}$.  Let  $J = \setof{ a \in E }{ \pi(a)\vert_{X} = 0}$.  Since $E$ is a trivial extension, $\hypon{E} \cong \setof{ (x ,y ) } { x \in \hypon{ \ideal{E}}, y \in \Z }$ and $V(E) \cong K_{0} ( \ideal{E} ) \oplus \Z_{\geq 0}$.  Let $f$ be a nonzero projection in $\ideal{E'}$ such that $\alpha( [ e_{11} ] ) = [ f ]$.  Let $\ftn{ \eta_{v} } { V(E) } { V( \ideal{E'} )}$ be $\eta_{v} (a, b) = \alpha (a, 0)$ and let $\ftn{\eta_{k} } { \hypon{E}}{\hypon{\ideal{E'}}}$ be $\eta_{k} (x,y) = \alpha(x, 0)$.  Note that $\eta_{v} \vert_{V(\ideal{E})} = \alpha \vert_{V(\ideal{E})}$ and $\eta_{k} \vert_{\hypon{\ideal{E}}} = \alpha \vert_{\hypon{\ideal{E}}}$.  It is easy to check that $\eta = ( \eta_{v} , \eta_{k})$ is a homomorphism from $\vstar{E}$ to $\vstar{ f \ideal{E' } f }$.  Therefore, by case (1), there exists a unital homomorphism $\ftn{\gamma}{E}{f \ideal{E'} f }$ such that $\gamma$ induces $\eta$.  

Let $B$ be the hereditary C*-subalgebra of $f \ideal{E'} f$ generated by $\gamma ( J )$.  By \cite[2.6]{heralgs} and \cite[1.2]{corona1}, $B$ is a stable purely infinite simple C*-algebra in \bootstrap.  Hence, there exist approximate identities $\{ e_{n} \}_{n = 1}^{\infty}$ and $\{ e_{n}' \}_{n = 1}^{\infty}$ for  $B$ and for $\ideal{E'}$ respectively such that $ e_{n}$ and $e_{n}'$ are projections with the following property:  $e_{0} = e_{0}' = 0$, $\{ e_{n} - e_{n-1} \}_{n = 1}^{\infty}$ and $\{ e_{n}' - e_{n-1}' \}_{n = 1}^{\infty}$ are collections of mutually orthogonal projections, and for all $n \in \Z_{\geq 0}$,  $[ e_{n} - e_{n - 1} ] = [ e_{n}' - e_{n-1}' ]$.  So, there exists $v_{n} \in \ideal{E'}$ such that $v_{n} ^{*} v_{n} = e_{n}' - e_{n - 1}'$ and $v_{n} v_{n}^{*} = e_{n} - e_{n - 1}$ for all $n \in \Z_{>0}$.  Set $V = \sum_{n = 1}^{\infty} v_{n}$, where the sum converges in the strict topology.  Then $V$ is an element in the double commutant of $\ideal{E'}$.  Note that $b V \in \ideal{E'}$ for all $ b \in B$ and $V^{*} a \in B$ for all $a \in \ideal{E'}$.  Therefore, the map $\ftn{\sigma}{ B } { \ideal{E'}}$ defined by $\sigma(b) = V^{*} b V$ is a homomorphism and $[ \sigma ( p  ) ] = [ p ] $ in $K_{0} (\ideal{E'} )$ for all projections $p \in B$.  Also, $[ \sigma( u) ] = [ u ]$ in $K_{1} ( \ideal{E'} )$ for all unitaries $u \in \unitize{B}$.  Let $\gamma_{1} = \left(\sigma \circ \gamma \right)\vert_{ J }$.  Then $\gamma_{1}$ maps an approximate identity of $J$ to an approximate identity of $\ideal{E'}$.  Hence, by \cite[3.2]{inftoeplitz}, there exists a unital extension $\ftn{ \widetilde{\gamma_{1}} } { \multialg{J} } { \multialg{\ideal{E'}}}$ of $\gamma_{1}$.  

Note that $J$ is an essential ideal of $E$ since $\ideal{E}$ is an essential ideal of $E$.  Let $w \in U(E)$ be as in the case $\spec ( \pi'(W) ) = S^{1}$.  It is now easy to check that $ \spec ( \pi' \circ \widetilde{\gamma_{1}} (w) ) = X = \spec ( \pi' (W) )$.  Therefore, $C^{*} ( \widetilde{\gamma_{1} }(w) , \ideal{E'} )$ and $C^{*} ( W, \ideal{E'})$ are unital essential trivial extensions of $C(X)$ by $\ideal{E'}$.  So, by \cite[8.3.1, pp. 125]{classentropy} (\cite[7(iii)]{kirchpureinf}), there exists a unitary $U \in \multialg{\ideal{E'}}$ such that $\norm{ U^{*} \widetilde{\gamma_{1} } ( w) U - W } < 1$ and $U^{*} \widetilde{\gamma_{1} } ( w) U - W \in \ideal{E'}$.  Let $\varphi = \left(\innerauto{U} \circ \widetilde{\gamma_{1}} \right)\vert_{E}$.  Then $\varphi$ is the desired homomorphism.
\end{proof} 

Let $A$ and $B$ be two separable amenable C*-algebras satisfying the UCT.  Denote the subgroup of $\Extab( K_{*} (A) , K_{1 - *}(B) )$ consisting of all pure extensions by $\pext( K_{*} (A) , K_{1 - *}(B) )$.  Set $\kl (A,B) = \kk(A,B)/ \pext( K_{*} (A) , K_{1 - *}(B) )$.  If we set $\ext(K_{*}(A), K_{1-*}(B))$ to be the quotient $\Extab( K_{*} (A) , K_{1 - *}(B) ) / \pext( K_{*} (A) , K_{1 - *}(B) )$,
then 
\begin{equation*}
\scalebox{.90}{
\xymatrix{ 0 \ar[r] & \ext(K_{*}(A), K_{1-*}(B)) \ar[r] & \kl(A,B) \ar[r]^(.38){\Gamma} & \Hom (K_{*}(A), K_{*}(B)) \ar[r] & 0}}
\end{equation*}
is an exact sequence.  See R\o rdam \cite[Section 5]{defKL}. 

\begin{theorem} (\cite[1.4]{multcoeff})\label{ex T : total kthy KL}
Let $A$ be a C*-algebra in \bootstrap and let $B$ be a $\sigma$-unital C*-algebra.  Then there is a short exact sequence 
\begin{equation*}
\scalebox{.9}{
\xymatrix{ 0   \ar[r] &  \pext( K_{*} (A) , K_{1 - *}(B) ) \ar[r]^(.63){d} & \kk(A,B) \ar[r]^(.39){\Gamma} & \Hom_{\Lambda} ( \totalk{A} , \totalk{B} ) \ar[r] & 0}}
\end{equation*}
which is natural in each variable.  Therefore, $\kl (A, B) \cong \Hom_{\Lambda} (\totalk{A} , \totalk{B} )$, where the isomorphism is natural.
\end{theorem} 

We are now ready to prove the main result of this section.

\begin{theorem}\label{ex T : existence thm}
Let $E = \bigoplus_{i = 1}^{n} E_{i}$ and $E' = \bigoplus_{j = 1}^{k} E_{j}'$ be two finite direct sums of \subealgs.  Suppose $\ftn{\alpha}{ \vtilde{E}} { \vtilde{E'}}$ is a homomorphism satisfying $\alpha_{v} ( [ \unit{E} ] ) = [ \unit{E'} ]$ and $\alpha\vert_{ V(\ideal{E_{i} }) } \neq 0$ for all $i = 1, 2, \dots, n$.  Suppose $\ftn{ \psi } { \quot{E} } {\quot{E'} }$ is a homomorphism as in Theorem \ref{ex T : existence thm V}.  Then there exists a unital homomorphism $\ftn{\varphi}{ E}{E'}$ such that $\varphi$ induces $\alpha$ and $\pi' \circ \varphi = \psi \circ \pi$.   
\end{theorem}

\begin{proof}
It is clear that we may assume $E'$ has only one summand.  Note that if $K_{*}(A)$ and $K_{*}(B)$ are finitely generated, then $\kk(A,B)$ is naturality isomorphic to $\kl(A,B)$.  

Suppose $E'$ is a unital purely infinite simple C*-algebra.  Then the existence of $\varphi$ follows from Theorem \ref{ex T : total kthy KL} and \cite[6.7]{sepBDF}.  Suppose that $E'$ is not a purely infinite simple C*-algebra.  Let $\{ \alpha_{i} \}_{i = 1}^{6}$ be as in Theorem \ref{ex T : existence thm V}.  By Theorem \ref{ex T : total kthy KL}, there exist $\beta \in \kk ( E , E'), \beta_{I} \in \kk ( \ideal{E} , \ideal{E'} ) $, and $\beta_{Q} \in \kk ( \quot{E}, \quot{E'})$ such that 

(1)  $\Gamma  ( \beta_{I} ) = ( \alpha_{1}, \alpha_{4}),\ \Gamma( \beta ) = ( \alpha_{2} , \alpha_{5}),$ and $\Gamma( \beta_{Q} ) = ( \alpha_{3} , \alpha_{6} ) = (\psi_{*, 1}, \psi_{*, 0} )$ and

(2)  $\beta_{I} = \alpha \vert_{\totalk{\ideal{E}}}, \beta = \alpha \vert_{\totalk{E}}, \beta_{Q} = \alpha \vert_{\totalk{\quot{E}}}$,
using the identification in Theorem \ref{ex T : total kthy KL}.  

Let $[ i ]$ be the element in $\kk ( \ideal{E} , E )$ induced by $\ftn{ i } { \ideal{E} } {E}$.  We define $[ i' ] , [ \pi' ], [ \pi ]$ in a similar fashion.  Since $\alpha$ is a homomorphism from $\vtilde{E}$ to $\vtilde{E'}$ and since the isomorphism in Theorem \ref{ex T : total kthy KL} is natural, we have $[ i ] \times \beta = \beta_{I} \times [ i' ]$ and $\beta \times [ \pi' ] = [ \pi ] \times \beta_{Q}$, where $\times$ represent the Kasparov product.  

Let $p$ be a nonzero projection $\ideal{E'}$ such that $ [ p ] = 0$ in $K_{0} (\ideal{E'})$.  It is easy to see that the embedding $\ftn{j}{ ( 1 - p )  E' ( 1 - p ) } { E ' }$ induces an isomorphism from $\vstar{(1 - p) E' (1 - p)}$ onto $\vstar{E'}$.  By \cite[1.17]{uct}, there exists $[j]^{-1} \in \kk ( E ' , ( 1 - p ) E' ( 1 - p) )$ such that $[ j ] \times [ j ] ^{-1} = [ \id_{( 1 - p ) E' ( 1 - p )} ] \ \mathrm{and} \ [ j ]^{ - 1} \times [ j ] = [ \id_{E'} ]$.  By Theorem \ref{ex T : existence thm V}, there exists a unital homomorphism $\ftn{\varphi'}{E}{ ( 1 -p )E' ( 1 -p ) }$ such that $\varphi'$ induces $j_{*}^{-1} \circ \alpha$ and $\pi' \circ j \circ \varphi' = \psi \circ \pi$.  Note that $\beta - [\varphi'] \times [ j ] \in \Extab(K_{0} (E) , K_{1} ( E') )$ and $\beta_{I} - [ \ideal{\varphi'} ] \times [ \ideal{j} ]  \in \Extab( K_{0} (\ideal{E} ) , K_{1} ( \ideal{E'}) )$.  Also, we have $\beta_{Q} = [ \psi ] = [ \quot{\varphi'} ]$ in $\kk(\quot{E}, \quot{E'})$.  

Note that the following diagram is commutative where the rows are split exact sequences.  
\begin{equation*}
\scalebox{.7}{
\xymatrix{
0 \ar[r] & \Extab ( K_{ 0 } ( E) , K_{1} ( \ideal{E'} ) ) \ar[r] ^{ i_{*} '} \ar[d]_{ i ^{*} } 
        				&  \Extab ( K_{ 0 } ( E ) , K_{1} ( E ' ) ) \ar[r] ^(.44){ \pi_{*} '} \ar[d] _{ i^{*} }
					& \Extab( ( K_{ 0 } ( E ), \ran (\pi')_{*,1} ) \ar[r] \ar[d]_{ i ^{*} }
						& 0 \\
0 \ar[r]   & \Extab ( K_{ 0 } ( \ideal{E} ) , K_{1} ( \ideal{E'} ) ) \ar[r] ^{ i_{*} '}    
		& \Extab ( K_{ 0 } ( \ideal{E} ) , K_{1} ( E ' ) ) \ar[r]^(.44){ \pi_{*} '} 
			& \Extab ( K_{ 0 } ( \ideal{E} ) ,  \ran (\pi')_{*,1} ) \ar[r]  
				& 0   }
}
\end{equation*}

\noindent Let $b = \pi_{*}'(\beta  - [ \varphi' ] \times [j ])$.  Since $( \beta_{I} - [ \ideal{\varphi'} ] \times [ \ideal{j} ] ) \times [ i ' ] =  [ i ] \times  ( \beta - [\varphi'] \times [ j ] )$, by a diagram chase in the above diagram, $i^{*} (b) = 0$.  Note that $K_{0} ( E ) \cong K_{0} (\ideal{E}) / \ran \indmap{E}  \oplus \Z$.  By considering the long exact sequence between Hom and Ext induced by 
\begin{equation*}
0 \to \ran \indmap{E} \to K_{0} (\ideal{E}) \to K_{0} (\ideal{E}) / \ran \indmap{E} \to 0,
\end{equation*}
we see that $\ftn{ i^{*} } { \Extab( K_{ 0 } ( E ), \ran (\pi')_{*,1} ) } {\Extab ( K_{ 0 } ( \ideal{E} ) ,  \ran (\pi')_{*,1} )}$ is injective.  So, $b = 0$.  Therefore, $\beta_{E} - [ \varphi' ] \times [ j ] = y \times [ i ']$ for some $y \in \Extab ( K_{0} (E) , K_{1} ( \ideal{E'} ))$.  By \cite[6.7]{sepBDF}, there exists a unital homomorphism $\ftn{\varphi''} { E } { p \ideal{E'} p}$ such that $[ \varphi'' ] = y$.  Then $\varphi = \varphi' + \varphi''$ is a unital homomorphism such that $\varphi$ induces $\alpha$ and $\pi' \circ \varphi = \psi \circ \pi$.
\end{proof}

\section{A CLASSIFICATION RESULT}\label{classification}

\begin{theorem}\emph{(Classification Theorem)}\label{res T : classification result}
Let $E$ and $E'$ be unital \subaealgs with real rank zero.  Suppose $\ftn{\alpha}{\vtilde{E}}{\vtilde{E'}}$ is an isomorphism such that $\alpha_{v}( [ 1_{E} ] ) = [ 1_{E'} ]$.  Then there exists a unital isomorphism $\ftn{\varphi}{E}{E'}$ such that $\varphi$ induces $\alpha$.

By Proposition \ref{def P : induce hom aealg}, the converse is also true.
\end{theorem}

\begin{proof} 
Let $E = \dirlim ( E_{i} , \varphi_{i, i + 1})$ and let $E' = \dirlim ( E_{i}' , \varphi_{i, i+1}' )$.  Since $E$ and $E'$ are unital C*-algebras, we may assume all maps are unital.   Since $\Hom_{\Lambda   } ( \totalk{A} ,\totalk{B} )$ is naturally isomorphic to $KL (A, B)$ and $KL ( A ,  \dirlim B_{n} ) =  \dirlim  KL (A , B_{n} )$ whenever $K_{*} (A)$ is finitely generated (see \cite[7.13]{uct} and \cite[5.1]{schIII}), by Proposition \ref{inv P : interwine diag invariant V} and by passing to subsystems and reindexing, we get the following commutative diagram:
\begin{equation*}
\scalebox{.85}{\xymatrix{
\vtilde{E_{1}}  \ar[r] \ar[d]_{\alpha^{(1)}}& \vtilde{E_{2}}  \ar[r] \ar[d]_{\alpha^{(2)}}& \cdots  \ar[r] & \vtilde{E} \ar@<1ex>[d]^{\alpha} \\
\vtilde{E_{1}'} \ar[r]	\ar[ru]^{\beta^{(1)}}			  & \vtilde{E_{2}'}  \ar[ru]^{\beta^{(2)}}	\ar[r]				&\cdots  \ar[r] & \vtilde{E'} \ar@<1ex>[u]^{\beta}		      }}
\end{equation*}
where $\alpha_{i} ([\unit{E_{i} } ] ) = [ \unit{E_{i}'} ]$ and $\beta_{i} ( [ \unit{E_{i}'} ] ) = [ \unit{E_{i+1}} ]$.  Recall from Section \ref{invariant} that $\sixk(E_{i})$ represent the six-term exact sequence in \kthy induced by the extension $E_{i}$.  Therefore, the above diagram induces the following commutative diagram:
\begin{equation*}
\scalebox{.85}{\xymatrix{
\sixk(E_{1})  \ar[r] \ar[d]_{\alpha^{(1)}}& \sixk(E_{2})  \ar[r] \ar[d]_{\alpha^{(2)}}& \cdots  \ar[r] & \sixk(E) \ar@<1ex>[d]^{\alpha} \\
\sixk(E_{1}') \ar[r]	\ar[ru]^{\beta^{(1)}} & \sixk(E_{2}')  \ar[ru]^{\beta^{(2)}} \ar[r]	&\cdots  \ar[r] & \sixk(E') \ar@<1ex>[u]^{\beta}		      }}
\end{equation*}
Furthermore, the following diagram commutes:
\begin{equation*}
\scalebox{.85}{
\xymatrix{
K_{*}(\quot{E_{1}})  \ar[r] \ar[d] & K_{*} ( \quot{E_{2}} )  \ar[r] \ar[d] & \cdots  \ar[r] & K_{*}(\quot{E}) \ar@<1ex>[d]  \\
K_{*}(\quot{E_{1}'}) \ar[r]	\ar[ru] & K_{*}(\quot{E_{2}'})  \ar[ru] \ar[r]	&\cdots  \ar[r] & K_{*}(\quot{E'}) \ar@<1ex>[u]		      }}
\end{equation*}
where each homomorphism lifts to a homomorphism at the level of C*-algebras.  By Theorem \ref{ex T : existence thm} , there are unital homomorphisms $\ftn{\eta_{k}'}{E_{k}}{E_{k}'} $ and $\ftn{\gamma_{k}'}{ E_{k}'}{E_{k + 1}}$ such that $\eta_{k}'$ induces $\alpha^{(k)}$ and $\gamma_{k}'$ induces $\beta^{(k)}$.    

Let $\{ \epsilon_{n} \}_{n = 1}^{\infty} \subset \R_{>0}$ be a decreasing sequence such that $\sum_{n = 1}^{\infty} \epsilon_{n} < \infty$.  Let $\mc{F}_{n}$ be a finite subset of the unit ball of $E_{n}$ and let $\mc{F}_{n}'$ be a finite subset of the unit ball of $E_{n}'$ such that (1) $\varphi_{n , n + 1 }( \mc{F}_{n} ) \subset \mc{F}_{n+1}$ and $\varphi_{n,n + 1}' (\mc{F}_{n}') \subset \mc{F}_{n + 1}'$; (2) $\bigcup_{n = 1}^{\infty} \varphi_{n , n + 1 }( \mc{F}_{n} )$ is dense in the unit ball of $E$; and (3) $\bigcup_{n = 1}^{\infty} \varphi_{n , n + 1 }'( \mc{F}_{n}' )$ is dense in the unit ball of $E'$.  Let $\pi_{k}$ and $\pi_{k}'$ denote the quotient maps from $E_{k}$ onto $\quot{E_{k}}$ and from $E_{k}'$ onto $\quot{E_{k}'}$ respectively.  We may assume that for all $n \in \Z_{>0}$, $\pi_{n} (\mc{F}_{n})$ and $\pi_{n}'( \mc{F}_{n}')$ contain the standard generators for $\quot{E_{n}}$ and $\quot{E_{n}'}$ respectively.   

Let $\{ E_{1,i} \}_{i=1}^{l(1)}$ be the summands of $E_{1}$.  Let $\delta_{1,i} > 0$ be the positive number given in Theorem \ref{uniq T : uniq 1} corresponding to $\frac{\epsilon_{1}}{2}$, $E_{1,i}$, and the image of $\mc{F}_{1}$ in $E_{1,i}$.  Let $0 < \delta_{1} < \min\{ \delta_{1,i} \}$, where the minimum is taken over all summand of $E_{1}$.  By \cite[1.4.14]{classrealrankzero}, there exists $l_{1} \in \Z_{>0}$ such that $\quot{\varphi_{1,l_{1}}} ( \pi_{1} ( \mc{F}_{1} ) )$ is weakly approximately constant to within $\delta_{1}/140$.  By \cite[1.4.14]{classrealrankzero}, there exists $n_{1} > l_{1}$ such that $\quot{\varphi_{l_{1} , n_{1}}} ( \pi_{l_{1}} ( \varphi_{1, l_{1} } (\mc{F}_{1} ) \cup \mc{F}_{l_{1}} ))$ is weakly approximately constant to within $\delta_{1}/ 140$.  

Let $\{ E_{n_{1},i}' \}_{i=1}^{l(n_{1})}$ be the summands of $E_{n_{1}}$.  Let $\lambda_{1,i}$ be the positive number given in Theorem \ref{uniq T : uniq 1} corresponding to $\frac{\epsilon_{2}}{2}$, $E_{n_{1},i}'$, and the image of $\mc{F}_{n_{1}}' \cup (\eta_{n_{1}}' \circ \varphi_{1,n_{1}})(\mc{F}_{1})$ in $E_{n_{1},i}'$.  Let $0 < \lambda_{1} < \min\{ \lambda_{1,i}\}$, where the minimum is taken over all summand of $E_{n_{1}}'$.  By \cite[1.4.14]{classrealrankzero}, there exists $l_{1}' \in \Z_{>0}$ such that $\quot{\varphi_{n_{1} , l_{1}'}'} (\pi' ( \mc{F}_{n_{1}}' \cup {\eta_{n_{1}}'} \circ \quot{\varphi_{1, n_{1}}} ( \mc{F}_{1} ) ) )$ is weakly approximately constant to within $\lambda_{1} / 140$.  By \cite[1.4.14]{classrealrankzero}, there exists $k_{1} > l_{1}'$ such that $\quot{\varphi_{l_{1}', k_{1}}'} ( \pi_{l_{1}}' ( \mc{F}_{l_{1}'}') \cup \quot{\varphi_{n_{1}, l_{1}'}'}( \mc{F}_{n_{1}}' \cup \quot{\eta_{n_{1}}'} \circ \quot{\varphi_{1, n_{1}}} ( \mc{F}_{1} ) ) )$ is weakly approximate constant to within $\lambda_{1}/ 140$. 

By \cite[2.29 and 3.25]{classrealrankzero}, there exists $m_{2} ' > k_{1} + 1$ such that $\quot{\varphi_{l_{1} , m}}$ and $\quot{\varphi_{k_{1}+1, m} \circ \gamma_{k_{1}}' \circ \varphi_{n_{1} , k_{1}} ' \circ \eta_{n_{1}}' \circ \varphi_{l_{1} , n_{1}} }$ are approximately unitarily equivalent to within $\delta_{1}/2$ on $\pi_{l_{1}} ( \varphi_{1, l_{1}} (\mc{F}_{1} ) \cup \mc{F}_{l_{1}} )$ for any $m \geq m_{2}'$.  In particular, there exists $v_{1} \in U(\quot{E_{m_{2}'}})$ such that if $ \quot{\gamma_{k_{1} , m_{2} '}'} = \innerauto (v_{1}) \circ \quot{\varphi_{k_{1} + 1, m_{2}'}} \circ \quot{\gamma_{k_{1}}'}$, then we have the following diagrams 
\begin{equation*}
\scalebox{.50}{\xymatrix{
\quot{E_{1}} \ar[r] & \quot{E_{l_{1}}} \ar[r] & \quot{E_{n_{1}}} \ar[r] \ar[d] & \quot{E_{k_{1}}} \ar[r]                                     & \quot{E_{k_{1} + 1}} \ar[r] & \quot{E_{m_{2}'}} \\
                    &                         & \quot{E_{n_{1} }'} \ar[r]      & \quot{E_{k_{1}}'} \ar[rru]_{ \quot{\gamma_{k_{1} , m_{2} '}'}}  &                          &    }
\quad\quad
\xymatrix{
 \quot{E_{l_{1}}} \ar[r] & \quot{E_{n_{1}}} \ar[r] \ar[d] & \quot{E_{k_{1}}} \ar[r]                                     & \quot{E_{k_{1} + 1}} \ar[r] & \quot{E_{m_{2}'}} \\
                                           & \quot{E_{n_{1} }'} \ar[r]      & \quot{E_{k_{1}}'} \ar[rru]_{ \quot{\gamma_{k_{1} , m_{2} '}'}}  &                          &    }
  }                                         
\end{equation*}
where the first diagram is approximately commutative on $\pi_{1} (\mc{F}_{1} )$ to within $\delta_{1}$ and the second diagram is approximately commutative on $\pi_{l_{1}}(\varphi_{1,l_{1}}(\mc{F}_{1}) \cup \mc{F}_{l_{1}})$ to within $\delta_{1}$.  By Theorem \ref{ex T : existence thm}, there exists a unital homomorphism $\gamma_{k_{1}}''$ from $E_{k_{1}}'$ to $E_{m_{2}'}$ such that $\gamma_{k_{1}}''$ and $\varphi_{k_{1} + 1 , m_{2}'} \circ \gamma_{k_{1}}'$ induce the same map on $\vtilde{E_{k_{1}}'}$ and $\quot{\gamma''_{k_{1}}} =  \quot{\gamma_{k_{1} , m_{2} '}'}$.

Let $m_{1} = 1$.  We will show that there exists $m_{2} \geq m_{2}'$ such that $\varphi_{m_{2}' , m_{2}} \circ \varphi_{1 , m_{2}'}$ and $\varphi_{m_{2}' , m_{2}} \circ \gamma_{k_{1}}'' \circ \varphi_{n_{1} , k_{1}}' \circ \eta_{n_{1}}' \circ \varphi_{1 , n_{1}}$ are approximately unitarily equivalent on $\mc{F}_{1}$ within $\epsilon_{1}$.  Let $E_{m_{2}'}^{j}$ be the $j$th summand of $E_{m_{2}'}$ and let $\ftn{P_{j}}{ E_{m_{2}'}}{E_{m_{2}'}^{j}}$ be the projection onto the $j$th summand.  Suppose that $\ftn{P}{ E_{1}}{H}$ is the projection of $E_{1}$ onto one of its summands.  

(a)  Suppose $P_{j} \circ \varphi_{1 , m_{2}'}$ is injective on $H$ and $\quot{H} \neq 0$.  By the choice of $\delta_{1} > 0$ and by Theorem \ref{uniq T : uniq 1} (1), there exists a unitary $W_{j} \in (P_{j} \circ \varphi_{1 , m_{2}'}) ( \unit{E_{1}} ) E_{m_{2}'}^{j}  (P_{j} \circ \varphi_{1 , m_{2}'} )( \unit{E_{1}} )$ such that $P_{j} \circ \varphi_{1, m_{2}'}$ and $\innerauto (W_{j} ) \circ P_{j} \circ \gamma_{k_{1}}'' \circ \varphi_{n_{1} , k_{1} }' \circ \eta_{n_{1}}' \circ \varphi_{1, n_{1}}$ are approximately the same on  $P (\mc{F}_{1} )$ to within $\epsilon_{1}$.  Set $m(j) = m_{2}'$.

(b) Suppose $\quot{H} = 0$.  In this case, $P \circ \varphi_{1, l_{1} } (\mc{F}_{1} )$ is a finite subset of a corner of $e \ideal{E_{n_{1}}} e$ for some nonzero projection $e \in \ideal{E_{n_{1}}}$.   Note that $e \ideal{E_{n_{1}}} e$ will be mapped to $e' \ideal{E_{m_{2}'}^{j}} e'$ for some nonzero projection $e' \in \ideal{E_{m_{2}'}^{j}}$.  So, from the commutative diagrams involving $\vtilde{\cdot}$ and by \cite[4.10]{stableapprox}, there is a unitary $W_{j} \in (P_{j} \circ \varphi_{1 , m_{2}'} )( \unit{E_{1}} )  E_{m_{2}'}^{j} ( P_{j} \circ \varphi_{1 , m_{2}'} )( \unit{E_{1}} )$ such that $P_{j} \circ \varphi_{1, m_{2}'}$ and $\innerauto (W_{j} ) \circ P_{j} \circ \gamma_{k_{1}}'' \circ \varphi_{n_{1} , k_{1} }' \circ \eta_{n_{1}}' \circ \varphi_{1, n_{1}} $ are approximately the same on $P (\mc{F}_{1} )$ to within $\epsilon_{1}$.  Set $m(j) = m_{2}'$.

(c)  Suppose $P_{j} \circ \varphi_{1 , m_{2}'}$ is not injective on $H$ but $\quot{P_{j} \circ \varphi_{1 , m_{2}'}}$ is injective on $\quot{H}$.  By Theorem \ref{uniq T : uniq 1} (2), there exists a unitary $W_{j} \in (P_{j} \circ \varphi_{1 , m_{2}'}) ( \unit{E_{1}} )E_{m_{2}'}^{j}  (P_{j} \circ \varphi_{1 , m_{2}'} )( \unit{E_{1}})$ such that $P_{j} \circ \varphi_{1, m_{2}'}$ and $\innerauto (W_{j} ) \circ P_{j} \circ \gamma_{k_{1}}'' \circ \varphi_{n_{1} , k_{1} }' \circ \eta_{n_{1}}' \circ \varphi_{1, n_{1}}$ are approximately the same on $P (\mc{F}_{1} )$ to within $\epsilon_{1}$.  Set $m(j) = m_{2}'$.

(d)  Suppose $P_{j} \circ \varphi_{1 , m_{2}'}$ is not injective on $H$ and $\quot{P_{j} \circ \varphi_{1 , m_{2}'}} \vert_{\quot{H}} = 0$.  Then by Theorem \ref{uniq T : uniq 2} (2), there exists a unitary $W_{j} \in (P_{j} \circ \varphi_{1 , m_{2}'}) ( \unit{E_{1}} )E_{m_{2}'}^{j}  (P_{j} \circ \varphi_{1 , m_{2}'} )( \unit{E_{1}})$ such that $P_{j} \circ \varphi_{1, m_{2}'}$ and $\innerauto (W_{j} ) \circ P_{j} \circ \gamma_{k_{1}}'' \circ \varphi_{n_{1} , k_{1} }' \circ \eta_{n_{1}}' \circ \varphi_{1, n_{1}}$ are approximately the same on $P (\mc{F}_{1} )$ to within $\epsilon_{1}$.  Set $m(j) = m_{2}'$.

(e) Suppose  $P_{j} \circ \varphi_{1, m_{2}'}$ is not injective on $H$ and $\quot{P_{j} \circ \varphi_{1 , m_{2}'}} \vert_{\quot{H}} \neq 0$ (but $\quot{H} \neq 0$).  Therefore, $H$ is not a purely infinite simple C*-algebra.  Let $u$ be the canonical unitary generator of $\quot{H}$.  Since $P_{j} \circ \varphi_{1 , m_{2}'}$ is not injective on $H$,  $P_{j} \circ \varphi_{1 , m_{2}'}$ factors through $\quot{H}$.  Let $\kappa$ be the homomorphism from $\quot{H}$ to $E_{m_{2}'}^{j}$ induced by $P_{j} \circ \varphi_{1 , m_{2}'}$.  The commutativity at the level of $\vtilde{\cdot}$ shows that $P_{j} \circ \gamma_{k_{1}}'' \circ \varphi_{n_{1}, k_{1} }' \circ \eta_{n_{1}}' \circ \varphi_{1, n_{1}} \vert_{H}$ also factors through $\quot{H}$.  Denote the induced map by $\kappa'$.  Let $u_{1} = \kappa(u)$ and $u_{2} = \kappa' (u)$.  Note that 
$[ u_{1} ] = [ u_{2} ] \ \mathrm{in} \  K_{1} ( (P_{j} \circ \varphi_{1 , m_{2}'}) ( \unit{E_{1}} ) E_{m_{2}'}^{j}  (P_{j} \circ \varphi_{1 , m_{2}'}) ( \unit{E_{1}} ))$ and $u_{1}, u_{2}$ are map to the same element in $\quot{ (P_{j} \circ \varphi_{1 , m_{2}'} )( \unit{E_{1}} ) E_{m_{2}'}^{j}  (P_{j} \circ \varphi_{1, m_{2}'}) ( \unit{E_{1}} )}$.  If $\spec (u_{1}) = \spec( u_{2} ) = S^{1}$, by Lemma \ref{pert L : approx unitarily eq unitaries}, there exists a unitary 
$W_{j} \in \ideal{ (P_{j} \circ \varphi_{1 , m_{2}'}) ( \unit{E_{1}}) E_{m_{2}'}^{j}  (P_{j} \circ \varphi_{1 , m_{2}'}) ( \unit{E_{1}})}$
such that $\norm{ W_{j} u_{2} W_{j}^{*} - u_{1} } < \epsilon_{1}$.  Set $m(j) = m_{2}'$.  

Suppose $\spec(u_{1}) \neq S^{1}$ or $\spec(u_{2}) \neq S^{1}$.  Then $u_{1}$ and $u_{2}$ are connected to the identity in the unitary group of $(P_{j} \circ \varphi_{1 , m_{2}'} )( \unit{E_{1}}) E_{m_{2}'}^{j} ( P_{j} \circ \varphi_{1, m_{2}'} )( \unit{E_{1}})$.  Since $RR(E)=0$, by \cite[5]{exprank}, there exists $m(j) > m_{2}'$ such that 
$\norm{ \varphi_{m_{2}', m(j) } (u_{i} ) - v_{i} } < \frac{\epsilon_{1}}{2} \ \mathrm{for} \ i = 1, 2,$
where $v_{1}$, $v_{2}$ are unitaries in $(P_{j} \circ \varphi_{1 , m(j)} )( \unit{E_{1}} ) E_{m_{2}'}^{j}  (P_{j} \circ \varphi_{1, m(j)}) ( \unit{E_{1}})$ with finite spectrum.  By replacing $u_{i}$ by $v_{i}$ for $i = 1, 2$, we obtain two homomorphisms $\ftn{ h_{1}, h_{2} } { H } { (P_{j} \circ \varphi_{1 , m(j)} )( \unit{E_{1}} ) E_{m_{2}'}^{j}  (P_{j} \circ \varphi_{1 , m(j)} )( \unit{E_{1}} )}$ such that $h_{1}$ and $\varphi_{l_{1}, m(j)} \vert_{H}$ are approximately the same to within $\frac{\epsilon_{1}}{2}$ on $P (\mc{F}_{1} )$ and $h_{2}$ and $P_{j} \circ \gamma_{k_{1}}'' \circ \varphi_{n_{1} , k_{1}}' \circ \eta_{n_{1}}' \circ \varphi_{l_{1}, n_{1}} \vert_{H}$ are approximately the same to within $\frac{\epsilon_{1}}{2}$ on $P ( \mc{F}_{1} )$.  It follows from Theorem \ref{uniq T : uniq 2} (1) that there exists a unitary $W_{j} \in (P_{j} \circ \varphi_{1 , m(j)} )( \unit{E_{1}}) E_{m_{2}'}^{j}  (P_{j} \circ \varphi_{1 , m(j)} )( \unit{E_{1}})$ such that $P_{j} \circ \varphi_{1, m(j)}$ and $\innerauto (W_{j} ) \circ P_{j} \circ \circ \varphi_{m_{2}', m(j)} \circ \gamma_{k_{1}}'' \circ \varphi_{n_{1} , k_{1} }' \circ \eta_{n_{1}}' \circ \varphi_{1, n_{1}}$ are approximately the same on $P (\mc{F}_{1} )$ to within $\epsilon_{1}$.  Set $m(j) = m_{2}'$.

If $m_{2,H}'' = \max_{j} \{ m(j) \}$, then by the above observations, there exists a unitary $W \in (\varphi_{1, m_{2}''} \circ P)(\unit{E_{1}}) E_{m_{2}''} (\varphi_{1, m_{2}''} \circ P)(\unit{E_{1}})$ such that 
\begin{equation*}
\norm{\varphi_{1, m_{2}''} (f) - \innerauto(W) \circ \varphi_{m_{2}' , m_{2}''} \circ \gamma_{k_{1}}'' \circ \varphi_{n_{1}, k_{1}}' \circ \eta_{n_{1}}' \circ \varphi_{1 , n_{1}} (f) } < \epsilon_{1}
\end{equation*}
for all $f \in P (\mc{F}_{1})$.  Then let $m_{2} = \max\{m_{2,H}'' \}$, where the maximum is taken over all direct summands of $E_{1}$.  So $\varphi_{1, m_{2}}$ and  $\varphi_{m_{2}' , m_{2}} \circ \gamma_{k_{1}}'' \circ \varphi_{n_{1}, k_{1}}' \circ \eta_{n_{1}}' \circ \varphi_{1, n_{1}}$ are approximately unitarily equivalent on $\mc{F}_{1}$ to within $\epsilon_{1}$.  Hence, there exists $U \in U(E_{m_{2}})$ such that if $\eta_{1} = \eta_{n_{1}}' \circ \varphi_{1, n_{1}}$ and $\gamma_{1} = \innerauto(U) \circ \varphi_{m_{2}', m_{2}} \circ \gamma_{k_{1}}'' \circ \varphi_{n_{1}, k_{1}}'$, then $\varphi_{1, m_{2}}$ and $\gamma_{1} \circ \eta_{1}$ are approximately the same on $\mc{F}_{1}$ to within $\epsilon_{1}$.

Let $\mc{G}_{1}' = \mc{F}_{n_{1}}' \cup \eta_{1} (\mc{F}_{1})$, let $\mc{G}_{2} = \mc{F}_{m_{2}} \cup \varphi_{1, m_{2}} (\mc{F}_{1}) \cup \gamma_{1} ( \mc{G}_{1}')$, and let $\{ E_{m_{2},i} \}_{i=1}^{l(m_{2})}$ be the summands of $E_{m_{2}}$.  Let $ \delta_{2,i} > 0$ be the positive number that is given in Theorem \ref{uniq T : uniq 1} which corresponds to $\frac{\epsilon_{3}}{2}$, $E_{m_{2}}^{i}$, and the image of $\mc{G}_{2}$ in $E_{m_{2}}^{i}$.  Let $0 < \delta_{2} < \min\{\delta_{2,i}\}$, where the minimum is taken over all direct summands of $E_{m_{2}}$.  By \cite[1.4.14]{classrealrankzero}, there exists $l_{2} > m_{2}$ such that $\quot{\varphi_{m_{2}, l_{2}} }(\mc{G}_{2})$ is weakly approximately constant to within $\delta_{2}/140$.  By \cite[1.4.14]{classrealrankzero}, we get $n_{2}' > l_{2}$ such that $\quot{\varphi_{l_{2}, n_{2}'} }(\pi_{l_{2}} ( \mc{F}_{l_{2}} \cup \varphi_{m_{2}, l_{2}} (\mc{G}_{2})))$ is weakly approximately constant to within $\delta_{2}/140$.  

Note $\quot{\varphi_{l_{1}', k_{1}}'}(\pi' ( \varphi_{n_{1}, l_{1}'} ( \mc{G}_{1}')))$ is weakly approximately constant to within $\lambda_{1}/140$.  Hence, by \cite[2.29 and 3.25]{classrealrankzero}, there exists $n_{2}'' > n_{2}'$ such that $\quot{\varphi_{l_{1}', m} }$ and $\quot{\varphi_{n_{2}' , m} \circ \eta_{n_{2}}'  \circ \innerauto(U) \circ \varphi_{m_{2}, n_{2}'} \circ \gamma_{k_{1}}'' \circ \varphi_{l_{1}', k_{1}} }$ are approximately unitarily equivalent to within $\lambda_{1}/2$ on $\pi_{l_{1}'}' ( \varphi_{n_{1}, l_{1}'} (\mc{G}_{1}'))$ for any $m \geq n_{2}''$.  Set $\kappa_{1} = \innerauto(U) \circ \varphi_{m_{2}, n_{2}'} \circ \gamma_{k_{1}}'' \circ \varphi_{l_{1}', k_{1}}$.

By the choice of $l_{1}'$ and $k_{1}$, then using the same argument as above we get a positive integer $n_{2} \geq n_{2}''$  and a unitary $V \in E_{n_{2}}'$ such that $\varphi_{n_{1}, n_{2}}'$ and $\innerauto(V) \circ \varphi_{n_{2}', n_{2} } \circ \eta_{n_{2}'}' \circ \varphi_{m_{2}, n_{2}'} \circ \gamma_{1}$ are approximately the same on $f \in \mc{G}_{1}'$ to within $\epsilon_{2}$.

Set $\eta_{2} = \innerauto(V) \circ \varphi_{n_{2}', n_{2} } \circ \eta_{n_{2}'}' \circ \varphi_{m_{2}, n_{2}'}$ and $\kappa_{2} = \innerauto (V) \circ \varphi_{n_{2}', n_{2}}' \circ \eta_{n_{2}'}'$.  Then the following diagram:
\begin{equation*}
\scalebox{.85}{\xymatrix{
E_{1} \ar[d]_{\eta_{1}} \ar[r]			& E_{m_{2}}  \ar[d]_{\eta_{2}}   \\
E_{n_{1}}'  \ar[ru]^{\gamma_{1}} \ar[r]	& E_{n_{2}}'   				}}
\end{equation*}
is approximately commutative on $\mc{F}_{1}$ to within $\epsilon_{1}$ for the upper triangle and is approximately commutative on $\mc{F}_{1} \cup \eta_{1} ( \mc{F}_{1} )$ to within $\epsilon_{2}$ for the lower triangle.  Furthermore, the homomorphism $\eta_{2}$ can be rewritten as $\varphi_{n_{2}',n_{2}}' \circ \eta_{n_{2}'} \circ \varphi_{m_{2}, n_{2}'}$ or $\kappa_{2} \circ \varphi_{m_{2}, n_{2}'}$.

The choice of $l_{2}$ and $n_{2}'$ ensures that the above construction can continue.  Therefore, we get the following approximately intertwining diagram:
\begin{equation*}
\scalebox{.85}{\xymatrix{
E_{1} \ar[r] \ar[d]_{\eta_{1}} 			& E_{m_{2}} \ar[r] \ar[d]_{\eta_{2}} 		& \cdots \ar[r]		& E		\\
E_{n_{1}}' \ar[ru]^{\gamma_{1}} \ar[r]		& E_{n_{2}}' \ar[ru]^{\gamma_{2}} \ar[r]	&\cdots \ar[r]		& E' 	
	}}
\end{equation*}
Hence, there exists a unital isomorphism $\ftn{\varphi}{E}{E'}$ such that $\varphi$ induces $\alpha$.
\end{proof} 

\begin{corollary}\label{res C : k1ideal0}
Suppose $K_{1}(\ideal{E_{n}}) = 0 = K_{1}(\ideal{E_{n}'})$ for all $n \in \Z_{\geq 0}$ in Theorem \ref{res T : classification result}.  If $\ftn{\alpha}{\vstar{E}}{\vstar{E'}}$ is an isomorphism such that $\alpha_{v}([ \unit{E} ]) = [\unit{E'}]$, then there exists a unital isomorphism $\ftn{\varphi}{E}{E'}$ such that $\varphi$ induces $\alpha$.  
\end{corollary}

\begin{proof}
Since $K_{1}(\ideal{E}) = K_{1}(\ideal{E'}) = 0$, there exists an isomorphism $(\beta_{1}, \beta_{2})$ from \\
$\vtilde{E}$ to $\vtilde{E}$ such that $\beta_{1} = \alpha$.  Now apply Theorem \ref{res T : classification result}.
\end{proof}

\begin{remark}
Every unital \ATalg with real rank zero is an \subaealg with real rank zero.
\end{remark}

\begin{proof}
Let $Q_{n} = \bigoplus_{j = 1}^{l(n)} M_{k(j)} (C(S^{1}))$ and let $Q = \dirlim ( Q_{n}, \psi_{n, n+1})$ be an \ATalg.  Take $E'$ to be an \subealg such that $E'$ is a trivial extension and $\quot{E'} = C(S^{1})$.  Set $E_{n} = \bigoplus_{j = 1}^{l(n)} M_{k(j)} (E')$.  Let $\ftn{\pi_{n}}{E_{n}}{Q_{n}}$ be the quotient map.  Since $E_{n}$ is a trivial extension, there exists a unital monomorphism $\ftn{j_{n}}{Q_{n}}{E_{n}}$ such that $j_{n} \circ \pi_{n}$ is the identity map on $E_{n}$.  Define $\varphi_{n, n+1}$ from $E_{n}$ to $E_{n+1}$ by $\varphi_{n,n+1} = j_{n+1} \circ \psi_{n,n+1} \circ \pi_{n}$.  Then $E = \dirlim (E_{n}, \varphi_{n,n+1})$ is an \subaealg such that $E \cong Q$.  
\end{proof}

\begin{remark}
Every unital separable amenable purely infinite simple C*-algebra satisfying the UCT with torsion free $K_{1}$ is an \subaealg.  See \cite{kirchpureinf} and \cite{phillipspureinf}.
\end{remark}

%\bibliographystyle{siam}
%\bibliography{/Users/efrenruiz/analysis/references}

 \end{document}